\documentclass[12pt,leqno]{article}

\usepackage{amsmath,amssymb,amsfonts,theorem,latexsym,epsfig}
\usepackage{psfrag}
\usepackage{hyperref}

\newtheorem{thm}{Theorem}
\newtheorem{theo}{Theorem}[section]
\newtheorem{lem}[theo]{Lemma}
\newtheorem{cor}[theo]{Corollary}
\newtheorem{prop}[theo]{Proposition}
\newtheorem{defi}[theo]{Definition}
\newtheorem{claim}[theo]{Claim}

{\theorembodyfont{\rmfamily} \newtheorem{rem}[theo]{Remark}}
{\theorembodyfont{\rmfamily} }
{\theorembodyfont{\rmfamily} }
{\theorembodyfont{\rmfamily} \newtheorem{con}[theo]{Conjecture}}

\usepackage[all]{xy}

\newcommand{\Appendix}[1]{%
  \refstepcounter{section}%
  \addtocontents{toc}{\protect\setcounter{tocdepth}{1}}
  \addcontentsline{toc}{section}%
    {\bfseries\appendixname~\thesection:\ #1}%
    {\medskip\noindent \Large\bfseries\appendixname\ \thesection:\ #1}%
\sectionmark{#1}\smallskip\noindent
\renewcommand{\theequation}{{\bf 
{{\thesection}}.{\arabic{equation}}}}
}

\DeclareFontFamily{U}{rsf}{}
\DeclareFontShape{U}{rsf}{m}{n}{
  <5> <6> rsfs5 <7> <8> <9> rsfs7 <10->  rsfs10}{}
\DeclareMathAlphabet{\mathscr}{U}{rsf}{m}{n}
\newcommand{\mycal}[1]{\mathscr{#1}}

\DeclareMathAlphabet{\mathpzc}{OT1}{pzc}{m}{it}

\newcommand{\les}[9]{
\xymatrix{     
 0 \ar[r] & {#1} \ar[r]  &  {#2} \ar[r]  &  {#3} 
\ar@{->}`r/10pt[d] `[l] `^dl[dlll]  `^dr/10pt[dll]    [dll] \\
 &  {#4} \ar[r] & {#5} \ar[r] & {#6} 
\ar@{->}`r/10pt[d] `[l] `^dl[dlll]  `^dr/10pt[dll]    [dll] \\
 & {#7} \ar[r]  & {#8} \ar[r] & {#9}  
\ar@{->}`r/10pt[d] `[l] `^dl[dlll]  `^dr/10pt[dll]    [dll] \\
 & 0 \ar[r] & \cdots & }
}

\newcommand{\lestwo}[9]{
\xymatrix{     
 0 \ar[r] & {#1} \ar[r]  &  {#2} \ar[r]  &  {#3} 
\ar@{->}`r/10pt[d] `[l] `^dl[dlll]  `^dr/10pt[dll]    [dll] \\
 &  {#4} \ar[r] & {#5} \ar[r] & {#6} 
\ar@{->}`r/10pt[d] `[l] `^dl[dlll]  `^dr/10pt[dll]    [dll] \\
 & {#7} \ar[r]  & {#8} \ar[r] & {#9} }
}

\newcommand{\lesthree}[5]{
\xymatrix{     
 0 \ar[r] & {#1} \ar[r]  &  {#2} \ar[r]  &  {#3} 
\ar@{->}`r/10pt[d] `[l] `^dl[dlll]  `^dr/10pt[dll]    [dll] \\
 &  {#4} \ar[r] & {#5} & }
}

\newcommand{\lesfour}[8]{
\xymatrix{     
 0 \ar[r] & {#1} \ar[r]  &  {#2} \ar[r]  &  {#3} 
\ar@{->}`r/10pt[d] `[l] `^dl[dlll]  `^dr/10pt[dll]    [dll] \\
 &  {#4} \ar[r]^-{#8} & {#5} \ar[r] & {#6} 
\ar@{->}`r/10pt[d] `[l] `^dl[dlll]  `^dr/10pt[dll]    [dll] \\
 & {#7} \ar[r]  & \cdots  &  }
}

\def\punkt{\refstepcounter{subsubsection}
           \noindent{\bf \thesubsubsection.\ }}

\newcommand{\opi}[1]{\operatorname{{\it #1}}}
\newcommand{\op}[1]{\operatorname{#1}}
\newcommand{\lan}[1]{{}^{L}{#1}}
\newcommand{\lanab}[1]{{}^{L}_{\op{ab}}{#1}}
\newcommand{\Higgs}{\op{{\bf Higgs}}}
\newcommand{\gHiggs}{\op{\boldsymbol{\mathcal{H}iggs}}}
\newcommand{\gTors}{\op{\boldsymbol{\mathcal{T}\!\!\text{\it\bfseries ors}}}}

\newcommand{\gLoc}{\op{\boldsymbol{\mathcal{L}oc}}}

\newcommand{\sHeck}{\op{{\boldsymbol{\mathcal{H}}{\bf ecke}}}}
\newcommand{\trans}{\op{{\sf{Trans}}}}
\newcommand{\tens}{\op{{\sf{Tens}}}}
\newcommand{\sHom}{\opi{\mycal{H}\!\!om}}
\newcommand{\sIsom}{\opi{\mycal{I}\!\!som}}
\newcommand{\LHS}{\opi{LHS}}
\newcommand{\RHS}{\opi{RHS}}

\newcommand{\bpo}{\mathfrak{o}}
\newcommand{\bl}{\opi{{\sf{l}}}^{\op{base}}}
\newcommand{\myell}{\opi{{\sf{l}}}^{\op{cam}}}
\newcommand{\paiso}{\opi{{\sf{l}}}}
\newcommand{\gerbeiso}{\opi{\boldsymbol{\ell}}}
\newcommand{\bD}{\boldsymbol{D}}
\newcommand{\bw}{\boldsymbol{w}}
\newcommand{\rts}{\op{{\sf{root}}}}
\newcommand{\wts}{\op{{\sf{weight}}}}
\newcommand{\crts}{\op{{\sf{coroot}}}}
\newcommand{\cwts}{\op{{\sf{coweight}}}}
\newcommand{\chr}{\op{{\sf{char}}}}
\newcommand{\cchr}{\op{{\sf{cochar}}}}
\newcommand{\aj}{\op{{\sf{aj}}}}
\newcommand{\haj}[1]{{\aj}^{#1}}
\newcommand{\saj}{\op{{\mathfrak{aj}}}}
\newcommand{\shaj}[1]{{\saj}^{#1}}
\newcommand{\ab}{\op{{\sf{a}}}}
\newcommand{\sab}{\op{{\mathfrak{a}}}}
\newcommand{\Bun}{\op{{\bf Bun}}}
\newcommand{\sBun}{\op{\boldsymbol{\mathcal{B}un}}}
\newcommand{\Pic}{\op{{\bf Pic}}}
\newcommand{\gPic}{\op{\boldsymbol{\mathcal{P}ic}}}
\newcommand{\bh}{\boldsymbol{h}}
\newcommand{\bH}{\boldsymbol{\mathcal{H}}}
\newcommand{\sR}{\boldsymbol{\sf{R}}}
\newcommand{\bi}{\boldsymbol{\sf{i}}}
\newcommand{\bj}{\boldsymbol{\sf{j}}}
\newcommand{\bb}{\boldsymbol{\sf{b}}}
\newcommand{\sss}{\boldsymbol{\sf{s}}}
\newcommand{\sS}{\boldsymbol{\sf{S}}}
\newcommand{\ba}{\boldsymbol{\sf{a}}}

\newcommand{\dia}{\op{{\sf diag}}}
\newcommand{\hm}{\boldsymbol{\sf{e}}}
\newcommand{\bphi}{\boldsymbol{\phi}}
\newcommand{\bpr}{\boldsymbol{{\sf{pr}}}}
\newcommand{\bPhi}{{\boldsymbol{\Phi}}}
\newcommand{\bq}{{\boldsymbol{{\sf q}}}}

\newlength{\seqone}
\newlength{\sone}
\newlength{\seqtwo}
\newlength{\stwo}

\begin{document}

\title{Langlands duality for Hitchin systems}
\author{R.Donagi \and T.Pantev}
\date{}
\maketitle

\tableofcontents

\section{Introduction} \label{s:intro}

\subsection{Main results} \label{ss:main}

The purpose of this work is to prove the self duality of Hitchin's
integrable system: Hitchin's system for a complex reductive Lie group
$G$ is dual to Hitchin's system for the Langlands dual group
$\lan{G}$.  This statement can be interpreted at several levels.

\begin{itemize}
\item To start with, there is an isomorphism (depending on the choice
of an invariant bilinear pairing) between the bases of the Hitchin
systems for $G$ and $\lan{G}$, interchanging the discriminant
divisors.
\item The general fiber of the neutral connected component $\Higgs_0$ of
Hitchin's system for $G$ is an abelian variety. We show that it is
dual to the corresponding fiber of the neutral connected component
$\lan{\Higgs}_0$ of the Hitchin system for $\lan{G}$.
\item The non-neutral connected components $\Higgs_{\alpha}$ form
torsors over $\Higgs_0$. According to the general philosophy of
\cite{dp}, these are dual to certain gerbes. In our case, we identify
these duals as natural gerbes over $\lan{\Higgs}_0$. The gerbe
$\gHiggs$ of $G$-Higgs bundles was introduced and analyzed in
\cite{ron-dennis}.  This serves as a universal object: we show that
the gerbes involved in the duals of the non-neutral connected
components $\Higgs_{\alpha}$ are induced by $\gHiggs$.
\item More generally, we establish a duality over the complement of
the discriminant between the gerbe $\gHiggs$ of $G$-Higgs bundles and
the gerbe $\lan{\gHiggs}$ of $\lan{G}$-Higgs bundles, which
incorporates all the previous dualities. 
\item Finally, the duality of the integrable systems lifts to an
equivalence of the derived categories of $\gHiggs$ and
$\lan{\gHiggs}$.  A striking corollary is the construction of
eigensheaves for the natural Hecke operators on Higgs bundles.  These
can be viewed as "abelianized" versions, or classical limits, of the
Hecke eigensheaves predicted by the Geometric Langlands
correspondence.
\end{itemize}

\

Several special cases of our result are already known. 
The $GL(n)$ case can be traced back to Hitchin's original work \cite{hitchin}, 
which also includes some speculation that Langlands duality may explain 
the nature of the Prym varieties obtained for the other classical groups.
Hausel and Thaddeus \cite{hausel-thaddeus} considered the case $G=\op{SL}(n)$,
$\lan{G}= \mathbb{P}GL(n)$. They also showed the equality of stringy
Hodge numbers for these Langlands-dual Hitchin systems and discussed
the relationship to mirror symmetry for hyper-K\"{a}hler
manifolds. The general duality of Hitchin systems is the starting
point of Arinkin's approach \cite{arinkin} to the quasi-classical
geometric Langlands correspondence. This approach was recently
utilized by Bezrukavnikov and Braverman \cite{bb} who proved the
geometric Langlands correspondence for curves over finite fields for
$G = \lan{G} = \op{GL}_{n}$.  As explained in
\cite{bjsv}, \cite{kw}, the duality of the gerbes $\gHiggs$ and
$\lan{\gHiggs}$ proven in Theorem~\ref{thm:gerbes} was expected to
hold on physical grounds. The construction of the
abelianized Hecke eigensheaves, mostly in the case $G = \op{GL}(n)$,
was discussed by one of us (R.D.) in several talks circa 1990, 
but has not previously appeared in print.

\subsection{Outline} \label{ss:outline}

The Hitchin system $h : \Higgs \to B$ for the group $G$ and a curve
$C$ \cite{hitchin} is an integrable system whose total space is the
moduli space of semistable $K_{C}$-valued principal $G$-Higgs bundles
on $C$. The base $B$ parametrizes $K_{C}$-valued cameral covers, which
are certain \linebreak covers $p: \widetilde{C}_{b} \to C$ with an
action of $W$, the Weyl group of $G$ (we recall the definition of
cameral covers in Definition~\ref{defi:cameral.cover}).  The cameral
cover over a general $b \in B$ is a $W$-Galois cover.  For classical
simple groups $G$, the base $B$ also parametrizes appropriate spectral
covers $\bar{p}: \overline{C}_b \to C$. The Hitchin fiber $h^{-1}(b)$
can be described quite precisely,
\cite{hitchin,faltings,donagi,donagi-msri,scog,ron-dennis}. There is a
natural discriminant divisor $\Delta \subset B$ (see
section~\ref{ssec:hmap}) such that for $b \in B-\Delta$, the connected
component of $h^{-1}(b)$ is isomorphic to a certain Abelian variety
$P_b$ which can be described as a generalized Prym variety of
$\widetilde{C}_b$ (or of $\bar{C}_b$) over $C$.

\subsubsection{The base}

Our basic duality is stated in Theorem~\ref{thm:duality} below and in section
\ref{s:duality}.  The Hitchin base $B$ and the universal
cameral cover $\widetilde{\mycal{C}} \to C \times B$ depend on the
group $G$ only through its Lie algebra $\mathfrak{g}$. As a first step
towards the duality between the Hitchin system $h : \Higgs \to B$ for
$G$ and the Hitchin system $\lan{h} : \lan{\Higgs} \to \lan{B}$ for
$\lan{G}$, we note that the bases are isomorphic and the isomorphism
preserves discriminats.  The choice of a
$G$-invariant bilinear form on $\mathfrak{g}$ determines an
isomorphism $\bl : B \to \lan{B}$ between the Hitchin bases for the
Langlands-dual algebras $\mathfrak{g}$, $\lan{\mathfrak{g}}$. This
isomorphsim lifts to an isomorphism $\myell$ of the corresponding
universal cameral covers. (These isomorphisms are unique up to
automorphisms of $\widetilde{\mycal{C}} \to C \times B$: There is a
natural action of $\mathbb{C}^{\times}$ on $B$ which also lifts to an
action on $\widetilde{\mycal{C}} \to C \times B$.  The apparent
ambiguity we get in the choice of the isomorphisms $\bl, \myell$ is
eliminated by these automorphisms.)

Under these isomorphisms, cameral covers of type {\sf{B}} are
interchanged with those of type {\sf{C}}.  For the remaining simple
algebras, of types {\sf{ADEFG}}, we can choose an isomorphism of the
Lie algebras $\mathfrak{g}$ and $\lan{\mathfrak{g}}$.  We can also
choose the Cartan subalgebras to match under this isomorphism.  This
induces identifications: $\mathfrak{g} = \lan{\mathfrak{g}}$, $B=
\lan{B}$, $\widetilde{\mycal{C}}=\lan{\widetilde{\mycal{C}}}$, so we
can view the isomorphisms $\bl$, $\myell$ as automorphisms.  For the
simply laced Lie algebras (of types {\sf{ADE}}) we can then take
$\bl$, $\myell$ to be the identity. But for the Lie algebras of types
{\sf{F}}, {\sf{G}}, the natural isomorphism $\bl$ is {\em not} the
identity: it takes one cameral cover to another, interchanging short
and long roots. This phenomenon was recently noted in the
Kapustin-Witten work \cite{kw} on the geometric Langlands
Correspondence, and was used in the Argyres-Kapustin-Seiberg work
\cite{aks} on $S$-duality in $N=4$ gauge theories.

\subsubsection{The fibers}

The remainder of Theorem~\ref{thm:duality} concerns the fiberwise
duality. We need to show that the connected component $P_b$ of the
Hitchin fiber $h^{-1}(b)$ over $b \in B-\Delta$ is dual (as a
polarized abelian variety) to the connected component
$\lan{P}_{\bl(b)}$ of the corresponding fiber for the Langlands-dual
system. This is achieved by analyzing the cohomology of three group
schemes $\overline{\mathcal{T}} \supset \mathcal{T} \supset
\mathcal{T}^0$ over $C$ attached to a group $G$. Of these, the first
two were introduced in \cite{ron-dennis}, where it was shown that
$h^{-1}(b)$ is a torsor over $H^1(C, \mathcal{T})$. We recall the
definitions of these two group schemes and add the third,
$\mathcal{T}^0$, which is simply their maximal subgroup scheme all of
whose fibers are connected.  (These fibers are the connected
components of the original fibers.)  It was noted in
\cite{ron-dennis} that $\overline{\mathcal{T}}={\mathcal{T}}$ except
when $G=\op{SO}(2r+1)$ for $r \ge 1$. Dually, we note here that
$\mathcal{T} = \mathcal{T}^0$ except for $G=\op{Sp}(r)$, $r \ge 1$. In
fact, it turns out that the connected components of
$H^1(\mathcal{T}^0)$ and $H^1(\overline{\mathcal{T}})$ are dual to the
connected components of $H^1(\lan{\overline{\mathcal{T}}}),
H^1(\lan{{\mathcal{T}}}^0)$, and we are able to identify the
intermediate objects $H^1({\mathcal{T}}), H^1(\lan{{\mathcal{T}}})$
with enough precision to deduce that they are indeed dual to each
other.

Altogether we get the following theorem whose proof is discussed in
section~\ref{s:duality}.

\

\medskip

\noindent
{\bfseries Theorem~\ref{thm:duality}} \ {\em
Let $G$ be a simple complex  group, 
$\lan{G}$  the Langlands dual complex group, 
and $C$ a smooth, connected, compact curve of genus $g > 0$.
\begin{enumerate}
\item[{\bf (1)}] There is an isomorphism $\bl :
  B \to \lan{B}$, from the base of
  the $G$-Hitchin system to the base of the $\lan{G}$-Hitchin system,
  which is uniquely determined up to overall scalar, and is such that:
\begin{itemize}
\item $\bl$ preserves discriminants:
  $\bl(\Delta) = \lan{\Delta}$.
\item $\bl$ lifts to an isomorphism $\myell :
  \widetilde{\mycal{C}} \to
  \lan{\widetilde{\mycal{C}}}$ between the universal
  cameral covers of $C$.  
\end{itemize}
\item[{\bf (2)}] For $b \in B  - \Delta$, the corresponding  $G$ and $\lan{G}$ 
 Hitchin fibers are dual. The duality is given by an
  isomorphism of polarized abelian varieties
\[
\xymatrix@1@M+1pc@C+1pc{
\paiso_{b} : \hspace{-7pc} &  P_{b}
   \ar[r]^-{\cong} & \left(\lan{P}_{\bl(b)}\right)^{D},
}   
\]
where $P^{D}$ denotes the dual abelian variety of $P$.  The
isomorphism $\paiso_{b}$ is the restriction of a global duality
$\paiso$ of $\Higgs_{0}$ and $\lan{\Higgs}_{0}$ over $B -
\Delta$.
\end{enumerate}
}
\

\

\noindent
Several topological results that are needed in our proof
are collected in section~\ref{s:Prym_structure}. The main result of
that section is an explicit formula for the cocharacters of the
Hitchin Prym.

\subsubsection{Other components and duality of gerbes}

The connected components of $\Higgs$ are indexed by the fundamental
group $\pi_1(G)$. The component $\Higgs_0$ corresponding to the
neutral element parametrizes $G$-Higgs bundles which are induced from
$G_{\op{sc}}$-Higgs bundles on $C$, where $G_{\op{sc}}$ is the
universal cover group of $G$.  As shown by Hitchin \cite{Teich} and
reviewed here, the restriction of $h$ to this neutral component always
admits a section (determined by the choice of a theta characteristic,
or spin structure, on the curve $C$). 
In \ref{claim:components} we show that the group of
connected components of a Hitchin fiber $h^{-1}(b)$ is $\pi_1(G)$. In
particular, the connected components of $h^{-1}(b)$ are its
intersections with the connected components of $\Higgs$ itself.

We extend the basic duality to the non-neutral components in 
section \ref{ss:global}. The non-canonical isomorphism from 
non-neutral components of
the Hitchin fiber to $P_b$ can result in the absence of a section,
i.e. in a non-trivial torsor structure \cite{hausel-thaddeus,dp}. In
general, the duality between a family of abelian varieties $A \to B$
over a base $B$ and its dual family $A^{\vee} \to B$ is given by a
Poincare sheaf which induces a Fourier-Mukai equivalence of derived
categories. It is well known \cite{dp,bb,oren} that the Fourier-Mukai
transform of an $A$-torsor $A_{\alpha}$ is an
${\mathcal{O}}^{*}$-gerbe ${}_{\alpha}A^{\vee}$ on $A^{\vee}$. 

In our case there is indeed a natural stack mapping to $\Higgs$,
namely the moduli stack $\gHiggs^{ss}$ of semistable $G$-Higgs bundles
on $C$. Over the locus of stable bundles, the stabilizers of this
stack are isomorphic to the center $Z(G)$ of $G$ and so over the
stable locus $\gHiggs^{ss}$ is a gerbe.  The stack $\gHiggs^{ss}$ is
an open substack in the stack $\gHiggs$ of all (not necessarily
semistable) Higgs bundles. Another important stack is the stack
$\gHiggs^{\op{reg}}$ of regularized Higgs bundles which parametrizes
Higgs bundles together with a choice of a sheaf of regular
centralizers for the Higgs field. There is a forgetfull map
$\gHiggs^{\op{reg}} \to \gHiggs$. There are natural analogues of the
Hitchin map for these stacks.  The stack $\gHiggs^{\op{reg}}$ was
introduced and analyzed in {\cite{ron-dennis}} and the fibers of the
relevant Hitchin map were completely described in terms of spectral
data. The analysis in \cite{ron-dennis} shows that over $B-\Delta$ all
these stacks coincide. Since we will always work over $B-\Delta$ we
will use the clean notation $\gHiggs$ rather than the more cumbersome
$\gHiggs^{\op{reg}}$.

From {\cite{dp}} we know that every pair $\alpha \in
\pi_0(\Higgs)=\pi_1(G), \beta \in \pi_1(\lan{G}) = Z(G)^{\wedge}$
defines a $U(1)$-gerbe ${}_{\beta}\Higgs_{\alpha}$ on the connected
component $\Higgs_{\alpha}$ and that there is a Fourier-Mukai
equivalence of categories $D^b(_{\beta}\Higgs_{\alpha}) \cong
D^b({}_{\alpha}\lan{\Higgs}_{\beta})$. In our case we find that all
the $U(1)$-gerbes ${}_{\beta}\Higgs_{\alpha}$ are induced from the
single $Z(G)$-gerbe $\gHiggs$, restricted to component
$\Higgs_{\alpha}$, via the homomorphisms $\beta: Z(G) \to U(1)$. These
results culminate in Theorem~\ref{thm:gerbes}, which gives a duality
between the Higgs gerbes $\gHiggs$ and $\lan{\gHiggs}$. We also note
that the gerbe $\gHiggs_{|B-\Delta} \to \Higgs_{|B-\Delta}$ measures
the obstruction to lifting the universal $G_{\op{ad}}$-Higgs bundle to
a universal $G$-Higgs bundle.

\

In summary we get the following theorem whose proof
is given in section~\ref{ss:global}.

\

\medskip

\noindent
{\bfseries Theorem~\ref{thm:gerbes}} {\em Let $\gHiggs$ be the stack of
  $G$ Higgs bundles on a curve $C$ and let $\lan{\gHiggs}$
  be the stack of $\lan{G}$ Higgs bundles on $C$. Use
  the isomorphism $\bl : B \to \lan{B}$ from {\bf
  Theorem~\ref{thm:duality}(1)} to identify $B-\Delta$ with
  $\lan{B}-\lan{\Delta}$. Under this identification one has a canonical
  isomorphism
\begin{equation} 
\gerbeiso : \gHiggs_{|B-\Delta} \stackrel{\cong}{\longrightarrow}
(\lan{\gHiggs}_{|B-\Delta})^{D}  
\end{equation}
of commutative group stacks over $B-\Delta$. The isomorphism
$\gerbeiso$ intertwines the action of the translation operators
$\trans^{\lambda,\tilde{x}}$ on $\gHiggs_{|B-\Delta}$ with the action
of the tensorization operators $\tens^{\lambda,\tilde{x}}$ on
$(\lan{\gHiggs}_{|B-\Delta})^{D}$. 
}

\

\medskip

Here $(\lan{\gHiggs}_{|B-\Delta})^{D}$ denotes the dual of
$\lan{\gHiggs}_{|B-\Delta}$ viewed as a family of commutative group
stacks over $B - \Delta$, i.e. $(\lan{\gHiggs}_{|B-\Delta})^{D}$  is
the stack of all commutative group stack homomorphisms from
$\lan{\gHiggs}_{|B-\Delta}$ to the commutative group stack
$B\mathbb{G}_{m}$ over $B - \Delta$.

The key to the proof of Theorem~\ref{thm:gerbes} is in  our ability to
move freely among the components of $\gHiggs$ via the abelianized
Hecke correspondences. These abelianized Hecke correspondences are
carefully introduced in the Appendix, following a review of the
cameral cover yoga of \cite{ron-dennis}. Note that the abelianized
Hecke correspondences appear twice in this work. In section
\ref{ss:global} they act on $\gHiggs$ and are used for tying the
components together. On the other hand, they also occur in the
Classical Limit Conjecture \ref{con:GLC.classical}, where they act on
the Langlands dual spaces $\lan{\gHiggs}$. These are the same
correspondences, except that the group is $G$ in one case and
$\lan{G}$ in the other. The Appendix discusses $\lan{\gHiggs}$; for
the purposes of section \ref{ss:global}, the same discussion applies
but the $\lan{(\bullet)}$ superscript needs to be dropped.

\

Combining Theorem~\ref{thm:gerbes} with the abelian Fourier-Mukai
duality gives the main result of this work, an extension of
Theorem~\ref{thm:gerbes} to the case
of reductive groups:

\

\medskip

\noindent
{\bfseries Theorem~\ref{thm:reductive}} {\em Let $\mathbb{G}$ be a
connected complex reductive group, let $\lan{\mathbb{G}}$ be the
Langlands dual reductive group, and let $C$ be a smooth compact
complex curve. Write $\gHiggs_{\mathbb{G}}$ and
$\gHiggs_{(\lan{\mathbb{G}})}$ for the stacks of $K_{C}$-valued Higgs
bundles on $C$ with structure group $\mathbb{G}$ and
$\lan{\mathbb{G}}$ respectively. Then there is an isomorphism $\bl : B
\; \widetilde{\to} \; \lan{B}$ of the respective Hitchin bases which
gives an identification $B - \Delta \cong \lan{B} -
\lan{\Delta}$. Under this identification one has an isomorphism
\[ 
\gHiggs_{\mathbb{G}}  \cong \left( \gHiggs_{(\lan{\mathbb{G}})} \right)^{D}
\] 
of commutative group stacks over $B - \Delta$, intertwining the action
of translation and tensorization operators.
}

\

\noindent
The proof of this theorem is discussed in section~\ref{ss:reductive}.

\subsubsection{Derived categories and Hecke eigensheaves}

Finally, our Theorem~\ref{cor:hecke} allows one to view the
Fourier-Mukai duality in Theorems~\ref{thm:gerbes} and
\ref{thm:reductive} as a classical limit of the geometric Langlands
correspondence: under this duality, the structure sheaves of gerby
points on $\gHiggs_{0}$ are transformed into coherent sheaves on the
space $\lan{\Higgs}$ (or equivaently Higgs sheaves on $\lan{\op{{\bf
Bun}}}$) which are eigensheaves for the abelianized Hecke
correspondences. These abelianized Hecke correspondences (or
translation operators) and their action on
$\lan{\Higgs}$ are introduced in the Appendix.

\

\medskip

\noindent
{\bfseries Theorem~\ref{cor:hecke}} {\em A topologically trivial
  $G$-Higgs bundle $(V,\varphi)$ on $C$ determines an eigensheaf for
  the abelianized Hecke operators.
  Explicitly let $p : \widetilde{C} \to C$ be a cameral cover
  corresponding to a point in $B - \Delta$, and let
  $\mathcal{T}_{\widetilde{C}}$ be 
  the corresponding sheaf of regular centralizers on $C$. The choice of
  $(V,\varphi)$ gives:
\begin{itemize}
\item A  $\mathcal{T}$-torsor $\mycal{L}_{(V,\varphi)}$ on $\widetilde{C}$.
\item A representable structure morphism 
$\boldsymbol{\iota} : B\op{Aut}((V,\varphi)) \to \gHiggs_{0}$.
\end{itemize} 
Write $\mathfrak{o}_{(V,\varphi)} :=
\boldsymbol{\iota}_{*}\mathcal{O}_{B\op{Aut}((V,\varphi))}$ for the
corresponding sheaf on $\gHiggs_{0}$.  Then for every character $\mu \in
\Lambda^{\vee}$ we have a functorial isomorphism 
\[
\lanab{\mathbb{H}}^{\mu}\left(
\mathfrak{c}_{0}(\mathfrak{o}_{(V,\varphi)})\right) \cong  
\mathfrak{c}_{0}(\mathfrak{o}_{(V,\varphi)})\boxtimes \mu\left(
\mycal{L}_{(V,\varphi)}\right), 
\]
i.e. $\mathfrak{c}_{0}(\mathfrak{o}_{(V,\varphi)})$ is an abelianized Hecke
eigensheaf with eigenvalue 
$\mycal{L}_{(V,\varphi)}$. 
}

\subsection{Open problems and loose ends}

Our argument is non algebraic, in that we use cohomology of
constructible sheaves and Hodge theory to prove the duality between
families of abelian varieties. It would be nice to have a purely algebraic
argument which is local and universal in nature.

Our work deals with smooth cameral covers, establishing the Hitchin
duality over the complement of the discriminant. A major step forward
would be to formulate and prove the extension to the entire base. The
heart of the matter would presumably be an understanding of what
happens over the nilpotent cones of the dual systems.

A finer understanding of the classical limit of the Geometric
Langlands Conjecture and its relation to the abelianized version
proved here would be desirable. We discuss some of the relevant
issues, somewhat informally, in section \ref{s:cl}. The rest of the
paper does not depend on that section.

\subsection{Review and notation} \label{ss:notation}

\punkt Let $G$ be a simple complex algebraic group and let
$\lan{G}$ be the Langlads dual complex group.  The Lie algebras of $G$
and $\lan{G}$ will be denoted by $\mathfrak{g}$ and
$\lan{\mathfrak{g}}$. We fix maximal tori $T \subset G$ and $\lan{T}
\subset \lan{G}$ and denote the corresponding Cartan subalgebras by
$\mathfrak{t} \subset \mathfrak{g}$ and $\lan{\mathfrak{t}} \subset
\lan{\mathfrak{g}}$. We will also write $T_{\mathbb{R}} \subset
G_{\mathbb{R}}$ and $\lan{T}_{\mathbb{R}} \subset
\lan{G}_{\mathbb{R}}$ for the compact real forms of the complex groups
and $\mathfrak{t}_{\mathbb{R}} \subset \mathfrak{g}_{\mathbb{R}}$ and
$\lan{\mathfrak{t}}_{\mathbb{R}} \subset
\lan{\mathfrak{g}}_{\mathbb{R}}$ will denote the corresponding real
Lie algebras. We denote the space of $\mathbb{C}$-linear functions on
$\mathfrak{t}$ by $\mathfrak{t}^{\vee}$, and the space of
$\mathbb{R}$-linear functions on $\mathfrak{t}_{\mathbb{R}}$ by
$\mathfrak{t}^{\vee}_{\mathbb{R}}$. Langlands duality gives an
isomorphism $\mathfrak{t}^{\vee} = \lan{\mathfrak{t}}$ which is
compatible with the real structure. We fix this isomorphism once and
for all. We will also write $W$ for the isomorphic Weyl groups of $G$
and $\lan{G}$. 

We denote the natural pairing between $\mathfrak{t}$ and $\mathfrak{t}^{\vee}$ 
by $(\bullet, \bullet) : \mathfrak{t}^{\vee}\otimes \mathfrak{t} \to
\mathbb{C}$, while we write $\langle \bullet, \bullet \rangle :
\mathfrak{t}\otimes \mathfrak{t} \to \mathbb{C}$ for the Killing form on
$\mathfrak{t}$. For any group $G$ we have a natural collection  
of lattices
\[
\xymatrix@-2pc{
\rts_{\mathfrak{g}} & \subset & \chr_{G} & \subset &
\wts_{\mathfrak{g}} &  \subset & 
\mathfrak{t}^{\vee} \\
\crts_{\mathfrak{g}} & \subset & \cchr_{G} & \subset &
\cwts_{\mathfrak{g}} &  \subset & 
\mathfrak{t}.
}
\]
Here $\rts_{\mathfrak{g}} \subset \wts_{\mathfrak{g}} \subset
\mathfrak{t}^{\vee}$ are the root and weight lattice corresponding to the root
system on $\mathfrak{g}$ and $\chr_{G} =
\op{Hom}(T,\mathbb{C}^{\times})  = \op{Hom}(T_{\mathbb{R}},S^{1})$ is
the character lattice of $G$. Analogously, 
\[
\begin{split}
\crts_{\mathfrak{g}} &
 = \{ x \in \mathfrak{t} \, |  \, (\wts_{\mathfrak{g}},x) \subset
 \mathbb{Z} \} \cong  \wts_{\mathfrak{g}}^{\vee}  \\
\cwts_{\mathfrak{g}} &
 = \{ x \in \mathfrak{t} \, |  \, (\rts_{\mathfrak{g}},x) \subset
 \mathbb{Z} \} \cong  \rts_{\mathfrak{g}}^{\vee}
\end{split}
\]
are the coroot and coweight lattices of $\mathfrak{g}$, and
\[
\cchr_{G} = \op{Hom}(\mathbb{C}^{\times},T) =
\op{Hom}(S^{1},T_{\mathbb{R}}) = \{ x \in \mathfrak{t} \, |  \,
(\chr_{G},x) \subset 
 \mathbb{Z} \} \cong  \chr_{G}^{\vee}
\]
is the cocharacter lattice of $G$. 

The Langlands duality isomorphism $\lan{\mathfrak{t}}^{\vee} =
\mathfrak{t}$ identifies $\rts_{[\lan{\mathfrak{g}}]} =
\crts_{\mathfrak{g}}$, $\chr_{[\lan{G}]} = \cchr_{G}$, and
$\wts_{[\lan{\mathfrak{g}}]} = \cwts_{\mathfrak{g}}$. To every root
$\alpha \in \rts_{\mathfrak{g}}$ of $\mathfrak{g}$ one associates in a
standard way a coroot $\alpha^{\vee} \in \crts_{\mathfrak{g}}$, given
by the formula $(\bullet,\alpha^{\vee}) := 2\langle \alpha, \bullet
\rangle/\langle \alpha, \alpha \rangle$. Under the identification
$\rts_{[\lan{\mathfrak{g}}]} = \crts_{\mathfrak{g}}$ the root system
of $\lan{\mathfrak{g}}$ is mapped to the system of coroots of
$\mathfrak{g}$ so that the short and long roots get exchanged.
 
\

\punkt \label{ssec:hmap} 
 Let $C$ be a smooth compact complex curve of genus $g > 0$.  
Recall \cite{hitchin} that the total space
$\Higgs$ of Hitchin's system 
parametrizes semistable $K_{C}$-valued principal
$G$-Higgs bundles on $C$,
i.e. pairs $(V,\varphi)$ where $V$ is a principal $G$-bundle on $C$,
$ad(V)$ is the vector bundle associated to $V$ by the adjoint representation,
and $\varphi$ is a global section of $ad(V) \otimes K_C$.
We recall that a $G$-bundle $V$ on $C$ is semistable if 
for any parabolic subgroup $P \subset G$ 
and any $P$-subbundle $V' \subset V$,
the degree of $V'$ is $\leq 0$. 
Similarly, a Higgs bundle $(V,\varphi)$ is semistable if 
for any parabolic subgroup $P \subset G$ 
and any $P$-Higgs subbundle $(V',\varphi')$ of $(V,\varphi)$,
the degree of $V'$ is $\leq 0$.

Let $h : \Higgs \to B$ 
and 
$\lan{h} : \lan{\Higgs} \to \lan{B}$ 
denote the Hitchin integrable systems for $C$ and 
$G$ and $\lan{G}$ respectively.
The base $B$ can be identified with the
space of sections $H^{0}(C,(K_{C}\otimes \mathfrak{t})/W)$. Its points
$b \in B$ parametrize $K_{C}$-valued
cameral covers $\widetilde{C_b} \to C$ of $C$
\cite{faltings,donagi,donagi-msri,ron-dennis}.

\begin{defi} \label{defi:cameral.cover} (see \cite{ron-dennis}) 
A {\bfseries cameral cover of $C$} is a scheme $\widetilde{C}$ together with a
morphism $p : \widetilde{C} \to C$ and a $W$-action along the
fibers of $p$ satisfying: 
\begin{itemize} 
\item $p$ is finite and flat over $C$;
\item as an ${\cal{O}}_C$-module with $W$ action, 
$p_*({\cal{O}}_{\widetilde{C}})$
is locally isomorphic to 
${\cal{O}}_C \otimes {\bf C}[W]$;
\item locally with respect to the etale (or analytic) topology on $C$, 
$\widetilde{C}$ is a pull-back of the $W$-cover 
$\mathfrak{t} \to \mathfrak{t}/W$.
\end{itemize}
\end{defi}

\

\noindent
For any line bundle $L$ on $C$, the $L$-valued cameral covers are
those parametrized by \linebreak $H^{0}(C,(L \otimes \mathfrak{t})/W)$.

\

\noindent

We will say that a cameral cover 
$p_{b} : \widetilde{C}_{b} \to C$ 
has simple Galois ramification if all ramification points 
$x \in D_{b} \subset \widetilde{C}_{b}$  of $p$ 
have ramification index one. 
The ramification divisor $D_{b} \subset \widetilde{C}_{b}$ 
of a cameral cover 
$p_{b} : \widetilde{C}_{b} \to C$ 
with simple Galois ramification is a disjoint union 
$D_{b} = \coprod_{\alpha} D^{\alpha}_{b}$ of subdivisors
labeled by the roots of $\mathfrak{g}$ \cite{ron-dennis}. 
We will denote the universal
cameral cover $\widetilde{\mycal{C}} \to B\times C$. The discriminant
$\Delta \subset B$ is the locus of all $b$ for which $p_{b} :
\widetilde{C}_{b} \to C$ does not have simple Galois ramification.

\

\noindent
The Hitchin fiber $h^{-1}(b)$ for $b \in B-\Delta$ is, in general,
disconnected but all of its connected components are torsors over a
generalized Prym variety $P_{b}$ naturally associated with the cover
$p_{b} : \widetilde{C}_{b} \to C$
\cite{faltings,donagi,ron-dennis,ddp}. If $G$ is a classical group the
generalized Prym variety can also be attached to a (non-Galois)
spectral cover $\bar{p}_{b} : \bar{C}_{p} \to C$ \cite{hitchin}.

The connected components of the space $\Higgs$ are labeled
by the topological types of Higgs bundles, which in turn are labeled
by elements in $H^{2}(C,\pi_{1}(G)) = \pi_{1}(G)$.  The component
$\Higgs_{0}$ corresponding to the neutral element
parametrizes $G$-Higgs bundles which are induced from
$G_{\op{sc}}$-Higgs bundles on $C$, where $G_{\op{sc}}$ is the
universal covering group of $G$.  The restriction of $h$ to this neutral
component always admits a section, determined by the choice of a theta
characteristic on the curve $C$ and called the Hitchin section
\cite{Teich}. The construction of the Hitchin section  is also reviewed in
the proof of Lemma~\ref{lem:section}.

\

\punkt For a finitely generated abelian group $H$ we will write
$H_{\op{tors}} \subset H$ for the torsion subgroup of $H$;
$H_{\op{tf}} := H/H_{\op{tors}}$ for the maximal torsion free quotient
of $H$. 

Throughout the paper we will frequently use several duality
transformations. The most important ones are as follows:

\begin{itemize}
\item[$(\bullet)^{\vee}$:]  will denote the duality operation
  $\op{Hom}_{\mathbb{Z}}(\bullet,\mathbb{Z})$ 
  which we will be applying to free
  abelian groups of finite rank.
\item[$(\bullet)^{\wedge}$:] depending on the context will denote 
  the Pontryagin duality operation $\op{Hom}(\bullet,S^{1})$ on
  locally compact topological abelian groups 
  or the Cartier duality operation $\op{Hom}(\bullet,\mathbb{G}_{m})$
  on (complexes of) flat abelian group sheaves over a base scheme.
\item[$(\bullet)^{D}$:] depending on the context will denote the
  duality operation $R\op{Hom}(\bullet,\mathcal{O}^{\times}[1])$ on
  complexes of abelian groups or the duality operation
  $\underline{\op{Hom}}_{\op{grp-stack}}(\bullet,
  B\mathbb{G}_{m})$ on commutative group stacks over a base;
\item[$\lan{(\bullet)}$:] depending on the context will denote the
  Langlands duality on reductive groups or Lie algebras, the Langlands
  duality on maximal tori or Cartan algebras, or the induced Langlands
  duality operation on the various sheaves of regular centralizers.
\end{itemize}

\

\medskip

\noindent
{\bf Acknowledgments:} We would like to thank Edward Witten for
encouraging us to complete this work, and Emanuel Diaconescu for
helpful conversations in connection with the related project
\cite{ddp}. We thank Dima Arinkin, Dennis Gaitsgory, and Constantin
Teleman for patiently answering our questions and for providing
valuable technical advice. Finally we would like to thank the referees
whose comments and suggestions helped us greatly improve the
paper. The work of Ron Donagi was supported by the NSF grants DMS
0139799, DMS 0612992, and DMS 0908487. The work of Tony Pantev was
supported by the NSF grants DMS 0403884, DMS 0700446, and DMS
1001693. Both authors were supported by NSF Focused Research Grant DMS
0139799 and NSF Research Training Group Grant DMS-0636606.

\section{The classical limit} \label{s:cl}

Our duality of the Hitchin integrable systems for a pair of Langlands
dual groups can be interpreted as a specialization or a ``classical limit''
of the Deligne, Laumon, Beilinson-Drinfeld geometric version of the Langlands
conjecture \cite{beilinson-drinfeld-langlands}. 
In this section we will explain how this interpretation works.
We will suppress the technicalities required to
make the statements precise. The material of this section is the
motivation behind much of what we do, but it will not be needed
explicitly anywhere in the paper. The proof of the main results begins in section \ref{s:duality}.

First we will need to introduce some notation. To avoid subtleties
requiring rigidification or derived structures we will only discuss
the case of semisimple groups $G$ and $\lan{G}$.  Let as before
$\lan{\sBun}$ denote the moduli stack of principal $\lan{G}$-bundles
on $C$. We will also let $\gLoc$ denote the moduli stack of algebraic
$G$-local systems on $C$, i.e. the moduli stack of pairs $(V,\nabla)$
where $V$ is a principal algebraic $G$-bundle on $C$ and $\nabla$ is a
flat algebraic connection on $C$.  For a sheaf of algebras
$\mathcal{A}$ on an algebraic stack $X$ we will write
$D_{\op{coh}}(X,\mathcal{A})$ for the derived category of complexes of
$\mathcal{A}$-modules whose cohomology sheaves are coherent
$\mathcal{A}$-modules. The sheaves of algebras $\mathcal{A}$ that we
will be primarily interested in will be $\mathcal{A} =
\mathcal{O}_{X}$ - the structure sheaf of $X$, or $\mathcal{A} =
\mathcal{D}_{X}$ - the sheaf of algebraic differential operators on
$X$, or $\mathcal{A} = \op{Sym}^{\bullet} T_{X}$ - the symmetric
algebra on the tangent sheaf of $X$.  According to
\cite{beilinson-drinfeld-langlands} a form of the following conjecture
must hold:

\begin{con} \label{con:GLC} There exists a canonical equivalence 
of categories ({\bfseries the geometric Langlands correspondence}):
\[
\mathfrak{c} : D_{\op{coh}}(\gLoc,\mathcal{O})
\stackrel{\cong}{\longrightarrow} D_{\op{coh}}(\lan{\sBun},\mathcal{D}),
\]
which intertwines the action of the tensorization functors on
$D_{\op{coh}}(\gLoc,\mathcal{O})$ with the action of the Hecke
functors on $D_{\op{coh}}(\lan{\sBun},\mathcal{D})$. 
\end{con}

\

\noindent
The tensorization functors $W^{\mu,x} :
D_{\op{coh}}(\gLoc,\mathcal{O}) \to D_{\op{coh}}(\gLoc,\mathcal{O})$,
and the Hecke functors ${}^{L}H^{\mu,x} :
 D_{\op{coh}}(\lan{\sBun},\mathcal{D})\to
 D_{\op{coh}}(\lan{\sBun},\mathcal{D})$, 
are endofunctors of the respective categories of sheaves labeled by
the same data: pairs $(x,\mu)$, where $x \in C$ is a closed point and
$\mu \in \cchr^{+}(\lan{G}) = \chr^{+}(G)$ is a dominant cocharacter
for $\lan{G}$, or equivalently a dominant character for $G$.

Given such a pair $(x,\mu)$ one defines the tensorization functor
$W^{\mu,x}$ as
\[
\xymatrix@R-1.5pc{
W^{\mu,x} : & \hspace{-0.3in} D_{\op{coh}}(\gLoc,\mathcal{O}) \ar[r] &
D_{\op{coh}}(\gLoc,\mathcal{O})  \\
& \mycal{F} \ar@{|->}[r] & \mycal{F} {\otimes}
\rho^{\mu}\left(\mycal{V}_{|\gLoc\times \{ x\}}\right),
}
\]
and the Hecke functor ${}^{L}H^{\mu,x}$ as
\[
\xymatrix@R-1.5pc{ {}^{L}H^{\mu,x} : & \hspace{-0.3in}
D_{\op{coh}}(\lan{\sBun},\mathcal{D}) \ar[r] &
D_{\op{coh}}(\lan{\sBun},\mathcal{D})\\ & \mycal{M} \ar@{|->}[r] &
q^{\mu,x}_{!}\left((p^{\mu,x})^{*}\mycal{M}\otimes
\lan{I}^{\mu,x}\right) }
\]
i.e. as the 
integral transform on $\mathcal{D}$-modules  with kernel
$\lan{I}^{\mu,x} \in
D_{\op{coh}}(\lan{\sHeck}^{\mu,x},\mathcal{D})$.
Here:
\begin{itemize}
\item $\rho^{\mu}$ is the irreducible representation of $G$ with
  highest weight $\mu$, $\mycal{V} \to \gLoc\times C$ is the principal
  $G$-bundle underlying the universal local system, and
  $\rho^{\mu}\left(\mycal{V}_{|\gLoc\times \{ x\}}\right)$ is the
  vector bundle on $\gLoc$ associated with $\mycal{V}_{|\gLoc\times \{
  x\}}$ via the representation $\rho^{\mu}$.
\item $\lan{\sHeck}^{\mu,x}$ is the moduli stack of triples
$(V,V',\beta)$, where $V$ and $V'$ are principal $\lan{G}$-bundles on
$C$, and $\beta$ is an isomorphism of $V_{|C-\{x\}}$ with
$V'_{|C-\{x\}}$, such that for every irreducible representation $\rho$
 of $\lan{G}$, the isomorphism:
\[
\xymatrix@1@M+0.5pc{
\rho(\beta) : & \hspace{-0.3in} \rho(V)_{|C-\{x\}} \ar[r]^-{\cong} & 
\rho(V')_{|C-\{x\}}
}
\]
of vector bundles away from $x \in C$ extends to an injection
of locally free sheaves on $C$ with pole at $x$ of order bounded by
$\mu$. More precisely, on $C$ we have 
\[
\xymatrix@1@M+0.5pc{ \rho(\beta) : & \hspace{-0.3in} \rho(V) \ar@{^{(}->}[r] &
\rho(V')\left(\left\langle \mu, \lambda^{\rho}\right\rangle\cdot
x\right),
}
\]
where $\lambda^{\rho} \in \chr^{+}(\lan{G})$ denotes the highest weight of
$\rho$.
\item The stack $\lan{\sHeck}^{\mu,x}$ is equipped with two projections
\begin{equation} \label{eq:hecke}
\xymatrix{
& \lan{\sHeck}^{\mu,x} \ar[dl]_-{p^{\mu,x}} \ar[dr]^-{q^{\mu,x}} & \\
\lan{\sBun} & & \lan{\sBun}
}
\end{equation}
where $p^{\mu,x}((V,V',\beta)) := V$ and  $q^{\mu,x}((V,V',\beta)) :=
V'$. Both maps $p^{\mu,x}$ and $q^{\mu,x}$ are proper representable
locally trivial fibrations. 
\item $\lan{I}^{\mu,x}$ is the Goresky-MacPherson middle perversity
  extension $j_{!*}\left(\mathbb{C}\left[\dim
  \lan{\sHeck}^{\mu,x}\right]\right)$ of the trivial
  rank one local system on the smooth part \linebreak $j :
  \left(\lan{\sHeck}^{\mu,x}\right)^{\op{smooth}} \hookrightarrow
  \lan{\sHeck}^{\mu,x}$ of the Hecke stack.
\end{itemize}

\

\begin{rem} \label{rem:one-loop} In the recent work of Kapustin-Witten
  \cite{kw} the geometric Langlands correspondence $\mathfrak{c}$ is
  interpreted physically in two different ways. On one hand it is
  argued that Conjecture~\ref{con:GLC} is a mirror symmetry statement
  relating the $A$ and $B$-type branes on the hyper-K\"{a}hler moduli
  spaces of Higgs bundles. On the other hand Kapustin and Witten use a
  gauge theory/string duality to show that Conjecture~\ref{con:GLC}
  can be thought of as an electric-magnetic duality between
  supersymmetric four-dimensional Gauge theories with structure groups
  $G$ and $\lan{G}$ respectively. In this interpretation the functors
  $W^{\mu,x}$ and $\lan{H}^{\mu,x}$ are viewed as natural symmetry
  operations in gauge theory associated to closed loops in an
  appropriately chosen four-manifold. In this context $W^{\mu,x}$
  appear as the so called Wilson loop operators, and $\lan{H}^{\mu,x}$
  as the {}'tHooft loop operators.
\end{rem}

\

\begin{rem}
 The derived category $D_{\op{coh}}(\gLoc,\mathcal{O})$ has a natural
orthogonal spanning class of objects: the structure sheaves of all
closed (stacky) points. Each of these skyscraper sheaves is an
eigensheaf for all of the tensorization functors: if $\mathbb{V} =
(V,\nabla)$ is a $G$-local system on $C$ and if
$\mathcal{O}_{\mathbb{V}}$ is the structure sheaf of the corresponding
stacky point of $\gLoc$, then $W^{\mu,x}\left(
\mathcal{O}_{\mathbb{V}}\right) = \mathcal{O}_{\mathbb{V}}\otimes
\rho^{\mu}(\mathbb{V})_{x}$. The geometric Langlands correspondence
$\mathfrak{c}$ in Conjecture~\ref{con:GLC} therefore sends structure
sheaves of points on $\gLoc$ to Hecke eigen-$\mathcal{D}$-module on
$\lan{\sBun}$: for every $G$-local system $\mathbb{V}$, 
\[
\lan{H}^{\mu,x}\left(\mathfrak{c}\left(
\mathcal{O}_{\mathbb{V}} \right) \right) = \mathfrak{c}\left(
\mathcal{O}_{\mathbb{V}} \right)\otimes \rho^{\mu}(\mathbb{V})_{x}.
\]
The Hecke eigen-$\mathcal{D}$-module on $\lan{\sBun}$ for $G = GL_{n}$
and an irreducible local system $\mathbb{V}$, was constructed by
Drinfeld \cite{drinfeld} for $n =2$ and by Frenkel, Gaitsgory and
Vilonen \cite{fgv-glc,dennis-glc} for all $n$. The general form of the
conjecture and the case of other groups is still open.
\end{rem}

\

\noindent
An interesting feature of Conjecture~\ref{con:GLC} is that it
specializes naturally in a one parameter family. The special fiber of
the specialization is (an appropriate version of) the stack of
$K_{C}$-valued Higgs bundles on $C$ and so is naturally related to the
Hitchin system. Recall \cite{hitchin} that a $K_{C}$-valued $G$-Higgs
bundle on $C$ is a pair $(V,\theta)$, where $V$ is a principal
$G$-bundle and $\theta \in H^{0}(C,\op{ad}(V)\otimes K_{C})$ is a
$K_{C}$-valued section of the adjoint bundle of $V$. 

The specialization of the geometric Langlands conjecture has a
different nature on the two sides of the conjecture. On the left hand
side it comes from a geometric specialization of the stack $\gLoc$,
whereas on the right side it comes from specializing the
filtered sheaf of non-commutative rings $\mathcal{D}$ to its associated
sheaf of graded commutative rings:

\begin{itemize}
\item On the left hand side of the conjecture there is a natural
  geometric one parameter jump deformation of the stack $\gLoc$ of $G$-local
  systems to the stack $\gHiggs$ of  $K_{C}$-valued $G$-Higgs bundles
  \cite{simpson-hodge.filtration}.  More precisely there is a family of
  stacks $\mathcal{H} \to \mathbb{C}$ parametrized by the affine line
  $\mathbb{C}$ such that $\gLoc$ is isomorphic to the general fiber,
  and the fiber over $0 \in \mathbb{C}$ is equal to
  $\gHiggs$. Explicitly \cite{simpson-hodge.filtration} $\mathcal{H}$
  is the moduli stack of Deligne's $z$-connections on $C$,
  i.e. $\mathcal{H}$ parametrizes triples $(V,\nabla,z)$, where $\pi
  : V \to C$
  is a principal $G$-bundle on $C$, $z \in \mathbb{C}$ is a complex
  number, and $\nabla$ is a differential operator satisfying the Leibnitz 
  rule up to a factor of $z$. Equivalently, $\nabla$ is a
  $z$-splitting of the Atiyah sequence for $V$:
\[
\xymatrix{
0 \ar[r] & \op{ad}(V) \ar[r] & \mathcal{E}(V) \ar[r]^-{\sigma} & T_{C}
\ar[r] \ar@/^1pc/[l]^-{\nabla} &  
  0.
}
\]
Here $\op{ad}(V) = V\times_{\op{ad}} \mathfrak{g}$ is the adjoint
bundle of $V$, $\mathcal{E}(V) = (\pi_{*}T_{V})^{G}$ is the Atiyah
algebra of $V$, $\sigma : \mathcal{E}(V) \to T_{C}$ is the map induced
from $d\pi : T_{V} \to \pi^{*}T_{C}$, and $\nabla$ is a map of vector
bundles satisfying $\sigma\circ \nabla = z\cdot \op{id}_{T_{C}}$. 

The map $\mathcal{H} \to \mathbb{C}$ assigns to $(V,\nabla,z)$ the
complex number $z$ and is equivariant under the action of $\mathbb{C}^{\times}$
which rescales $z$. This $\mathbb{C}^{\times}$-action trivializes
$\mathcal{H}_{|\mathbb{C}^{\times}}$ and so we have a specialization family:
\[
\xymatrix@C-2pc{
\gHiggs \ar[d] &\subset & \mathcal{H} \ar[d] & \supset &
\mathcal{H}_{|\mathbb{C}^{\times}} \ar[d] & \cong & \gLoc \times
\mathbb{C}^{\times}  \\
0 & \in & \mathbb{C} & \supset & \mathbb{C}^{\times} &
}
\]
Passing to derived categories of coherent sheaves one gets a
one parameter deformation of $D_{\op{coh}}(\gLoc,\mathcal{O})$ to
$D_{\op{coh}}(\gHiggs,\mathcal{O})$ 
of the left hand side category in Conjecture~\ref{con:GLC}.
\item On the right hand side we have a natural one parameter
  jump deformation of the sheaf $\mathcal{D}$
    of rings of differential operators on $\lan{\sBun}$ to the
    symmetric algebra of the tangent sheaf on $\lan{\sBun}$. More
    precisely, the filtration by order of differential operators on
    $\mathcal{D}$ gives rise to an associated  Rees sheaf $\mathcal{R}
    \to  \lan{\sBun}\times \mathbb{C}$ - a sheaf of rings, flat over
    the affine line $\mathbb{C}$ yielding specializations
\[
\begin{split}
\mathcal{R}_{| \lan{\sBun}\times \{ z \}} & \cong \mathcal{D} \qquad
\text{for } z \neq 0 \\
\mathcal{R}_{| \lan{\sBun}\times \{ 0 \}} & \cong \op{Sym}^{\bullet} T
\end{split}
\]
Explicitly let 
$p_{1} : \lan{\sBun}\times \mathbb{C}^{\times} \to \lan{\sBun}$
be the projection on the first factor. Define a subsheaf $\mathcal{R}
\subset p_{1}^{*}\mathcal{D}$ as follows. A section of
$p_{1}^{*}\mathcal{D}$ is of the form $\sum z^{i} P_{i}$ for $P_{i}
\in \mathcal{D}$; by definition, the section is in $\mathcal{R}$ if
and only if the degree of $P_{i}$ is at most $i$, i.e. $P_{i} \in
\mathcal{D}^{\leq i}$.

Passing to derived categories of modules we get a one parameter
deformation of $D_{\op{coh}}(\lan{\sBun},\mathcal{D})$ to
$D_{\op{coh}}(\lan{\sBun},\op{Sym}^{\bullet} T)$ of the right hand
side category in Conjecture~\ref{con:GLC}. Furthermore the
abelian category of coherent $\op{Sym}^{\bullet} T$-modules on
$\lan{\sBun}$ is naturally equivalent to the category of coherent
sheaves of $\mathcal{O}$-modules on $\op{Spec}(\op{Sym}^{\bullet} T)$,
i.e. on the total space of the cotangent bundle of the smooth
algebraic stack $\lan{\sBun}$. It is well known (see
e.g. \cite{hitchin, beilinson-drinfeld-langlands}) that this total
space can be idenified naturally with the moduli stack $\lan{\gHiggs}$
of all $K_{C}$-valued $\lan{G}$-Higgs bundles on $C$. Using this
identification we can recast the  deformation of
$D_{\op{coh}}(\lan{\sBun},\mathcal{D})$ to 
$D_{\op{coh}}(\lan{\sBun},\op{Sym}^{\bullet} T)$  as a
one parameter deformation of $D_{\op{coh}}(\lan{\sBun},\mathcal{D})$
to  $D_{\op{coh}}(\lan{\gHiggs},\mathcal{O})$.
\end{itemize}

\

\noindent
One expects that the geometric Langlands correspondence $\mathfrak{c}$
also deforms along with the above of the two sides of
Conjecture~\ref{con:GLC}. In other words we expect to have a classical
limit correspondence $\mathfrak{cl} :
D_{\op{coh}}(\gHiggs,\mathcal{O}) \to
D_{\op{coh}}(\lan{\gHiggs},\mathcal{O})$ which is an equivalence of
categories and is a one parameter deformation of $\mathfrak{c} :
D_{\op{coh}}(\gLoc,\mathcal{O}) \to
D_{\op{coh}}(\lan{\sBun},\mathcal{D})$ Here the deformation from
$D_{\op{coh}}(\gLoc,\mathcal{O})$ to
$D_{\op{coh}}(\gHiggs,\mathcal{O})$ should be thought of as a
specialization from the general to the closed fiber in the stack of dg
enhanced derived categories of coherent sheaves along the fibers of
the family of spaces $\mathcal{H} \to \mathbb{C}$.  Similarly, the
deformation from $D_{\op{coh}}(\lan{\sBun},\mathcal{D})$ to
$D_{\op{coh}}(\lan{\gHiggs},\mathcal{O})$ should be viewed as a
specialization from the general to the closed fiber in the stack of dg
enhanced derived caregories of complexes of $\mathcal{R}$-modules
along the fibers of the family $\lan{\sBun}\times \mathbb{C}$

It is also expected that the  classical limit correspondence
$\mathfrak{cl}$  should intertwine suitably defined specializations of
the tensorization and Hecke functors. The specialization of the
tensorization functors is easy to describe: for each $x \in C$ and
$\mu \in \chr^{+}(G)$ as before we define a {\em classical limit}
tensorization functor as
\[
\xymatrix@R-1.5pc{ 
\mathbb{W}^{\mu,x} : & \hspace{-0.3in}
D_{\op{coh}}(\gHiggs,\mathcal{O}) \ar[r] &
D_{\op{coh}}(\gHiggs,\mathcal{O}) \\ & \mycal{F} \ar@{|->}[r] &
\mycal{F}{\otimes} \rho^{\mu}\left(\mycal{V}_{|\gHiggs\times \{
x\}}\right), }
\]
where $\mycal{V} \to \gHiggs\times C$ is the principal
  $G$-bundle underlying the universal Higgs bundle
  $(\mathcal{V},\theta) \in \Gamma\left(\op{ad}(\mathcal{V})\otimes
  p_{C}^{*}K_{C}\right)$.

\

The passage to the classical limit for Hecke functors is more
involved.  First notice that the spectral correspondence (see
e.g. \cite{donagi-msri}) gives an equivalence of the abelian category
of quasi-coherent sheaves on $\lan{\gHiggs} =
\op{tot}\left(T^{\vee}_{\lan{\sBun}}\right)$ with the abelian category
of $\Omega^{1}$-valued quasi-coherent Higgs sheaves on $\lan{\sBun}$,
that is with the abelian category of pairs $(\mathcal{E},\varphi)$,
where $\mathcal{E}$ is a quasi-coherent sheaf on $\lan{\sBun}$ and
$\varphi : \mathcal{E} \to \mathcal{E}\otimes \Omega^{1}$ is an
$\mathcal{O}$-linear map satisfying $\varphi\wedge \varphi = 0$.  In
particular we can view $D_{\op{coh}}(\lan{\gHiggs},\mathcal{O})$ as a
full subcategory of the derived category $D\text{\sf Higgs}(\lan{\sBun})$ of
quasi-coherent Higgs sheaves on $\lan{\sBun}$.

Since the Hecke functors were defined as integral transforms for
$\mathcal{D}$ modules on $\lan{\sBun}$ we can use the Higgs sheaf
interpretation of $D_{\op{coh}}(\lan{\gHiggs},\mathcal{O})$ and define
the classical limit of the Hecke functor as an integral transform for
Higgs sheaves. There is one missing ingredient for such a definition
however: we must specify a specialization of the kernel
$\mathcal{D}$-module $\lan{I}^{\mu,x}$ to a quasi-coherent Higgs sheaf
on the Hecke stack $\lan{\sHeck}^{\mu,x}$. For this one can use
the same process that we used to define the classical limit of the
right hand side of Conjecture~\ref{con:GLC}, namely the Rees
deformation of a filtered object to its associated graded. 

More
precisely, suppose that we can find a quasi-coherent  sheaf
$\lan{\mathbb{I}}^{\mu,x}$ on $\lan{\sHeck}^{\mu,x}\times C$ so that:
\begin{itemize}
\item $\lan{\mathbb{I}}^{\mu,x}$ is a module over the Rees sheaf for
  the sheaf of 
  differential operators on $\lan{\sHeck}^{\mu,x}$;
\item The restriction of $\lan{\mathbb{I}}^{\mu,x}$ to
  $\lan{\sHeck}^{\mu,x}\times \{ 
  1 \}$ is isomorphic to $\lan{I}^{\mu,x}$ as a $\mathcal{D}$-module. 
\end{itemize}
 Typically such an extension $\lan{\mathbb{I}}^{\mu,x}$ of
$\lan{I}^{\mu,x}$ will come from choosing a good filtration on
$\lan{I}^{\mu,x}$, since for any good filtration we can take
$\lan{\mathbb{I}}^{\mu,x}$ 
to be the Rees module associated with the filtration. Thus
one strategy for finding the classical limit will be to equip
$\lan{I}^{\mu,x}$ with a functorial good filtration.

The restriction $\lan{\mathfrak{I}}^{\mu,x} :=
\lan{\mathbb{I}}^{\mu,x}/z\cdot 
\lan{\mathbb{I}}^{\mu,x}$ of $\lan{\mathbb{I}}^{\mu,x}$ to
$\lan{\sHeck}^{\mu,x}\times \{ 0 \}$ is then naturally a module over
the associated graded ring $\op{\sf gr} \mathcal{D}$ ($\cong
\op{Sym}^{\bullet} T$), i.e. it is a Higgs sheaf on
$\lan{\sHeck}^{\mu,x}$ which can be viewed as the {\em classical
limit} Hecke kernel. This immediately 
gives rise to a {\em classical limit} Hecke functor
\begin{equation} \label{eq:limithecke}
\xymatrix@R-2pc{ 
\lan{\mathbb{H}}^{\mu,x} : & \hspace{-0.3in}
D_{\op{coh}}\left(\lan{\gHiggs},\mathcal{O}\right) \ar[r]  &
D_{\op{coh}}\left(\lan{\gHiggs},\mathcal{O}\right)  \\ 
& \cap & \cap \\
& D\text{\sf Higgs}\left(\lan{\sBun}\right) \ar[r] & D\text{\sf
  Higgs}\left(\lan{\sBun}\right) \\ 
& & \\
& (\mathcal{E},\varphi) \ar@{|->}[r] &
q^{\mu,x}_{!}\left((p^{\mu,x})^{*}(\mathcal{E},\varphi)  {\otimes}
(\lan{\mathfrak{I}}^{\mu,x})
\right), }
\end{equation}

\

\noindent
In general one expects that the correct filtration on the intersection
cohomology sheaf 
$\lan{I}^{\mu,x}$ comes from mixed Hodge theory. By definition
$\lan{I}^{\mu,x}$ is the middle perversity extension of the trivial
rank one local system from the smooth locus of $\lan{\Bun}$. In
particular, Saito's theory \cite{saito-mhm} implies that
$\lan{I}^{\mu,x}$ has a 
canonical structure of a mixed Hodge module. It is natural to expect
that the Hodge filtration of this mixed Hodge module will provide the
correct classical limit of $\lan{I}^{\mu,x}$. This is trivially the case
for minuscule $\mu$'s and suggests the following 

\begin{defi} \label{defi-parameter}  The
    {\bfseries classical limit Hecke kernel} $\lan{\mathfrak{I}}^{\mu,x}$ is
    the associated graded of $\lan{I}^{\mu,x}$ with respect to the
    Hodge filtration in Saito's mixed Hodge module structure on
    $\lan{I}^{\mu,x}$.
\end{defi}

Using the classical limit Hecke kernel we can now define the classical
limit Hecke functor by formula
\eqref{eq:limithecke}.  With these definition of
$\mathbb{W}^{\mu,x}$ and $\lan{\mathbb{H}}^{\mu,x}$ we can now
formulate the full version of the  classical limit geometric Langlands
conjecture: 

\begin{con} \label{con:GLC.classical} 
There exists a canonical equivalence 
of categories ({\bfseries the geometric Langlands correspondence}):
\[
\mathfrak{cl} : D_{\op{coh}}(\gHiggs,\mathcal{O})
\stackrel{\cong}{\longrightarrow} D_{\op{coh}}(\lan{\gHiggs},\mathcal{O}),
\]
which intertwines the action of the classical limit tensorization
functors $\mathbb{W}^{\mu,x}$ with the action of the classical limit
Hecke functors $\lan{\mathbb{H}}^{\mu,x}$.
\end{con}

\

This conjecture is one of the principal motivations for our
results. We will prove a version of this conjecture in
section~\ref{ss:hecke}. In the Appendix we review the abelianization
of the stack of Higgs bundles
\cite{hitchin, ron-dennis} and introduce the abelianized Hecke functors
on $D_{\op{coh}}(\gHiggs,\mathcal{O})$ which again generate a
commutative algebra of endofunctors. In section~\ref{ss:hecke} 
we show that, away from the discriminant, there exists a
Fourier-Mukai kernel on $\gHiggs\times (\lan{\gHiggs})$ which
gives an equivalence of $D_{\op{coh}}(\gHiggs,\mathcal{O})$ with
$D_{\op{coh}}(\lan{\gHiggs},\mathcal{O})$ and intertwines 
the algebra of $\mathfrak{W}^{\mu,x}$'s and
the algebra of abelianized Hecke functors. 
The precise comparison of the classical limit Hecke
functors and our abelianized Hecke functors is somewhat subtle and 
will be addressed in a forthcoming work of Arinkin and Bezrukavnikov.
Away from the discriminant they identify the algebras of 
classical and abelianized Hecke functors. This lends support to 
\ref{defi-parameter} as the correct definition of the classical limit.

\section{Duality for cameral Prym varieties} \label{s:duality}

In this section we formulate and prove the main duality result for
Hitchin Pryms.

\begin{thm} \label{thm:duality} 
Let $G$ be a simple complex  group, 
$\lan{G}$  the Langlands dual complex group, 
and $C$ a smooth, connected, compact curve of genus $g > 0$.
\begin{enumerate}
\item[{\bf (1)}] There is an isomorphism $\bl :
  B \to \lan{B}$, from the base of
  the $G$-Hitchin system to the base of the $\lan{G}$-Hitchin system,
  which is uniquely determined up to overall scalar, and is such that:
\begin{itemize}
\item $\bl$ preserves discriminants:
  $\bl(\Delta) = \lan{\Delta}$.
\item $\bl$ lifts to an isomorphism $\myell :
  \widetilde{\mycal{C}} \to
  \lan{\widetilde{\mycal{C}}}$ between the universal
  cameral covers of $C$.  
\end{itemize}
\item[{\bf (2)}] For $b \in B  - \Delta$, the corresponding  $G$ and $\lan{G}$ 
 Hitchin fibers are dual. The duality is given by an
  isomorphism of polarized abelian varieties
\[
\xymatrix@1@M+1pc@C+1pc{
\paiso_{b} : \hspace{-7pc} &  P_{b}
   \ar[r]^-{\cong} & \left(\lan{P}_{\bl(b)}\right)^{D},
}   
\]
where $P^{D}$ denotes the dual abelian variety of $P$.  The
isomorphism $\paiso_{b}$ is the restriction of a global duality
$\paiso$ of of $\Higgs_{0}$ and $\lan{\Higgs}_{0}$ over $B -
\Delta$.
\end{enumerate}
\end{thm}
\

\

\begin{rem} Given two isomorphic Lie algebras $\mathfrak{g}, \mathfrak{g}'$,
there is a canonical isomorphism $W \cong W'$ between their Weyl
groups and an isomorphism $ \mathfrak{t} \cong \mathfrak{t}'$ between
their Cartan subalgebras, taking roots to roots and intertwining the
Weyl actions.  This isomorphism is unique up to the action of $W$. For
Langlands self-dual algebras, this gives a canonical choice 
of an invariant bilinear form such that the composition 
$ \mathfrak{t} \to \mathfrak{t}^{\vee} = \lan{\mathfrak{t}} \cong \mathfrak{t}$ 
sends short roots to long roots. The resulting automorphism of $\mathfrak{t}$ is
in $W$ if $\mathfrak{g}$ is simply laced (i.e.  of type {\sf{ADE}}) but not
otherwise (types {\sf{FG}}), since it sends long roots to a multiple (greater
than $1$) of the short roots. The induced automorphism of the base $B$
of the Hitchin system then will be the identity for types {\sf{ADE}} but not
for types {\sf{FG}}. The action of these non-trivial automorphisms of the
Hitchin space was recently identified \cite{aks} as an $S$-duality
transformation in $N=4$ gauge theories compatible with a $T$-duality
transformation upon embedding in string theory.
\end{rem}

\

\noindent
{\bf Proof of Theorem~\ref{thm:duality}.} {\bf (1)} \  Recall that $B
= H^{0}(C,(K_{C}\otimes 
\mathfrak{t})/W)$ and similarly $\lan{B} = H^{0}(C,(K_{C}\otimes
\lan{\mathfrak{t}})/W)$. The choice of an invariant scalar product gives an
isomorphism $\kappa: \mathfrak{t} \stackrel{\cong}{\to} \mathfrak{t}^{\vee} =
\lan{\mathfrak{t}}$ compatible with the $W$-action and taking
reflection hyperplanes in $\mathfrak{t}$ to reflection hyperplanes in
$\lan{\mathfrak{t}}$. Therefore the isomorpism $\bl :  H^{0}(C,(K_{C}\otimes
\mathfrak{t})/W) \to H^{0}(C,(K_{C}\otimes
\lan{\mathfrak{t}})/W)$ induced from $\kappa$ will preserve
discriminants.

The isomorphism $\kappa$ globalizes to a commutative diagram of
bundles over $C$:
\[
\xymatrix@C+2pc{
K_{C}\otimes
\mathfrak{t} \ar[r]^-{\op{id}_{K_{C}}\otimes \kappa} \ar[d] & K_{C}\otimes
\lan{\mathfrak{t}} \ar[d] \\
(K_{C}\otimes\mathfrak{t})/W \ar[r]^-{\op{id}_{K_{C}}\otimes \kappa}
& (K_{C}\otimes 
\lan{\mathfrak{t}})/W
}
\]
The universal cameral covers $\widetilde{\mycal{C}}$ and
$\lan{\widetilde{\mycal{C}}}$ are the pullbacks of the columns of this
diagram by the natural evaluation maps 
\[
\xymatrix@R-1pc{
H^{0}(C,(K_{C}\otimes
\mathfrak{t})/W)\times C \ar[r] & \op{tot}(K_{C}\otimes
\mathfrak{t})/W \\
H^{0}(C,(K_{C}\otimes
\lan{\mathfrak{t}})/W)\times C \ar[r] &  \op{tot}(K_{C}\otimes
\lan{\mathfrak{t}})/W.
}
\] 
The isomorphism $\myell : \widetilde{\mycal{C}} \to
\lan{\widetilde{\mycal{C}}}$ is then induced by the isomorphism in the
top row of the diagram.

\

\medskip

\noindent
{\bf (2)} \ Very roughly our argument will proceed as follows. As a
corollary of Claim~\ref{claim:cochar} below we establish an
isomorphism of lattices $\cchr(P) \cong
\cchr({\lan{P}}^{D})$. Tensoring with $S^{1}$ we get a
diffeomorphism of the underlying real tori. From the Leray spectral
sequence we get a compatible identification of the universal covers of
$P$ and ${\lan{P}}^{D}$ showing that this diffeomorphism is a
complex analytic isomorphism. Some direct topological calculations in
section~\ref{s:Prym_structure} then show that this isomorphism
respects the natural polarizations.

In order to avoid too many exceptional cases in the
exposition, we will from now on exclude the case when $G$ is of type
${\sf{A}}_{1}$. This case is well understood and recorded in the
literature, see e.g. \cite{hausel-thaddeus}, \cite{d3hp}.

Let $b \in B - \Delta$ and let $\widetilde{C} =
\mycal{C}_{b}$ be the corresponding cameral cover of $C$, which from
now on we identify with the cameral cover
$\lan{\widetilde{\mycal{C}}}_{\bl(b)}$ via the isomorphism $\myell$. We will
denote the covering map by $p : \widetilde{C} \to C$.  Let $T \subset G$
be the maximal torus of $G$ and let $\Lambda := \cchr_{G} =
\op{Hom}(\mathbb{C}^{\times},T)$ be the corresponding cocharacter
lattice.

The two sheaves of commutative groups 
\[
\begin{split}
\overline{\mathcal{T}} & =  p_{*}(\Lambda\otimes
\mathcal{O}_{\widetilde{C}}^{\times})^{W} \\
\mathcal{T} & = \left\{ t \in \overline{\mathcal{T}}
\left|
\begin{minipage}[c]{1.5in} for every root $\alpha$ of $\mathfrak{g}$ we
  have $\alpha(t)_{|D^{\alpha}} = 1$ \end{minipage}\right. \right\}
\end{split}
\]
were introduced 
in \cite{ron-dennis}. In the above formula we identify
$\Lambda \otimes \mathbb{C}^{\times}$ 
with $T$ and we view a root $\alpha$ as a
homomorphism $\alpha : T \to
\mathbb{C}^{\times}$. The divisor $D^{\alpha} \subset
\widetilde{C}$ is the fixed divisor for the reflection
$\rho_{\alpha} \in W$ corresponding to $\alpha$. To these we now add a
third group scheme $\mathcal{T}^{o}$, the connected component of
$\mathcal{T}$. By definition, this is the maximal group subscheme of 
$\mathcal{T}$ all of whose fibers are connected.

It will be convenient to introduce real forms
$\mathcal{T}^{o}_{\mathbb{R}}$, $\mathcal{T}_{\mathbb{R}}$, and
$\overline{\mathcal{T}}_{\mathbb{R}}$ which are defined in the same
way but with the holomorphic sheaf
$\mathcal{O}^{\times}_{\widetilde{C}}$ replaced by the constant real
sheaf $S^{1}$. 

By definition we have sheaf inclusions
$\mathcal{T}^{o}_{\mathbb{R}} \subset \mathcal{T}_{\mathbb{R}} \subset
\overline{\mathcal{T}}_{\mathbb{R}}$. At any $x \in C$, which is not a
branch point of $p$, the fibers of the three sheaves are equal to each
other and non-canonically isomorphic to the compact torus $T_{\mathbb{R}} :=
\Lambda\otimes S^{1}$. At a simple branch point $s \in C$ sitting
under a ramification point in  $D^{\alpha}
\subset \widetilde{C}$, the fibers are:
\begin{equation} \label{eq:TR}
\begin{split}
\overline{\mathcal{T}}_{\mathbb{R},s} & =  \left\{ \lambda\otimes z \,
  \left| \,
  \alpha^{\vee}\left(z^{(\alpha, \lambda)}\right) = 1  \text{ in }
  T_{\mathbb{R}} 
  \right. \right\}
\\
\mathcal{T}_{\mathbb{R},s} & =  \left\{ \lambda\otimes z  \, \left| \,
  z^{(\alpha, \lambda)} = 1  \text{ in } S^{1}
 \right. \right\}
\\
\mathcal{T}^{o}_{\mathbb{R},s} & = \left\{ \lambda\otimes z  \, \left| \,
  (\alpha, \lambda) = 0  \text{ in } \mathbb{Z}
 \right. \right\}
\end{split}
\end{equation}

One of the two main results in \cite{ron-dennis} was that the Hitchin
fiber over $b$ is a torsor over the group $H^{1}(C,\mathcal{T})$
(computed in the etale or in the analytic topology).  To carry out the
comparison with the Hitchin fiber for $\lan{G}$, we will also make use
of the complex algebraic groups $H^{1}(C,\mathcal{T}^{o})$ and
$H^{1}(C,\overline{\mathcal{T}})$.  We do not know how to complete the
argument algebraically, so we resort to a topological argument which
assures us that the sheaves $\mathcal{T}^{o}$, $\mathcal{T}$,
$\overline{\mathcal{T}}$ have the same first cohomology as their real
forms $\mathcal{T}^{o}_{\mathbb{R}}$, $\mathcal{T}_{\mathbb{R}}$, and
$\overline{\mathcal{T}}_{\mathbb{R}}$.  We briefly recall the
argument, which was first observed in \cite{ddp}.

\begin{lem} \label{lem:real.forms} The natural inclusion maps of
  sheaves $\mathcal{T}^{o}_{\mathbb{R}} \subset
  \mathcal{T}^{o}$, $\mathcal{T}_{\mathbb{R}} \subset \mathcal{T}$,
  $\overline{\mathcal{T}}_{\mathbb{R}} \subset \overline{\mathcal{T}}$
induce isomorphisms on first cohomology (in the analytic topology).
\end{lem}
{\bf Proof.} \ 
The inclusion of groups $S^{1}
\subset \mathbb{C}^{\times}$ induces a natural inclusion of sheaves
\begin{equation} \label{eq-inclusion}
\nu : \Lambda \otimes S^{1} \hookrightarrow \Lambda\otimes
\mathcal{O}^{\times}_{\widetilde{C}}.
\end{equation} 
We claim that $\nu$ induces an isomorphism of commutative Lie groups
\[
\xymatrix@R-1pc{
h^{1}(\nu) : \hspace{-2pc} & H^{1}(C, (p_{*}(\Lambda\otimes S^{1}))^{W})
\ar[r] \ar@{=}[d] &  H^{1}(C, (p_{*}(\Lambda\otimes
\mathcal{O}^{\times}_{\widetilde{C}}))^{W}) \ar@{=}[d] \\
& H^{1}(C,\overline{\mathcal{T}}_{\mathbb{R}}) &
H^{1}(C,\overline{\mathcal{T}}). 
}
\]
Indeed, observe that $H^{1}(C,(p_{*}(\Lambda\otimes
S^{1}))^{W})$ is isogenous to
$H^{1}(\widetilde{C},\Lambda\otimes S^{1})^{W}$ and
similarly $H^{1}(C, (p_{*}(\Lambda\otimes
\mathcal{O}^{\times}_{\widetilde{C}}))^{W})$ is isogenous to
$H^{1}(\widetilde{C},\Lambda\otimes
\mathcal{O}_{\widetilde{C}}^{\times})^{W}$. Under these isogenies
the map $h^{1}(\nu)$ is compatible with the map
\[
H^{1}(\widetilde{C},\Lambda\otimes S^{1})^{W} \to
H^{1}(\widetilde{C},\Lambda\otimes
\mathcal{O}_{\widetilde{C}}^{\times})^{W}
\]
and so $h^{1}(\nu)$ has at most a finite kernel and a finitely
 generated
 cokernel.

Let ${\sf{Q}}$ be the cokernel of the injective map of sheaves
\eqref{eq-inclusion}.  Since 
the constant sheaf $\mathbb{C}^{\times}_{\widetilde{C}}$ has a
resolution
$\mathbb{C}^{\times}_{\widetilde{C}} \to
\mathcal{O}^{\times}_{\widetilde{C}} \to
\Omega^{1}_{\widetilde{C}}$,
and since $\mathbb{C}^{\times} = S^{1}\times \mathbb{R}$, it follows
from the snake lemma
that ${\sf{Q}}$ is isomorphic to a  sheaf of 
$\mathbb{R}$-vector
spaces on $\widetilde{C}$ which is an extension  of $\Lambda\otimes
\Omega^{1}_{\widetilde{C}}$ (considered as a sheaf of
$\mathbb{R}$-vector spaces) by the constant sheaf $\Lambda\otimes
\mathbb{R}$. Consider the push-forward of the exact sequence of
$W$-equivariant abelian sheaves
\[
1 \to \Lambda\otimes S^{1} \to \Lambda\otimes
\mathcal{O}^{\times}_{\widetilde{C}} \to {\sf{Q}} \to 0
\]
by the finite map $p : \widetilde{C} \to C$. We get a short exact
sequence 
\[
1 \to p_{*}(\Lambda\otimes S^{1}) \to p_{*}(\Lambda\otimes
\mathcal{O}^{\times}_{\widetilde{C}}) \to p_{*}{\sf{Q}} \to 0
\]
of $W$-equivariant sheaves on $C$. Taking the derived functors of $W$
invariants with coefficients in these sheaves we get a long
exact sequence of cohomology sheaves:
\[
\lesthree{(p_{*}(\Lambda\otimes S^{1}))^{W}}{(p_{*}(\Lambda\otimes
\mathcal{O}^{\times}_{\widetilde{C}}))^{W}}{(p_{*}{\sf{Q}})^{W}}
{\mycal{H}^{1}(W,p_{*}(\Lambda\otimes S^{1}))}
{\mycal{H}^{1}(W,p_{*}(\Lambda\otimes
\mathcal{O}^{\times}_{\widetilde{C}})}
\]
Since in the analytic topology  $W$ acts properly discontinuously on
$\widetilde{C}$ and since $p : \widetilde{C} \to C$ is assumed to have
simple Galois ramification, it follows  that the sheaves
$\mycal{H}^{1}(W,p_{*}(\Lambda\otimes S^{1}))$ and
$\mycal{H}^{1}(W,p_{*}(\Lambda\otimes 
\mathcal{O}^{\times}_{\widetilde{C}})$ are supported at branch points
of $p : \widetilde{C} \to C$ and that their stalks at a branch point
are the first cohomologies of $\mathbb{Z}/2$ with coefficients in
$T_{\mathbb{R}}$ and $T$ respectively (see e.g
\cite[Theorem~5.3.1]{grothendieck-tohoku}). However $T$ is the product
of $T_{\mathbb{R}}$ with an $\mathbb{R}$-vector space and so 
  the natural map $H^{1}(\mathbb{Z}/2,T_{\mathbb{R}}) \to
    H^{1}(\mathbb{Z}/2,T)$ is injective (in fact is an
    isomorphism). Therefore the natural map 
    $\mycal{H}^{1}(W,p_{*}(\Lambda\otimes S^{1})) \to
    \mycal{H}^{1}(W,p_{*}(\Lambda\otimes  
\mathcal{O}^{\times}_{\widetilde{C}})$ is also injective and so we
have a short exact sequence 
\[
1 \to \overline{\mathcal{T}}_{\mathbb{R}} \to  \overline{\mathcal{T}}
\to (p_{*}{\sf{Q}})^{W} \to 0 
\]
of sheaves on $C$. From the associated long exact sequence in cohomology
\[
\lesfour
{H^{0}(C,\overline{\mathcal{T}}_{\mathbb{R}})}
{H^{0}(C,\overline{\mathcal{T}})}
{H^{0}(C,(p_{*}{\sf{Q}})^{W})}
{H^{1}(C,\overline{\mathcal{T}}_{\mathbb{R}})}
{H^{1}(C,\overline{\mathcal{T}})}
{H^{1}(C,(p_{*}{\sf{Q}})^{W})}
{H^{2}(C,\overline{\mathcal{T}}_{\mathbb{R}})}
{h^{1}(\nu)}
\]
it follows that $H^{0}(C,(p_{*}{\sf{Q}})^{W})$ surjects onto $\ker
h^{1}(\nu)$. Since $\ker h^{1}(\nu)$ is finite and
$H^{0}(C,(p_{*}{\sf{Q}})^{W})$ is an $\mathbb{R}$-vector space it
follows that $\ker h^{1}(\nu) = 0$. Also, note that
$H^{2}(C,\overline{\mathcal{T}}_{\mathbb{R}}) =
H^{2}(C,(p_{*}(\Lambda\otimes S^{1}))^{W})$ is isogenous to
$H^{2}(C,p_{*}(\Lambda\otimes S^{1}))^{W} =
H^{2}(\widetilde{C},\Lambda\otimes S^{1})^{W} \cong (\Lambda\otimes
S^{1})^{W} $. Since $G$ is assumed simple it follows that
$H^{2}(C,\overline{\mathcal{T}}_{\mathbb{R}})$ is finite and since
$H^{1}(C,(p_{*}{\sf{Q}})^{W})$ is an $\mathbb{R}$-vector space we get
that $H^{1}(C,(p_{*}{\sf{Q}})^{W}) \to
H^{2}(C,\overline{\mathcal{T}}_{\mathbb{R}})$ is the zero map.
But from the long exact sequence this
kernel is equal to the cokernel of $h^{1}(\nu)$ which is finitely
generated as an abelian group. Hence $\op{coker}(h^{1}(\nu)) = 0$ and
$h^{1}(\nu)$ is an isomorphism.
  
Next note that by
our assumption of simple Galois ramification for $p : \widetilde{C}
\to C$ and from the definitions of $\mathcal{T}$ and
$\mathcal{T}_{\mathbb{R}}$, it follows that
$\overline{\mathcal{T}}/\mathcal{T}$ is a sheaf of groups, which is
supported at the branch points of $p$, and whose stalk at a branch
point $s$ is representable by the finite group
$\overline{\mathcal{T}}_{\mathbb{R},s}/\mathcal{T}_{\mathbb{R},s}$.
Using this fact and the isomorphism $h^{1}(\nu)$, we can compare the
long exact cohomology sequences associated with $0 \to \mathcal{T} \to
\overline{\mathcal{T}} \to \overline{\mathcal{T}}/\mathcal{T} \to 0$
and \linebreak $0 \to \mathcal{T}_{\mathbb{R}} \to
\overline{\mathcal{T}}_{\mathbb{R}} \to
\overline{\mathcal{T}}_{\mathbb{R}}/\mathcal{T}_{\mathbb{R}} \to 0$,
to conclude that $H^{1}(C,\mathcal{T}_{\mathbb{R}}) \cong
H^{1}(C,\mathcal{T})$. The same reasoning also yields the
identification $H^{1}(C,\mathcal{T}^{o}_{\mathbb{R}}) \cong
H^{1}(C,\mathcal{T}^{o})$. \ \hfill $\Box$

\

\medskip

\noindent
Recall that a root $\alpha$ for $\mathfrak{g}$  determines a homomorphism 
$(\alpha,\bullet) : \Lambda_{G} \to \mathbb{Z}$. We let $\varepsilon =
\varepsilon_{\alpha,G}$ be the positive generator of the image. We
also define $\varepsilon^{\vee} = \varepsilon^{\vee}_{\alpha,G} :=
\varepsilon_{\alpha^{\vee},\lan{G}}$.

\begin{lem} \label{lem:epsilon} 
\begin{itemize}
\item[{\bf (a)}] $\varepsilon_{\alpha,G} = 2$ when $G = \op{Sp}(r)$ and
  $\alpha$ is a long root, and $\varepsilon_{\alpha,G} = 1$ in all
  other cases. Dually $\varepsilon^{\vee}_{\alpha,G} = 2$ when $G =
  \op{SO}(2r+1)$ and $\alpha$ is a short root and
  $\varepsilon^{\vee}_{\alpha,G} = 1$ in all other cases.
\item[{\bf (b)}] $\varepsilon_{\alpha,G}$ is characterized by the property
  that $\alpha/\varepsilon_{\alpha,G}$ is a primitive vector in
  $\Lambda^{\vee}$. 
\end{itemize}
\end{lem}
{\bf Proof.} {\bf (a)} This is standard and in fact the statement for
$\varepsilon^{\vee}_{\alpha,G}$ was already noted in
\cite{ron-dennis}. The explicit argument goes as follows. Without loss
of generality we may assume that $\alpha$ is a simple root. For any
group $G$ we have $\Lambda_{G} \supset \op{coroot}_{\mathfrak{g}}$,
so for any root $\beta$ we get $(\alpha,\beta^{\vee}) \in
(\alpha,\Lambda_{G})$. When $\beta$ is simple we can read this
number from the Dynkin diagram:
\[
(\alpha,\beta^{\vee}) = \frac{2\langle \alpha, \beta \rangle}{\langle
  \beta,  \beta \rangle} = 
\begin{cases} 
2 & \text{ when } \alpha = \beta, \\
-n & \text{ when } \beta \text{ is short and $n$ edges connect
  $\alpha$ and $\beta$}, \\
-1 & \text{ when $\beta$ is long and connected to $\alpha$}, \\
0 & \text{ otherwise}. 
\end{cases}
\]
This shows that $\varepsilon_{\alpha,G} = 1$ unless all roots $\beta$
connected to $\alpha$ are short and connect to $\alpha$ by an even
number of edges. This happens only when $\alpha$ is a long root and
$\mathfrak{g}$ is of type $\text{\sf{C}}_{r}$. In the latter case we
compute that $(\alpha,\op{coroot}_{\mathfrak{g}}) = 2\mathbb{Z}$, while
$(\alpha,\op{coweight}_{\mathfrak{g}}) = \mathbb{Z}$.

\

\bigskip

\noindent
{\bf (b)} Clearly if $k$ is an integer and $\alpha/k \in
\Lambda^{\vee}$, then $k$ divides $\varepsilon_{\alpha,G}$. So we
only need to check that for a long root $\alpha$ of $Sp(r)$ we have
that $\alpha/2$ is in the weight lattice. But the root lattice for
type $\text{\sf{C}}_{r}$ has generators $e_{1}-e_{2}, \ldots, e_{r-1} - e_{r},
2e_{r}$ and the weight lattice has generators $e_{1}, e_{2}, \ldots,
e_{r}$. Again, up to a $W$-action, we may assume that $\alpha$ is a
simple root, i.e. that $\alpha = 2 e_{r}$. Thus $\alpha/2 = e_{r}$
which is indeed in the weight lattice. \ \hfill $\Box$

\

\bigskip

Consider the cover $p : \widetilde{C} \to C$. We will denote the
branch locus of this cover by $S = \{ s_{1}, \ldots, s_{b} \} \subset
C$. For each $i = 1, \ldots, b$ we will write $\alpha_{i}$ for the
root of $\mathfrak{g}$ determined (up to $W$ action) by $s_{i}$. Let
$\varepsilon_{i} := \varepsilon_{\alpha_{i},G}$ and
$\varepsilon_{i}^{\vee} := \varepsilon^{\vee}_{\alpha_{i},G}$. 
 We write $\jmath : U 
\hookrightarrow C$ for the inclusion of the complement, and $p^{o} :
p^{-1}(U) \to U$ for the unramified part of $p$. Define a local system
$A$ on $U$ by $A := (p^{o}_{*}\Lambda)^{W}$.  Note that the fibers of
$A$ are non-canonically isomorphic to $\Lambda$. We can also consider
$\lan{A} = (p^{o}_{*}(\lan{\Lambda}))^{W}$. The canonical
identification $\lan{\Lambda} = \Lambda^{\vee} =
\op{Hom}(\Lambda,\mathbb{Z})$ gives also an identification $\lan{A} =
A^{\vee} = \underline{\op{Hom}}(\Lambda,\mathbb{Z})$.

\begin{lem} \label{lem:diagram} There  are natural isomorphisms of
  sheaves $\mathcal{T}_{\mathbb{R}}^{o} \cong (\jmath_{*}A)\otimes
  S^{1}$ and \linebreak $\overline{\mathcal{T}}_{\mathbb{R}} \cong
  \jmath_{*}(A\otimes S^{1})$, while $\mathcal{T}_{\mathbb{R}}$ is
  determined by the commutative diagram:
\[
\xymatrix@-1pc{
& & 0 \ar[d] & 0 \ar[d] &  \\
0 \ar[r] & \mathcal{T}_{\mathbb{R}}^{o} \ar[r] \ar@{=}[d] &
  \mathcal{T}_{\mathbb{R}} \ar[r] \ar[d] & \oplus_{i = 1}^{b}
  \mathbb{Z}/\varepsilon_{i} \ar[d] \ar[r] & 0 \\
0 \ar[r] & \mathcal{T}_{\mathbb{R}}^{o} \ar[r]  &
  \overline{\mathcal{T}}_{\mathbb{R}} \ar[r]^-{\xi}  \ar[d] &
  \oplus_{i = 1}^{b} 
  \mathbb{Z}/\varepsilon_{i}\varepsilon^{\vee}_{i}
  \ar[d]^-{\epsilon^{\vee}}  \ar[r] & 0 \\
& & \oplus_{i = 1}^{b} 
  \mathbb{Z}/\varepsilon^{\vee}_{i} \ar[d] \ar@{=}[r] & \oplus_{i = 1}^{b} 
  \mathbb{Z}/\varepsilon^{\vee}_{i} \ar[d] & \\
& & 0 & 0 &
}
\]
\end{lem}
{\bf Proof.} Clearly the sheaves $\mathcal{T}_{\mathbb{R}}^{o}$,
$\overline{\mathcal{T}}_{\mathbb{R}}$, $(\jmath_{*}A)\otimes S^{1}$
and $\jmath_{*}(A\otimes S^{1})$ coincide on $U$. Since $A$ is a local
system we have (see Section~\ref{ss:ls_cohomology})
$(\jmath_{*}A)_{s_{i}} \cong 
\Lambda^{\rho_{i}}$, where $\rho_{i}(\lambda) = \lambda -
(\alpha_{i},\lambda)\alpha_{i}^{\vee}$ is the
reflection corresponding to $\alpha_{i}$. Similarly
$(\jmath_{*}(A\otimes S^{1}))_{s_{i}} \cong (\Lambda \otimes
S^{1})^{\rho_{i}}$. The formula for the reflection $\rho_{i}$ and
\eqref{eq:TR}  now imply that
$\mathcal{T}_{\mathbb{R},s_{i}}^{o} = (\jmath_{*}A)_{s_{i}}\otimes
S^{1}$ and $\overline{\mathcal{T}}_{\mathbb{R},s_{i}} =
(\jmath_{*}(A\otimes S^{1}))_{s_{i}}$. 

On the stalk at $s_{i}$ the map $\xi$ is given by
\[
\xi(\lambda\otimes z) :=
z^{\left(\alpha_{i}/\varepsilon_{i},\lambda\right)} \in
\boldsymbol{\mu}_{\varepsilon_{i}\varepsilon^{\vee}_{i}} \subset
S^{1}. 
\]
Here $\boldsymbol{\mu}_{\varepsilon_{i}\varepsilon^{\vee}_{i}} \subset
S^{1}$ denotes the roots of unity of order
$\varepsilon_{i}\varepsilon^{\vee}_{i}$. Since
$\varepsilon_{i}\varepsilon^{\vee}_{i}$ divides $2$, we have a natural
identification
$\boldsymbol{\mu}_{\varepsilon_{i}\varepsilon^{\vee}_{i}} =
\mathbb{Z}/\varepsilon_{i}\varepsilon^{\vee}_{i}$. 

>From \eqref{eq:TR} we now deduce that $\mathcal{T}_{\mathbb{R}}^{o} =
\ker(\xi)$, and that $\mathcal{T}_{\mathbb{R}} =
\ker(\epsilon^{\vee}\circ \xi)$, where 
\[
\epsilon^{\vee} : \oplus_{i =
  1}^{b} \mathbb{Z}/\varepsilon_{i}\varepsilon^{\vee}_{i} \to \oplus_{i =
  1}^{b} \mathbb{Z}/\varepsilon_{i}
\]
 is the map which multiplies
the $i$-th summand by $\varepsilon_{i}^{\vee}$.  \ \hfill $\Box$

\

\medskip

\begin{claim} \label{claim:components}
\begin{itemize}
\item[{\bf (i)}]
The connected components $P^{o}$, $P$, $\overline{P}$ of
  $H^{1}(C,\mathcal{T}^{o})$, $H^{1}(C,\mathcal{T})$,
  $H^{1}(C,\overline{\mathcal{T}})$ are abelian varieties. The natural maps 
$H^{1}(C,\mathcal{T}^{o}) \to H^{1}(C,\mathcal{T}) \to
  H^{1}(C,\overline{\mathcal{T}})$ and $P^{o} \to P \to \overline{P}$
  are surjective. 
\item[{\bf (ii)}] The group of connected components of
  $H^{1}(C,\mathcal{T}^{o})$ is $\mathbb{Z}/2$ for $G = \op{Sp}(r)$
  and is $\pi_{1}(G)$ otherwise. 
\item[{\bf (iii)}] The group of connected components of
  $H^{1}(C,\mathcal{T})$ is always $\pi_{1}(G)$, so the components of
  the fiber of $h : \Higgs \to B$ are in one-to-one
  correspondence with the components of the $G$-Hitchin $\op{{\bf
  Higgs}}$ system itself.
\end{itemize}
\end{claim} 
{\bf Proof.} {\bf (i)} We already noted that the connected component
of $H^{1}(C,\overline{\mathcal{T}})$ is an abelian variety, and that 
$\overline{\mathcal{T}}/\mathcal{T}$ and $\mathcal{T}/\mathcal{T}^{o}$
have finite supports and fibers which are finite groups. It follows
that 
$H^{1}(C,\mathcal{T}^{o})$ and  $H^{1}(C,\mathcal{T})$ map to
$H^{1}(C,\overline{\mathcal{T}})$ surjectively with finite kernels. In
particlular the connected components of $H^{1}(C,\mathcal{T}^{o})$ and
$H^{1}(C,\mathcal{T})$ are also abelian varieties. 

\

\medskip

\noindent
{\bf (ii)} The case $g=1$ is elementary: the cameral cover is a product of $C$ with a finite scheme, so there is no ramification, $U=C$, and the result is straightforward. So assume $g>1$. We consider the exponential sequence 
$0 \to \mathbb{Z} \to \mathbb{R} \to S^{1} \to 0$ 
of constant sheaves on $C$. Tensoring with $\jmath_{*}A$ gives 
\[
\xymatrix@1{
\underline{\op{Tor}}_{1}(\jmath_{*}A,S^{1})
\ar[r]^-{\partial}  &
\jmath_{*}A \ar[r] & (\jmath_{*}A)\otimes \mathbb{R} \ar[r] & 
\mathcal{T}^{o}_{\mathbb{R}} \ar[r] & 0.
}
\]
The sheaf
$\underline{\op{Tor}}_{1}(\jmath_{*}A,S^{1})$ is
supported on $S$ while $\jmath_{*}A$ has no compactly supported
sections. Therefore $\partial = 0$ and we get a short exact sequence
\[
\xymatrix@R-1pc{
0 \ar[r] & H^{1}(\jmath_{*}A\otimes \mathbb{R})/H^{1}(\jmath_{*}A) \ar[r]
\ar@{=}[d] & 
  H^{1}(C,\mathcal{T}^{o}_{\mathbb{R}}) \ar[r] & H^{2}(C,\jmath_{*}A)
  \ar[r] & 0. \\ 
& H^{1}(C,\mathcal{T}^{o}_{\mathbb{R}})^{o} & & &
}
\]
The group of connected components of
$H^{1}(C,\mathcal{T}^{o}_{\mathbb{R}})$ is therefore
$H^{2}(C,\jmath_{*}A)$, which can be identified (see
Lemma~\ref{lem:Hij}) with $H^{1}(U,A^{\vee})_{\op{tor}}^{\wedge}$.

The calculation performed in Corollary~\ref{cor:torsion} below gives that:
\[
H^{1}(U,A^{\vee})_{\op{tor}} = \left(
\frac{(\Lambda^{\vee})^{b}}{(1-\rho_{1},
  1-\rho_{2},\ldots, 
  1-\rho_{b})\Lambda^{\vee}}\right)_{\op{tor}}.
\]
Now for any inclusion of lattices $N \subset M$, the torsion in $M/N$ is
equal to the quotient $N'/N$ where $N' := \{ m \in M | k\cdot m \in N
\text{ for some } k \neq 0 \in \mathbb{Z} \}$ is the saturation of $N$
in $M$. In our case $N = \Lambda^{\vee}$, while the saturation is
\[
N' = \left\{ \xi \in \mathfrak{t}^{\vee} \, \left| \,
(\xi,\alpha^{\vee}) \cdot \alpha \in 
\Lambda^{\vee} \text{ for every root }
\alpha \right. \right\}.  
\]
This holds since our genericity assumption on $\widetilde{C}$ implies
that $D^{\alpha} \subset \widetilde{C}$ is non-empty for every root
$\alpha$. Using the characterization of $\varepsilon_{\alpha,G}$ in
Lemma~\ref{lem:epsilon}{\bf (b)}, we see that $\xi \in N'$ if and only if
$\varepsilon_{\alpha,G}(\xi,\alpha^{\vee}) \in \mathbb{Z}$ for all
roots $\alpha$. In case $G =
\op{Sp}(r)$ we see that $N'$ contains the weight lattice as a
sublattice of index two. Explicitly the weight lattice is generated by
$e_{1}, \ldots, e_{r}$ and $N'$ is spanned by the $e_{i}$'s and the
additional vector $\frac{1}{2} \sum_{i = 1}^{r} e_{i}$. For all other
$G$, all $\varepsilon$'s are $1$, so $N'$ is the weight lattice. We
conclude that
\[
H^{1}(U,A^{\vee})_{\op{tor}} = \begin{cases}
\mathbb{Z}/2 & \text{ when } G = \op{Sp}(r) \\
\op{wts}_{\mathfrak{g}}/\Lambda^{\vee}_{G} = \pi_{1}(G)^{\wedge} &
\text{ for all other } G.
\end{cases}
\]
This completes the proof of {\bf (ii)}. 

\

\medskip

\noindent 
{\bf (iii)} As we saw in Lemma~\ref{lem:diagram} we have
$\mathcal{T}^{0} = \mathcal{T}$ except when $G = \op{Sp}(r)$. In the latter
case the fiber of the Hitchin map was shown to be connected in
\cite{hitchin}, using an interpretation via spectral covers.  
\ \hfill $\Box$

\

\bigskip

\noindent
In general for any compact torus $H$ we define the cocharacter lattice
$\cchr(H)$ as the lattice of homomorphisms from the circle
$S^{1}$ to $H$. We recover $H$ as $\cchr(H)\otimes S^{1}$.

\begin{claim} \label{claim:cochar}
\begin{itemize}
\item[{\bf (i)}] There is a natural isomorphism
$\cchr(P^{o}) = H^{1}(C,\jmath_{*}A)_{\op{tf}}$.
\item[{\bf (ii)}] There is a natural isomorphism
  $\cchr(\overline{P})  = H^{1}(C,\jmath_{*}A^{\vee})^{\vee}$. 
\item[{\bf (iii)}]  The map $\zeta : \cchr(\overline{P}) \to
  \oplus_{i = 1}^{b} \mathbb{Z}/\varepsilon_{i}\varepsilon_{i}^{\vee}$
  induced from the map $\xi$ in Lemma~\ref{lem:diagram}
  satisfies
\[
\begin{split}
\ker(\zeta) & = \cchr(P^{o}) \\
\ker(\epsilon^{\vee}\circ \zeta) & = \cchr(P).
\end{split}
\]
\end{itemize}
\end{claim}
{\bf Proof.} {\bf (i)} By Lemma~\ref{lem:diagram}  we know
that $\mathcal{T}^{o}_{\mathbb{R}} = (\jmath_{*}A)\otimes S^{1}$. As
in the proof of Claim~\ref{claim:components}~{\bf (ii)}, we 
tensor the exponential sequence for $S^{1}$ by $\jmath_{*}A$ and we get 
\[
\xymatrix@1{
\underline{\op{Tor}}_{1}(\jmath_{*}A,S^{1})
\ar[r]^-{\partial}  &
\jmath_{*}A \ar[r] & \jmath_{*}A\otimes \mathbb{R} \ar[r] & 
\mathcal{T}^{o}_{\mathbb{R}} \ar[r] & 0.
}
\]
Again the sheaf
$\underline{\op{Tor}}_{1}(\jmath_{*}A,S^{1})$ is
supported on $S$ while $\jmath_{*}A$ has no compactly supported
sections. Therefore $\partial = 0$ and we get a short exact sequence
\[
\xymatrix@R-1pc{
0 \ar[r] & H^{1}(\jmath_{*}A)\otimes S^{1} \ar[r] &
  H^{1}(C,\mathcal{T}^{o}) \ar[r] & H^{2}(C,\jmath_{*}A) \ar[r] & 0. 
}
\]
Since $H^{2}(C,\jmath_{*}A)$ is finite and $H^{1}(\jmath_{*}A)\otimes
S^{1}$ is connected, it follows that $P^{o} = H^{1}(\jmath_{*}A)\otimes
S^{1}$, or equivalently $\cchr(P^{o}) =
H^{1}(C,\jmath_{*}A)_{\op{tf}}$.  

\

\medskip

\noindent
{\bf (ii)} \ Start with the Leray spectral sequence (aka
Mayer-Vietoris) for the 
inclusion $\jmath : U 
\subset C$ and the sheaf $A$. It gives
\[
0 \to H^{1}(C,\jmath_{*}A) \to H^{1}(U,A) \to Q \to 0,
\]
where $Q = \ker(H^{0}(R^{1}\jmath_{*}A) \to H^{2}(\jmath_{*}A))$ (see
Section~\ref{ss:ls_cohomology} for details). We tensor this sequence
with $S^{1}$ and map to the 
Leray sequence for $\jmath$ and $A\otimes S^{1}$:
\[
\xymatrix@R-1pc{
Q_{\op{tor}} \ar[d] \ar[r] & H^{1}(C,\jmath_{*}A)\otimes S^{1} \ar[r]
\ar[d] & H^{1}(U,A)\otimes S^{1} \ar[r] \ar@{^{(}->}[d] & Q\otimes
S^{1} \ar[r] 
\ar@{=}[d]  
& 0 \\
0 \ar[r] & H^{1}(C, \overline{\mathcal{T}}_{\mathbb{R}}) \ar[r] &
  H^{1}(U,A\otimes S^{1}) \ar[r] & Q\otimes S^{1} \ar[r] & 0.
}
\]
Recall that $\overline{P}$ is the connected component of
$H^{1}(C, \overline{\mathcal{T}}_{\mathbb{R}})$ and that by 
  Lemma~\ref{lem:diagram} we have
$\overline{\mathcal{T}}_{\mathbb{R}} = \jmath_{*}(A\otimes S^{1})$. It
  follows that $\overline{P}$ can be identified with the image
$\op{im}\left[ P^{o} \to H^{1}(U,A)\otimes S^{1} \right]$. In
  particular, on character lattices we get 
\[
\op{char} \overline{P} = \op{im} \left[ H^{1}(U,A)^{\vee} \to
  \op{char} P^{o} = H^{1}(C,\jmath_{*}A)^{\vee} \right]  =
  H^{1}(C,\jmath_{*}A^{\vee})_{\op{tf}}, 
\]
where the last equality follows from Corollary~\ref{cor:image}.

\

\medskip

\noindent
{\bf (iii)} This is immediate from parts {\bf (i)} and {\bf (ii)} and
the commutative diagram of sheaves in Lemma~\ref{lem:diagram}. \
\hfill $\Box$  

\

\bigskip

\noindent
The statement of the previous claim can be organized in a diagram:
\[
\xymatrix@R-1pc@C-2pc{
H^{1}(C,\jmath_{*}A)_{\op{tf}} \ar[d] & \subset & \cchr(P) \ar[d] &
\subset  & H^{1}(C,\jmath_{*}A^{\vee})^{\vee} \ar[d] \\
0 & \subset & \oplus_{i = 1}^{b} \mathbb{Z}/\varepsilon_{i}^{\vee} &
\subset &
\oplus_{i = 1}^{b} \mathbb{Z}/\varepsilon_{i}\varepsilon_{i}^{\vee}. 
}
\]
Writing the analogous diagram for $\lan{G}$ and dualizing gives
\[
\xymatrix@R-1pc@C-2pc{
H^{1}(C,\jmath_{*}A)_{\op{tf}} \ar[d] & \subset &
\cchr(\lan{P})^{\vee} \ar[d] & 
\subset  & H^{1}(C,\jmath_{*}A^{\vee})^{\vee} \ar[d] \\
0 & \subset & \oplus_{i = 1}^{b} \mathbb{Z}/\varepsilon_{i}^{\vee} &
\subset & 
\oplus_{i = 1}^{b} \mathbb{Z}/\varepsilon_{i}\varepsilon_{i}^{\vee}. 
}
\]
This gives the desired isomorphism of lattices: $\cchr(P) \cong
\cchr(\lan{P})^{\vee}$. Tensoring with $S^1$, we get a diffeomorphism 
of the underlying real manifolds.
The Leray spectral
sequence for $p : \widetilde{C} \to C$ allows us to identify the
universal covers of $P$ and $\lan{P}$ with the complex vector spaces
$H^{1}(\widetilde{C},\Lambda\otimes \mathcal{O}_{\widetilde{C}})^{W}$
and $H^{1}(\widetilde{C},\Lambda^{\vee}\otimes
\mathcal{O}_{\widetilde{C}})^{W}$, showing that the diffeomorphism 
is an isomorphism of complex manifolds.
Finally, we need to check that the two polarizations correspond 
to each other. This amounts to the compatibility of our isomorphism with 
the Poincare duality map for the cohomologies of $A$ and $A^{\vee}$ on $U$.
We defer this calculation to Lemma~\ref{lem:Hij} and
Corollary~\ref{cor:image} in section 
\ref{s:Prym_structure}. We have thus produced an isomorphism
between the polarized abelian varieties $P$ and ${\lan{P}}^{D}$. 
This completes the proof of Theorem A.
\
\hfill $\Box$

\section{Duality for Higgs gerbes} \label{s:full}

In this section we extend the duality of cameral Pryms established in
Theorem~\ref{thm:duality} to a more general duality for the stacks of
Higgs bundles, considered as families of stacky groups over the
Hitchin base. 

\subsection{Triviality} \label{ss:triviality}

The moduli stack of $G$-Higgs bundles on $C$ was defined and studied in
detail in \cite{ron-dennis}.  We briefly
recall the highlights of that discussion. Let
\[
\xymatrix{
\widetilde{\mycal{C}} \ar[rr]^-{p} \ar[dr]_-{\tilde{\pi}} & & B\times
C \ar[dl]^-{\pi}  \\
& B & }
\] 
denote  the universal cameral cover. A sheaf  $\mathcal{T}$ of abelian
groups on $B\times C$ was introduced in \cite{ron-dennis}. It was
defined as:
\begin{equation} \label{eq:sheafT}
\mathcal{T}(U)  = \left\{ t \in \Gamma\left(p^{-1}(U),
\Lambda\otimes \mathcal{O}_{\widetilde{\mycal{C}}}^{\times}\right)^{W}
\left|
\begin{minipage}[c]{1.5in} for every root $\alpha$ of $\mathfrak{g}$ we
have $\alpha(t)_{|D^{\alpha}} = 1$ \end{minipage}\right. \right\}.
\end{equation}
It was shown in \cite{ron-dennis} that relatively over the Hitchin
base the stack $\gHiggs$ of Higgs bundles is a banded
$\mathcal{T}$-gerbe. Informally, this means that the sheaf of groups
for which $\gHiggs$ is a gerbe is $\mathcal{T}$ itself rather than a
more general sheaf of groups which is only locally isomorphic to
$\mathcal{T}$. Equivalently, when viewed as a stack over the Hitchin
base, $\gHiggs$ is a torsor over the commutative group stack
$\gTors_{\mathcal{T}}$ parametrizing $\mathcal{T}$-torsors along the
fibers of $\pi : B\times C \to B$. This description is valid for Higgs
bundles over a base variety of arbitrary dimension.  When the base is
a compact curve, the picture can be made even more precise.

\begin{lem} \label{lem:section} Let $C$ be a smooth compact curve and
  let $\gHiggs$ be  
  the moduli stack of $G$-Higgs bundles on $C$. 
\begin{itemize}
\item[{\bf (a)}] The commutative group stack $\gTors_{\mathcal{T}}$
  parametrizing the $\mathcal{T}$-torsors along the fibers of $\pi$
  is isomorphic to the Picard stack
associated (see \cite[Section 
1.4 of Expos\'{e} XVIII]{sga4} and \cite[Section 6]{laumon}) 
with the amplitude one complex
$R^{\bullet}\pi_{*}\mathcal{T}[1]$ of abelian sheaves on $B$.
\item[{\bf (b)}] There exists an isomorphism $\gHiggs \cong
  \gTors_{\mathcal{T}}$ of stacks over $B$. 
\end{itemize}
\end{lem}
{\bf Proof.} {\bf (a)} This follows from the fact that $\pi : B\times C \to
B$ is smooth of relative dimension one, the standard description of
torsors in terms of \v{C}ech cocylces, and the definition (see
\cite[Section 1.4 of Expos\'{e} XVIII]{sga4} and \cite[Section
6]{laumon}) of a Picard stack associated with an amplitude one complex
of abelian sheaves.

\

\noindent
{\bf (b)} By \cite[Theorem~4.4]{ron-dennis} the stack $\gHiggs$ is a
torsor over the commutative group stack $\gTors_{\mathcal{T}}$. Thus
to get the isomorphism $\gHiggs \cong
  \gTors_{\mathcal{T}}$, it suffices to show that the stacky Hitchin
  fibration $\bh : \gHiggs \to B$ admits a section. This is due
  to Hitchin who in \cite{Teich}  constructed a family of holomorphic
  sections of $\bh : \gHiggs \to B$ induced from a Kostant section of
  the Chevalley map $\mathfrak{g} \to \mathfrak{g}/G \cong
  \mathfrak{t}/W$.  Since these sections play a prominent role in what
  follows, we briefly recall Hitchin's construction. 

Let $\{ \boldsymbol{e}, \boldsymbol{f}, \boldsymbol{g} \} \subset
  \mathfrak{g}$ be any prinicipal $\mathfrak{s}\mathfrak{l}_{2}$
  triple in $\mathfrak{g}$. This means that $\boldsymbol{e},
  \boldsymbol{f}, \boldsymbol{g}$ span a Lie subalgebra in
  $\mathfrak{g}$ isomorphic to
  $\mathfrak{s}\mathfrak{l}_{2}(\mathbb{C})$, and that
  $\boldsymbol{e}$ and $\boldsymbol{f}$ are regular nilpotent elements
  of $\mathfrak{g}$. Let $\mathfrak{C}(\boldsymbol{e}) \subset
  \mathfrak{g}$ be the centralizer of the element $\boldsymbol{e}$ in
  the algebra $\mathfrak{g}$. Consider the linear coset
  $\boldsymbol{f} + \mathfrak{C}(\boldsymbol{e}) = \{ \boldsymbol{f} +
  x | \, x \in \mathfrak{C}(\boldsymbol{e}) \}$. In \cite{kostant}
  Kostant showed that the map $\mathfrak{g} \to \mathfrak{t}/W$
  becomes an isomorphism, when restrictited to $\boldsymbol{f} +
  \mathfrak{C}(\boldsymbol{e})$. Thus $\boldsymbol{f} +
  \mathfrak{C}(\boldsymbol{e})$ is a section for the Chevalley
  projection. By construction this section consists of regular
  elements in $\mathfrak{g}$ and is a generalizaton of the rational
  canonical form of a matrix. 

Fix a Kostant section $\boldsymbol{k} : \mathfrak{t}/W \to
\mathfrak{g}$, corresponding to an
$\mathfrak{s}\mathfrak{l}_{2}$-triple in $\mathfrak{g}$.  The
inclusion of the $\mathfrak{s}\mathfrak{l}_{2}$-triple in
$\mathfrak{g}$ induces a group homomorphism $\dia :
\op{SL}_{2}(\mathbb{C}) \to G$.  Let $\zeta \in \op{Pic}^{g-1}(C)$,
$\zeta^{\otimes 2} = K_{C}$ be a theta characteristic on $C$. Consider
the frame bundle $\underline{\op{Isom}}(\zeta\oplus
\zeta^{-1},\mathcal{O}^{\oplus 2})$ of the vector bundle $\zeta\oplus
\zeta^{-1}$ on $C$. This is a principal
$\op{SL}_{2}(\mathbb{C})$-bundle which via $\dia$ gives rise to an
associated principal $G$ bundle $P :=
\underline{\op{Isom}}(\zeta\oplus \zeta^{-1},\mathcal{O}^{\oplus
2})\times_{\dia} G$ on $C$. Recall that the Hitchin base $B$ is the
space of sections of the bundle $(K_{C}\otimes \mathfrak{t})/W$ on
$C$. Let $\boldsymbol{U}$ denote the total space of the bundle
$(K_{C}\otimes \mathfrak{t})/W$, and let $u : \boldsymbol{U} \to C$ be the
natural projection. We have
\[
\op{ad}(u^{*}P) = u^{*}\op{ad}(P) =
u^{*}\underline{\op{Isom}}(\zeta\oplus
\zeta^{-1},\mathcal{O}^{\oplus 2})\times_{\op{ad}(\dia)}
\mathfrak{g}.
\]
In \cite{Teich} Hitchin checked that the Kostant section
$\boldsymbol{k} : \mathfrak{t}/W \to \mathfrak{g}$ induces a well
defined section $\varphi \in H^{0}(\boldsymbol{U},
\op{ad}(u^{*}P)\otimes u^{*}K_{C})$ and hence a $u^{*}K_{C}$-valued
Higgs bundle $(u^{*}P,\varphi)$ on $\boldsymbol{U}$. Pulling back this
Higgs bundle by the sections $b \in B = H^{0}(C,\boldsymbol{U})$, one
gets a family of Higgs bundles on $C$, parametrized by $B$. We will
call the resulting section of $\bh : \gHiggs \to B$ the {\em Hitchin
  section} and denote it by $\mathfrak{v} : B \to \gHiggs$. \   \hfill $\Box$

\

\subsection{Stabilizers, components, and universal bundles}
\label{ss:universal} 

From now on, we restrict our attention to the open substack of
$\gHiggs$ consisting of stable $G$-Higgs bundles whose authomorphism
group is the minimal possible, i.e. coincides with the center of $G$.
The Hitchin fiber for a cameral cover in $B - \Delta$ consists only of
stable Higgs bundles, and in fact each Higgs bundle in such fiber has
minimal automorphism group:

\begin{lem} \label{lem:auto}  \begin{itemize}
\item[{\bf (i)}] $\gHiggs_{|B-\Delta}$ is a smooth Deligne-Mumford
stack with a coarse moduli space $\Higgs_{|B-\Delta}$.  If we view
$\gHiggs_{|B-\Delta}$ as a group stack, then the group of \linebreak 
connected
components, as well as the connected components of each fiber of
\linebreak $\bh :
\gHiggs_{|B-\Delta} \to (B-\Delta)$ are canonically isomorphic to
$\pi_{1}(G)$.
\item[{\bf (ii)}] $\gHiggs_{|B-\Delta}$ is a banded $Z(G)$-gerbe over
$\Higgs_{|B-\Delta}$ which is locally trivial over $B - \Delta$. In
particular the restriction of $\gHiggs_{|B-\Delta}$ to a Hitchin fiber
is a trivial gerbe.
\item[{\bf (iii)}] The gerbe $\gHiggs_{|B-\Delta} \to
  \Higgs_{|B-\Delta}$  measures the obstruction to lifting the
  universal $G_{\op{ad}}$-Higgs bundle to a universal $G$-Higgs bundle.
\end{itemize}
\end{lem}
{\bf Proof.} {\bf (i)} It is well known
 \cite{simpson-moduli1,simpson-moduli2} that the stack $\gHiggs$ of
 $G$-Higgs bundles is an Artin algebraic stack with an affine diagonal
 which is locally of finite type. The substack of semistable Higgs
 bundles is of finite type and has a quasi-projective moduli
 space. The statement about the connected components is now automatic,
 since $\pi_{0}(\gHiggs_{|B - \Delta})    = \pi_{0}(\Higgs_{|B -
 \Delta})$ and the connected components of  $\Higgs_{|B -
 \Delta}$ and the corresponding Hitchin fibers were already described
 in Claim~\ref{claim:components}{\bf (ii)}.

It is also known
 \cite{simpson-moduli2,beilinson-drinfeld-langlands} that the stack
 $\gHiggs$ is a local complete intersection which is smooth at all
 points with finite stabilizers. Therefore it suffices to show that
 $\gHiggs_{|B-\Delta}$ parametrizes Higgs bundles with minimal
 automorphism group. In \cite[Theorem~III.2]{faltings}  Faltings
 showed that $B$ contains a Zariski open and dense subset $B^{o}
 \subset B$, such that $\gHiggs_{|B^{o}}$ parametrizes only stable
 Higgs bundles with automorphism group $Z(G)$. We will give a direct
 argument for this over $B-\Delta$, i.e. we will show that $B^{o}
 \supset B-\Delta$. 

Fix a point $b \in (B - \Delta)$ and let $p : \widetilde{C} \to C$ be
the corresponding cameral cover. We must show that every object in the
groupoid $\gHiggs_{\widetilde{C}} := \bh^{-1}(b)$ has automorphism
group $Z(G)$. As explained in Lemma~\ref{lem:section},
$\gHiggs_{\widetilde{C}}$ is the groupoid of
$\mathcal{T}_{\widetilde{C}} := \mathcal{T}_{|\{b \}\times
  C}$-torsors, and hence the automorphism 
group of any object in $\gHiggs_{\widetilde{C}}$ is isomorphic to the
cohomology group $H^{0}(C,\mathcal{T}_{\widetilde{C}})$.
By the argument we used in the proof of Theorem~\ref{thm:duality}{\bf
(2)} we have isomorhisms $H^{0}(C,\mathcal{T}_{\widetilde{C}}) \cong
H^{0}(C,\mathcal{T}_{\widetilde{C},\mathbb{R}})$,
$H^{0}(C,\mathcal{T}^{o}_{\widetilde{C}}) \cong
H^{0}(C,\mathcal{T}_{{\widetilde{C}},\mathbb{R}}^{o})$, 
$H^{0}(C,\overline{\mathcal{T}}_{\widetilde{C}}) \cong
H^{0}(\widetilde{C},\overline{\mathcal{T}}_{{\widetilde{C}},\mathbb{R}})$.
Thus it 
suffices to compute
$H^{0}(C,\mathcal{T}_{{\widetilde{C}},\mathbb{R}})$. We start by
calculating the global sections of $\mathcal{T}_{\mathbb{R}}^{o} =
(\jmath_{*}A)\otimes S^{1}$. As in the proof of
Claim~\ref{claim:components}{\bf (ii)} we get a short exact sequence
of sheaves on $\widetilde{C}$:
\[
0 \to \jmath_{*}A \to (\jmath_{*}A)\otimes \mathbb{R} \to
(\jmath_{*}A)\otimes S^{1} \to 0.
\]
Passing to cohomology, and taking into account that
$H^{0}(C,\jmath_{*}A) = 0$ and
$H^{0}(C,(\jmath_{*}A)\otimes \mathbb{R}) =
H^{0}(C,\jmath_{*}A)\otimes \mathbb{R} = 0$, we get that 
\[
\begin{split}
H^{0}(C,\mathcal{T}_{\widetilde{C},\mathbb{R}}^{o}) & = \ker\left[
  H^{1}(C,\jmath_{*}A) \to
  H^{1}(C,\jmath_{*}A\otimes \mathbb{R})\right] \\
& = H^{1}(C,\jmath_{*}A)_{\op{tor}}.
\end{split}
\]
The latter group can be calculated explicitly from
Corollary~\ref{cor:torsion}. In the notation of
Corollary~\ref{cor:torsion}, let $N$ denote the saturation of
$(1-\rho_{1}, \ldots, 1-\rho_{b})\Lambda$ inside $\oplus_{i =1}^{b}
\mathbb{Z}\varepsilon_{i}\alpha_{i}^{\vee}$. Then  
\[
\begin{split}
H^{1}(C,\jmath_{*}A)_{\op{tor}} & = N/(1-\rho_{1}, \ldots,
1-\rho_{b})\Lambda \\
& = \left\{ \xi \in \mathfrak{t} \, \left| \,
(\xi,\alpha) \in \varepsilon_{\alpha,G}\mathbb{Z} \text{ for every root }
\alpha \right. \right\}/\Lambda \\
& = \begin{cases} 
Z(G) & \text{ if } G \neq \op{Sp}(r) \\
0 & \text{ if }  G = \op{Sp}(r).
\end{cases}
\end{split}
\]
From Lemma~\ref{lem:diagram} we know that as long as $G \neq
\op{Sp}(r)$ we have 
$\mathcal{T}_{\widetilde{C},\mathbb{R}}^{o} =
\mathcal{T}_{\widetilde{C},\mathbb{R}}$.  This proves
our claim for $G \neq
\op{Sp}(r)$. For $G = \op{Sp}(r)$, Lemma~\ref{lem:diagram} gives  a
short exact sequence 
\[
0 \to \mathcal{T}_{\widetilde{C},\mathbb{R}}^{o} \to
\mathcal{T}_{\widetilde{C},\mathbb{R}} \to 
\oplus_{i = 1}^{b} \mathbb{Z}/\varepsilon_{i} \to 0,
\]
and after passing to cohomology we get 
\[
\begin{split}
H^{0}(C,\mathcal{T}_{\widetilde{C},\mathbb{R}}) & = \ker 
\left[ \oplus_{i = 1}^{b} \mathbb{Z}/\varepsilon_{i} \to
  H^{1}(C, \mathcal{T}_{\mathbb{R}}^{o}) \right] \\
& = Z(\op{Sp}(r)) \cong \mathbb{Z}/2,
\end{split} 
\]
where $Z(\op{Sp}(r))\cong \mathbb{Z}/2$ maps diagonally in $\oplus_{i
= 1}^{b} \mathbb{Z}/\varepsilon_{i}$.  This proves our assertion about
the automorphisms of objects in $\gHiggs_{\widetilde{C}}$ and finishes
the proof of {\bf (i)}.

\

\noindent
{\bf (ii)} As we saw above, the $\mathcal{T}$ torsors in $\gHiggs_{|B
  - \Delta}$ all have automorphism groups isomorphic to $Z(G)$, and so
$\gHiggs_{|B - \Delta}$ is a $Z(G)$-gerbe on $\Higgs_{|B
  - \Delta}$. In particular  $R^{0}\pi_{*}\mathcal{T}$ is a local
system on $B - \Delta$ with fiber $Z(G)$. However $Z(G) = T^{W}
\subset T$  and so, by the definition of $\mathcal{T}$ we have a
canonical inclusion of the constant sheaf $Z(G)$ into
$\mathcal{T}$. Thus, every element in  $Z(G)$ gives rise to a global
section of $\mathcal{T}$ on $B\times C$, and hence to a global section
of $R^{0}\pi_{*}\mathcal{T}$ on $B-\Delta$. This shows that
$R^{0}\pi_{*}\mathcal{T}$ is the constant sheaf and so $\gHiggs_{|B
  - \Delta}$ is banded as a gerbe over $\Higgs_{|B
  - \Delta}$. Finally, note that locally over $B - \Delta$, the
universal cameral cover $\widetilde{\mycal{C}}$ admits a section. The stack of
$\mathcal{T}$-torsors which are framed along such a section is
isomorhic to the space $\Higgs$, which shows that the gerbe $\gHiggs_{|B
  - \Delta}$ is locally trivial over $(B
  - \Delta)$. 

\

\noindent
{\bf (iii)} This is completely analogous to the $\op{P}\op{SL}(r)$
argument in \cite{hausel-thaddeus}. Let $G_{\op{ad}}$ be the adjoint
form of $G$. From part {\bf (ii)} it follows that the stack of
$G_{\op{ad}}$-Higgs bundles that have cameral cover in $B - \Delta$ is
actually a space, i.e.  $\gHiggs_{G_{\op{ad}}|(B - \Delta)} =
\Higgs_{G_{\op{ad}}|(B - \Delta)}$. Since the stack always has a
universal bundle, we have a universal
$G_{\op{ad}}$-Higgs bundle $(\boldsymbol{V},\boldsymbol{\varphi})$ on
$\Higgs_{G_{\op{ad}}|(B - \Delta)} \times C$. The natural map $G \to
G_{\op{ad}}$ induces a morphism of spaces $q : \Higgs_{|B-\Delta} \to
\Higgs_{G_{\op{ad}}|(B - \Delta)}$ and we can consider the pullback
$G_{\op{ad}}$-Higgs bundle
$q^{*}(\boldsymbol{V},\boldsymbol{\varphi})$ on
$\Higgs_{|B-\Delta}$. Since $\ker[G \to G_{\op{ad}}] = Z(G)$, it follows
that the obstruction to lifting
$q^{*}(\boldsymbol{V},\boldsymbol{\varphi})$ to a $G$-Higgs bundle is 
simply the obstruction to the existence of an universal $G$-Higgs
bundle on $\Higgs_{|(B - \Delta)}\times C$, i.e it is the gerbe 
$\gHiggs_{|(B- \Delta)}\times C$. In particular, restricting to
$\Higgs_{|(B - \Delta)}\times \{ \op{pt} \}$ we get the statement {\bf
  (iii)}. \ \hfill $\Box$

\

\subsection{Global duality} \label{ss:global}

We are now ready to state the main result of this section. For any
commutative group stack $h : \boldsymbol{\mycal{X}} \to S$ over a scheme $S$
(for us $S$ will always be $B-\Delta$)  with
zero section $\mathfrak{v} : S \to \boldsymbol{\mycal{X}}$,
and group law $\ba :
\boldsymbol{\mycal{X}}\times_{S}\boldsymbol{\mycal{X}} \to
\boldsymbol{\mycal{X}}$,   we
define the dual commutative group stack as the stack of homomorphisms
of commutative group stacks  from $\mycal{X}$ to $B\mathbb{G}_{m}$:
\[
\boldsymbol{\mycal{X}}^{D} :=
\underline{\op{Hom}}_{\op{grp-stack}}(\boldsymbol{\mycal{X}}, 
\mathcal{O}^{\times}[1])w.
\]
Geometrically $\boldsymbol{\mycal{X}}^{D}$ is the stack of group
extensions of $\boldsymbol{\mycal{X}}$ by $\mathcal{O}^{\times}$, or
equivalently, the stack parametrizing triples
$(\mathfrak{L},\mathfrak{m},\mathfrak{f})$, where $\mathfrak{L}$ is a
line bundle 
on $\boldsymbol{\mycal{X}}$, and $(\mathfrak{m},\mathfrak{f})$ is a
'theorem of the 
square structure' on $\mathfrak{L}$
\cite{mumford-abelian,breen-cube}. Concretely the square structure
cosists of an
isomorphism of line bundles $\mathfrak{m} :
p_{1}^{*}\mathfrak{L}\otimes p_{2}^{*}\mathfrak{L} \to \ba^{*}
\mathfrak{L}$ on
$\boldsymbol{\mycal{X}}\times_{S}\boldsymbol{\mycal{X}}$ and a framing 
$\mathfrak{f} : \mathfrak{v}^{*}\mathfrak{L} \stackrel{\cong}{\to} 
\mathcal{O}_{S}$
which satisfy the normalizations $(\mathfrak{v}\circ h\times
\op{id})^{*}\mathfrak{m} = 
h^{*}e\otimes \op{id}_{\mathfrak{L}}$  and $(\op{id}\times
\mathfrak{v}\circ h)^{*}\mathfrak{m}  = 
\op{id}_{\mathfrak{L}}\otimes  h^{*}e$ on $\mycal{X}$ and the cocycle
condition $(p_{1+2}\times p_{3})^{*}\mathfrak{m}\circ
(p_{12}^{*}\mathfrak{m}\otimes \op{id}) = (p_{1}\times
p_{2+3})^{*}\mathfrak{m}\circ 
(\op{id}\otimes p_{23}^{*}\mathfrak{m})$ on the triple product
$\boldsymbol{\mycal{X}} 
\times_{S} \boldsymbol{\mycal{X}} \times_{S}
\boldsymbol{\mycal{X}}$. 

Note that if $\boldsymbol{\mycal{X}}$ is an abelian scheme over $S$,
then  $\boldsymbol{\mycal{X}}^{D}$ is the usual dual abelian scheme.
Also for well behaved group stacks $\boldsymbol{\mycal{X}}$ the duality
operation $(\bullet)^{D}$ converts disconnectedness into gerbiness and
vice versa \cite{dp}.  In particular if the relative group
$\pi_{0}(\boldsymbol{\mycal{X}}/S)$ of connected components
is a finite flat group scheme over $S$, then
$\boldsymbol{\mycal{X}}^{D}$ will be a banded
$\pi_{0}(\boldsymbol{\mycal{X}}/S)^{\wedge}$-gerbe. Here as usual
$\pi_{0}(\boldsymbol{\mycal{X}}/S)^{\wedge}$ denotes the Pontryagin
dual group
$\underline{\op{Hom}}_{\mathbb{Z}}(\pi_{0}(\boldsymbol{\mycal{X}}/S),
\mathcal{O}^{\times}_{S})$.

\

\noindent
With this notation we now have:

\begin{thm} \label{thm:gerbes} Let $\gHiggs$ be the stack of
  $G$ Higgs bundles on a curve $C$ and let $\lan{\gHiggs}$
  be the stack of $\lan{G}$ Higgs bundles on $C$. Use
  the isomorphism $\bl : B \to \lan{B}$ from {\bf
  Theorem~\ref{thm:duality}(1)} to identify $B-\Delta$ with
  $\lan{B}-\lan{\Delta}$. Under this identification one has a canonical
  isomorphism
\begin{equation} \label{eq:gerbeD}
\gerbeiso : \gHiggs_{|B-\Delta} \stackrel{\cong}{\longrightarrow}
(\lan{\gHiggs}_{|B-\Delta})^{D}  
\end{equation}
of commutative group stacks over $B-\Delta$.  The isomorphism
$\gerbeiso$ intertwines the action of the translation operators
$\trans^{\lambda,\tilde{x}}$ on $\gHiggs_{|B-\Delta}$ with the action
of the tensorization operators $\tens^{\lambda,\tilde{x}}$ on
$(\lan{\gHiggs}_{|B-\Delta})^{D}$. 
\end{thm}

\

\begin{rem}
Here we only consider the parts of the stacks of Higgs
bundles sitting over $B-\Delta$. To simplify notation throughout the
proof we will write $\gHiggs$ and $\lan{\gHiggs}$ instead of
$\gHiggs_{|B-\Delta}$ and $\lan{\gHiggs}_{|B-\Delta}$. 
\end{rem}

\

\begin{rem} \label{rem:isogenies} It will be clear form the proof of
Theorem~\ref{thm:gerbes} that the canonical isomorphism  we construct
is functorial with respect to isogenies $G \to G'$ between
simple groups. 
\end{rem}

\

\subsection{Strategy of the proof} \label{ssec:strategy}
Before we present any details of the proof we outline the strategy
that we will folow.  
Fix a cameral cover $\widetilde{C}$ corresponding to a point in $B -
\Delta$. We want to extend our isomorphism 
\begin{equation} \label{eq:paiso}
\xymatrix@1@M+1pc{
\paiso_{\widetilde{C}} : \hspace{-5pc} &  
P_{\widetilde{C}} \ar[r]^-{\cong} & 
\widehat{\lan{P}}_{\widetilde{C}} = \lan{P}_{\widetilde{C}}^{D}
}
\end{equation}
from {\bf  Theorem~\ref{thm:duality}(2)} to  
 a natural isomorphism of group stacks 
\begin{equation} \label{eq:global.iso}
\gHiggs_{\widetilde{C}}  \stackrel{\cong}{\to}
\left(\lan{\gHiggs}_{\widetilde{C}}\right)^{D}.
\end{equation}
Naturality will imply in particular that the isomorphisms
\eqref{eq:global.iso} of individual Hitchin fibers will globalize to
an isomorphism
\[
\gHiggs  \stackrel{\cong}{\to}
\left(\lan{\gHiggs}\right)^{D}
\]
over $B - \Delta$.

We will construct \eqref{eq:global.iso} in several steps.

\begin{itemize}
\item[(i)] We will construct actions of the algebra of abelianized
  Hecke operators $\trans^{\lambda,\tilde{x}}$ on
  $\gHiggs_{\widetilde{C}}$ and of the abelianized tensorization
  operators $\tens^{\lambda,\tilde{x}}$ on
  $\left(\lan{\gHiggs}_{\widetilde{C}}\right)^{D}$.  Both types of
  operators are labeled by a point $\tilde{x} \in \widetilde{C}$ and a
  cocharacter $\lambda \in \Lambda = \cchr{T}$. We will use the
  abelianized Hecke operators to construct a Higgs version of
  the Abel-Jacobi map for cameral covers.
\item[(ii)] The stacks in \eqref{eq:global.iso} are (possibly
  disconnected) commutative group stacks. In
  Section~\ref{ss:universal} we saw that the groups of connected
  components of both $\gHiggs_{\widetilde{C}}$ and
  $\left(\lan{\gHiggs}_{\widetilde{C}}\right)^{D}$ are naturally
  identified with $\pi_{1}(G) = \Lambda/\crts_{\mathfrak{g}}$. The
  actions of $\trans^{\lambda,\tilde{x}}$ and
  $\tens^{\lambda,\tilde{x}}$ on $\pi_{1}(G)$ are induced from the
  translation action of $\lambda$ on $\Lambda$ and so the  abelianized
  Hecke and tensorization operators permute
  transitively the components of $\gHiggs_{\widetilde{C}}$ and
  $\left(\lan{\gHiggs}_{\widetilde{C}}\right)^{D}$.

Thus to construct the map \eqref{eq:global.iso} it suffices
to construct an isomorphism 
\begin{equation} \label{eq:connected.stacks.global.iso} 
\gHiggs_{0,\widetilde{C}} \stackrel{\cong}{\to}
\left(\lan{\Higgs}_{\widetilde{C}}\right)^{D} 
\end{equation}
of the connected components of the identity of our stacks so that:
\begin{itemize}
\item the isomorphism \eqref{eq:connected.stacks.global.iso}  extends
the isomorphism of abelian varieties \eqref{eq:paiso}.
\item the isomorphism \eqref{eq:connected.stacks.global.iso}
  intertwines the action of the Hecke and tensorization operators
  that preserve the connected components of the identity. That is, for
  every $\tilde{x} \in
  \widetilde{C}$, $\lambda \in \crts_{\mathfrak{g}}$,
  \eqref{eq:connected.stacks.global.iso}   the
  action of $\trans^{\lambda,\tilde{x}}$ on
  $\gHiggs_{0,\widetilde{C}}$ with the action of
  $\tens^{\lambda,\tilde{x}}$ on 
  $\left(\lan{\Higgs}_{\widetilde{C}}\right)^{D}$;
\end{itemize}
Indeed an isomorphism  \eqref{eq:connected.stacks.global.iso} with
these properties  will  extend uniquely to an isomorphism
\eqref{eq:global.iso} by Hecke equivariance.
\item[(iii)] The stacks $\gHiggs_{0,\widetilde{C}}$ and
$\left(\lan{\Higgs}_{\widetilde{C}}\right)^{D}$ are connected
commutative group stacks with moduli spaces $P_{\widetilde{C}}$ and
$\lan{P}_{\widetilde{C}}^{D}$ respectively.  We construct
\eqref{eq:connected.stacks.global.iso} by exhibiting explicit groupoid
scheme presentations for these stacks and showing that the isomorphism
$\paiso_{\widetilde{C}}$ extends to an isomorphism of the
presentations. We also check that the abelianized Hecke and
tensorization operators are induced from automorphisms of the
presentations, and that the lift of $\paiso_{\widetilde{C}}$
intertwines these automorphisms. To achieve this we study the moduli
spaces of neutralizations of $\gHiggs_{0}$ and
$\left(\lan{\Higgs}\right)^{D}$ along each Hitchin fiber.
\end{itemize}

\subsection{Abelianized Hecke operators and the Abel-Jacobi map}
\label{ssec:AJ}

We will use the
abelianization of Higgs bundles described in
\cite[Section~6]{ron-dennis} to define Hecke operators on the moduli
stack of Higgs bundles. 

\

\medskip

\punkt For every $\lambda \in \Lambda := 
\cchr_{G} = \op{Hom}(\mathbb{C}^{\times},T)$
and $\tilde{x} \in \widetilde{C} \subset \widetilde{\mycal{C}}$ we will
construct a canonical Hecke automorphism $\trans^{\lambda,\tilde{x}} :
\gHiggs_{\widetilde{C}} \to \gHiggs_{\widetilde{C}} $. These
automorphisms can be combined into a single map of stacks
\begin{equation} \label{eq:heckemap}
\trans : \gHiggs\times_{(B-\Delta)} \widetilde{\mycal{C}}\times \Lambda
\to \gHiggs,
\end{equation}
which we will call {\em the abelianized Hecke correspondence}. This
map induces a natural map on coarse moduli spaces which we will denote
again by $\trans$. 

Let $(V,\varphi)$ be a $K_{C}$-valued $G$-Higgs bundle with cameral
cover $p: \widetilde{C} \to C$, and let $\tilde{x} \in \widetilde{C}$,
$\lambda \in \Lambda$. Informally $\trans$ takes the data
$((V,\varphi),\tilde{x},\lambda)$ to a new Higgs bundle
$(V',\varphi')$ having the same cameral cover $p: \widetilde{C} \to
C$, an underlying $G$-bundle $V'$ which is a modification of $V$ at
$p(\tilde{x})$ in the direction $\lambda$, and a Higgs field
$\varphi'$ which agrees with the original $\varphi$ on $C - \{
p(\tilde{x}) \}$.

More formally, \cite[Theorem~6.4]{ron-dennis} establishes an
equivalence between the groupoid of $G$-Higgs bundles on $C$ with
cameral cover $\widetilde{C}$ and the groupoid of $G$-spectral data on
$\widetilde{C}$. Spectral data are collections $(\mycal{L},\bi,\bb)$
consisting of a ramification twisted $W$-invariant $T$-bundle
$\mycal{L} = \mycal{L}_{(V,\varphi)}$ on $\widetilde{C}$ equipped with
additional twisting and framing structures $\bi$ and $\bb$ satisfying
compatibility conditions. For the convenience of the reader we review
the precise description of $G$-spectral data and the abelianized Hecke
correspondences $\trans^{\lambda,\tilde{x}}$ in
Appendix~\ref{appendix:spectral}. (Note that the Appendix discusses
Heckes for $\lan{G}$ while here, dually, we need Heckes for $G$.)

As explained in Appendix~\ref{appendix:spectral} the abelianized Hecke
correspondence $\trans^{\lambda,\tilde{x}}$ is an automorphism of the
stack $\gHiggs_{\widetilde{C}}$ defined by tensoring with a particular
$\mathcal{T}_{\widetilde{C}}$ torsor $\mathcal{S}^{\lambda,\tilde{x}}$
on $C$. In the spectral picture $\trans^{\lambda,\tilde{x}}$ modifies
the bundle $\mycal{L}$ and the framing $\bb$ in the spectral
datum. Specifically (see Appendix~\ref{appendix:spectral})
$\trans^{\lambda,\tilde{x}}(\mycal{L},\bi,\bb) = (\mycal{L}\otimes
S^{\lambda,\tilde{x}},\bi,\bb\otimes \sss^{\lambda,\tilde{x}})$, where
$S^{\lambda,\tilde{x}}$ is the $W$-equivariant $T$-bundle:
\begin{equation} \label{eq:Sxlambda}
S^{\lambda,\tilde{x}} := \bigotimes_{w \in W} (w\lambda)\left(
\mathcal{O}_{\widetilde{C}}(w\tilde{x}) \right).
\end{equation}
The necessary compatibilities are automatic as long as
$\tilde{x} \in \widetilde{C}$ is not a ramification point of $p :
\widetilde{C} \to C$. Since these compatibilities are closed
conditions, they hold on all of $\widetilde{C}$ (or
$\widetilde{\mycal{C}}$).  

\

\medskip

\punkt \ \label{sss:hecke.action} 
By construction the abelianized Hecke operator
$\trans^{\lambda,\tilde{x}} : \gHiggs_{\widetilde{C}}\to
\gHiggs_{\widetilde{C}}$ shifts the components of
$\gHiggs_{\widetilde{C}}$ by the image of $\lambda$ in $\pi_{1}(G) =
\Lambda/\crts_{\mathfrak{g}}$. In particular, any component of
$\gHiggs_{\widetilde{C}}$ can be reached from the neutral one by
applying a suitable Hecke operator. As a side note, observe that the
$T$-bundles $S^{\lambda,\tilde{x}}$ all have the same topological type
and so the tensorization by these bundles preserves components of
$\sBun_{T}$:

\begin{lem} \label{lem:trivial.type} For every $\tilde{x} \in
  \widetilde{C}$ and every $\lambda \in \Lambda$ the $T$-bundle
  $S^{\lambda,\tilde{x}}$ is topologically trivial.
\end{lem}
{\bf Proof.} The topological type of a $T$-bundle on $\widetilde{C}$
is classified by a characteristic class in the lattice
$H^{2}(\widetilde{C},\pi_{1}(T)) = H^{2}(\widetilde{C},\Lambda) \cong
\Lambda$.  From \eqref{eq:Sxlambda} it follows that the charactersitic
class of $S^{\lambda,\tilde{x}}$ is a $W$-invariant element in
$\Lambda$. Since $G$ is assumed semisimple we have that the group of
$W$-invariants $\Lambda^{W}$ is trivial, which proves the lemma. \
\hfill $\Box$

\

\medskip

\punkt \ \label{sss:dual.hecke.action}
Let now $\mathfrak{v} : (B-\Delta) \to \gHiggs$ be a Hitchin
section. By applying the Hecke map $\trans$ to $\mathfrak{v}$ we get
a Higgs bundle version of the Abel-Jacobi map:
\[
\xymatrix@R-1.5pc{
\widetilde{\mycal{C}}\times \Lambda \ar[r]^-{\shaj{G}} & \gHiggs \\
(\tilde{x},\lambda) \ar@{|->}[r] &
\trans^{\lambda,\tilde{x}}(\mathfrak{v}(\bh(\tilde{x}))). 
}  
\]  
Composing this map with the map $\gHiggs \to \Higgs$ to the coarse
moduli space we get a moduli space version of the Abel-Jacobi map:
\[
\xymatrix@1{
\widetilde{\mycal{C}}\times \Lambda \ar[r]^-{\haj{G}} & \Higgs.
}  
\]  
For the construction of the groupoid presentations of
$(\lan{\Higgs})^{D}$ and $\gHiggs$ we will also need the Langlands
dual version of these maps. Fix a Hitchin section $\lan{\mathfrak{v}}
: (B-\Delta) \to \lan{\gHiggs}$ and let
\[
\xymatrix@R-1pc@C-1pc{
& \lan{\gHiggs} \ar[dd] \\
\widetilde{\mycal{C}}\times \lan{\Lambda} \ar[ru]^-{\shaj{\lan{G}}}
\ar[rd]_-{\haj{\lan{G}}}
 & \\
& \lan{\Higgs}.
}  
\] 
be the corresponding Abel-Jacobi maps. The restrictions of these maps
to slices of the form $\{\tilde{x} \} \times \lan{\Lambda}$ yields
group homomorphisms 
\[
\xymatrix@R-1pc@C+1pc{
& \lan{\gHiggs} \ar[dd] \\
\lan{\Lambda} \ar[ru]^-{\shaj{\lan{G}}(\tilde{x},-)}
\ar[rd]_-{\haj{\lan{G}}(\tilde{x},-)}
 & \\
& \lan{\Higgs}.
}  
\] 
For future reference we introduce the shortcut notation
$\hm_{\tilde{x}} : \Lambda \to \Higgs_{\widetilde{C}}$ and \linebreak 
$\lan{\hm}_{\tilde{x}} :
\lan{\Lambda} \to \lan{\Higgs}_{\widetilde{c}}$ for the homomorphisms
$\haj{G}(\tilde{x},-)$ and $\haj{\lan{G}}(\tilde{x},-)$.  We will also
write $\hm'_{\tilde{x}} : \crts_{\mathfrak{g}} \to P_{\widetilde{C}}$
and 
$\lan{\hm}'_{\tilde{x}} : \crts_{\lan{\mathfrak{g}}} \to
\lan{P}_{\widetilde{C}}$ for the restrictions of $\haj{G}(\tilde{x},-)$
and $\haj{\lan{G}}(\tilde{x},-)$ to $\crts_{\mathfrak{g}}$ and
$\crts_{\lan{\mathfrak{g}}}$ respectively.

\subsection{Abelianized tensorization operators} \label{ssec:dualHecke}

\

\noindent
Since $(\lan{\gHiggs})^{D}$ is the group stack of line bundles with
square structures on $\lan{\gHiggs}$ it follows that every line bundle
$\mathfrak{L}$ with square structure on $\lan{\gHiggs}$ will give rise
to a tensorization automorphism $(-)\otimes \mathfrak{L}$ of
$(\lan{\gHiggs})^{D}$.

The abelianized tensorization operator $\tens^{\lambda,\tilde{x}} :
\left(\lan{\gHiggs}_{\widetilde{C}}\right)^{D} \to
\left(\lan{\gHiggs}_{\widetilde{C}}\right)^{D}$ will be given
concretely as the tensorization  
$\tens^{\lambda,\tilde{x}}(-) := (-)\otimes
\mathfrak{L}^{\lambda,\tilde{x}}$ with a particular line bundle with
square structure $\mathfrak{L}^{\lambda,\tilde{x}}$ which we construct
next.

\noindent
To ensure the additivity in $\lambda$ of the operators
$\tens^{\lambda,\tilde{x}}$ we will construct the line bundle
$\mathfrak{L}^{\lambda,\tilde{x}}$ with square structure on
$\lan{\gHiggs}_{\widetilde{C}}$ so that the line bundle
$(\shaj{\lan{G}})^{*}\mathfrak{L}^{\lambda,\tilde{x}}$ on
$\widetilde{C}\times \lan{\Lambda}$, interpreted as a $T$-bundle on
$\widetilde{C}$, will equal $S^{\lambda,\tilde{x}}$.

Note that by Lemma~\ref{lem:section}, the choice of the Hitchin
section $\lan{\mathfrak{v}}$ identifies $\lan{\gHiggs}$ with the group stack
over $B - \Delta$ associated with
$R^{\bullet}\pi_{*}(\lan{\mathcal{T}})[1]$, and similarly identifies
$\lan{\Higgs}$ with the commutative group scheme on $B-\Delta$
representing the sheaf $R^{1}\pi_{*}(\lan{\mathcal{T}})$. The sheaf
inclusion $\lan{\mathcal{T}} \subset \lan{\overline{\mathcal{T}}}$
induces natural maps of complexes of abelian sheaves on $B-\Delta$:
\begin{equation} \label{eq:composite}
\xymatrix@R-1pc@C-1pc{ R^{\bullet}\pi_{*}(\lan{\mathcal{T}}) \ar[r]
   \ar@/_5pc/[rddddd]_-{\lan{\sab}} &
R^{\bullet}\pi_{*}(\lan{\overline{\mathcal{T}}}) \ar@{=}[d] \\ 
&
 R^{\bullet}\pi_{*}\left(\left[ p_{*}\left( \lan{\Lambda} \otimes
\mathcal{O}^{\times}_{\widetilde{\mycal{C}}}\right) \right]^{W}\right)
\ar[d] \\ & R^{\bullet}\pi_{*}\left( p_{*}\left( \lan{\Lambda} \otimes
\mathcal{O}^{\times}_{\widetilde{\mycal{C}}}\right) \right)
\ar[d]^-{\text{Leray}} \\ & R^{\bullet}\tilde{\pi}_{*}\left( \lan{\Lambda}
\otimes \mathcal{O}^{\times}_{\widetilde{\mycal{C}}}\right) \ar@{=}[d]
\\ &
\left(R^{\bullet}\tilde{\pi}_{*}
\mathcal{O}^{\times}_{\widetilde{\mycal{C}}}\right) 
\otimes \lan{\Lambda} 
\ar@{=}[d]
\\
& \gPic(\widetilde{\mycal{C}}/(B-\Delta))\otimes \lan{\Lambda}.
}
\end{equation}
The map labeled ``Leray'' comes  from the Leray spectral sequence for
$p : \widetilde{\mycal{C}} \to (B\times C)$ and is an isomorphism
because $p$ is finite. 

The composite map $\lan{\sab} :
R^{\bullet}\pi_{*}(\lan{\mathcal{T}}) \to
\gPic(\widetilde{\mycal{C}}/B-\Delta)\otimes \lan{\Lambda}$ induces a 
morphism of stacks
\[
\lan{\sab} : \lan{\gHiggs} \to
\gPic(\widetilde{\mycal{C}}/(B-\Delta))\otimes \lan{\Lambda} =
\sBun_{\, \lan{T}},
\]
where $\sBun_{\, \lan{T}}$ denotes the stack parametrizing
$\lan{T}$-bundles along the 
fibers of \linebreak $\tilde{\pi} : \widetilde{\mycal{C}} \to (B-\Delta)$. 

Similarly, if in diagram \eqref{eq:composite} we replace
$R^{\bullet}\pi_{*}$ with $R^{1}\pi_{*}$ we get a composite map \linebreak
$\lan{\ab} : R^{1}\pi_{*}(\lan{\mathcal{T}}) \to
\Pic(\widetilde{\mycal{C}}/B-\Delta)\otimes \lan{\Lambda}$ which
induces a morphism of spaces
\[
\lan{\ab} : \lan{\Higgs} \to
\Pic(\widetilde{\mycal{C}}/(B-\Delta))\otimes \lan{\Lambda} =
\Bun_{\, \lan{T}}.
\]
Combining these maps with the Abel-Jacobi maps for Higgs bundles we
get
commutative diagrams
\begin{equation} \label{eq:aj}
\xymatrix@C-1pc{
\widetilde{\mycal{C}}\times \lan{\Lambda} \ar[r]^-{\shaj{\lan{G}}}
\ar[d]_-{\saj\times \op{id}} & \lan{\gHiggs} \ar[dd]^-{\lan{\sab}} \\
\gPic(\widetilde{\mycal{C}}/(B-\Delta))\otimes \lan{\Lambda}
\ar@{=}[d] & \\
\sBun_{\, \lan{T}} \ar[r]_-{\sum_{w \in W}\,\dia_{w}}  & \sBun_{\, \lan{T}}
} \qquad \text{and} \qquad 
\xymatrix@C-1pc{
\widetilde{\mycal{C}}\times \lan{\Lambda} \ar[r]^-{\haj{\lan{G}}}
\ar[d]_-{\aj\times \op{id}} & \lan{\Higgs} \ar[dd]^-{\lan{\ab}} \\
\Pic(\widetilde{\mycal{C}}/(B-\Delta))\otimes \lan{\Lambda}
\ar@{=}[d] & \\
\Bun_{\, \lan{T}} \ar[r]_-{\sum_{w \in W}\, \dia_{w}}  & \Bun_{\, \lan{T}}
}
\end{equation}
Here as usual $\dia_{w}(\bullet) = w^{*}((\bullet)\times_{w}\, \lan{T})$ is
the diagonal action of $w$ on $\lan{T}$-bundles and \linebreak $\saj :
\widetilde{\mycal{C}} \to \gPic(\widetilde{\mycal{C}}/(B-\Delta))$
denotes the classical Abel-Jacobi map, sending a point $\tilde{x} \in
\widetilde{C}$ to the line bundle
$\mathcal{O}_{\widetilde{C}}(\tilde{x})$ of degree one on
$\widetilde{C}$.

We are now ready to construct the line bundles
$\mathfrak{L}^{\lambda,\tilde{x}}$ on $\lan{\gHiggs}$ given by
$\tilde{x} \in \widetilde{C}$, $\lambda \in \Lambda$. For this we will
use the well known fact \cite{laumon} that the Picard gerbe on
any smooth family of curves is self-dual. Note that this is precisely
our Theorem~\ref{thm:gerbes} in the abelian case $G =
\mathbb{C}^{\times}$.  More precisely, for any smooth compact complex
curve $\Sigma$, there exists a Poincare sheaf on $\gPic(\Sigma)\times
\gPic(\Sigma)$ which induces a canonical isomorphism
$(\gPic(\Sigma))^{D} = \gPic(\Sigma)$. In fact this isomorphism is
induced by the classical Abel-Jacobi map $\aj : \Sigma \to
\gPic(\Sigma)$:
\[
\aj^{*} : (\gPic(\Sigma))^{D}  \to \gPic(\Sigma).
\]
We apply this to the curve $\widetilde{C}$ and tensor with
$\lan{\Lambda}$ to get an induced isomorphism 
\[
(\aj\times \op{id})^{*} : (\gPic(\widetilde{C})\otimes
\lan{\Lambda})^{D} \; \stackrel{\cong}{\to} \; \sBun_{T}(\widetilde{C}).
\]
In particular, for every $\tilde{x} \in \widetilde{C}$ and any $\lambda \in
\Lambda$ we can find a canonical line bundle
$\mathbb{L}^{\lambda,\tilde{x}}$ on $\gPic(\widetilde{C})\times
\lan{\Lambda}$ such that $(\aj\times
\op{id})^{*}\mathbb{L}^{\lambda,\tilde{x}}  =
\lambda(\mathcal{O}_{\widetilde{C}}(\tilde{x}))$. To finish the construction
we set 
\[
\mathfrak{L}^{\lambda,\tilde{x}} := \sab^{*}\mathbb{L}^{\lambda,\tilde{x}}
\]
and invoke the commutative diagram \eqref{eq:aj} to get the desired identity
\begin{equation} \label{eq:compatible.pi0}
(\shaj{\lan{G}})^{*}\mathfrak{L}^{\lambda,\tilde{x}} =
(\shaj{\lan{G}})^{*}(\lan{\sab})^{*}\mathbb{L}^{\lambda,\tilde{x}} =
\left(\sum_{w \in W} \dia_{w}\right)^{*}(\aj\times
\op{id})^{*}\mathbb{L}^{\lambda,\tilde{x}} = 
S^{\lambda,\tilde{x}}.
\end{equation}
This concludes the construction of $\mathbb{L}^{\lambda,\tilde{x}}$
and $\tens^{\lambda,\tilde{x}}$. Again by construction the
operators $\tens^{\lambda,\tilde{x}}$ will act transitively on the
connected componets of $\left(
\lan{\gHiggs}_{\widetilde{C}}\right)^{D}$. For future reference we will
write $\hm_{\tilde{x}}^{D} : \Lambda \to
\left(\lan{\gHiggs}_{\widetilde{C}} \right)^{D}$ for the homomorphism
given by $\lambda \mapsto \mathfrak{L}^{\lambda,\tilde{x}}$. Similarly
we will write $\lan{\hm}_{\tilde{x}}^{D} : \lan{\Lambda} \to
\left(\gHiggs_{\widetilde{C}} \right)^{D}$ for the analogous
homomorphism associated with the Langlands dual group. Finally we will
denote the restrictions of these homomorphisms to
$\crts_{\mathfrak{g}} \subset \Lambda$ by $\hm_{\tilde{x}}^{D'} :
\crts_{\mathfrak{g}} \to \left(\lan{\Higgs}_{\widetilde{C}}
\right)^{D}$ and $\lan{\hm}_{\tilde{x}}^{D'} :
\crts_{\lan{\mathfrak{g}}} \to \left(\Higgs_{\widetilde{C}}
\right)^{D}$ respectively.

\

\begin{rem} \label{rem:ab.wlson}
Note that the operation of tensoring with the line bundle
$\mathfrak{L}^{\lambda,\tilde{x}}$ gives rise to an
abelianized tensorization (or Wilson) auto equivalence
\[
\lanab{\mathbb{W}}^{\lambda,\tilde{x}} :
D_{coh}(\lan{\gHiggs}_{\widetilde{C}},\mathcal{O})
\stackrel{\cong}{\longrightarrow}
D_{coh}(\lan{\gHiggs}_{\widetilde{C}},\mathcal{O}).
\]
Similarly to the Hecke operators, the algebra of these auto
equivalences can be related to the algebra of the classical limit
tensorization functors $\lan{\mathbb{W}}^{\lambda,x}$ discussed in
Section~\ref{s:cl}. As with the Hecke action we will not discuss the
precise relationship here but will work with the abelianized versions
only.
\end{rem}

\

\subsection{Proof of global duality} \label{ssec:proofThmB}

From the analysis in section~\ref{ss:universal} we know that the group
of connected components of $\gHiggs$ is naturally $\pi_{1}(G)$ while
the group of connected components of $\left(\lan{\gHiggs} \right)^{D}$
is \linebreak $Z(\lan{G})^{\wedge} =
\op{Hom}(Z(\lan{G}),\mathbb{C}^{\times})$. Thus the groups of connected
components of $\gHiggs$ and $\left(\lan{\gHiggs} \right)^{D}$ are
naturally identified by the isomorphism $\lan{\Lambda} \cong
\Lambda^{\vee}$.  In the previous section we defined the Hecke
action on the stack $\gHiggs$ and the tensorization action on the stack
$\left(\lan{\gHiggs} \right)^{D}$, and we noted that these actions induce
transitive actions on sets of connected components of $\gHiggs$.
Furthermore, the compatibility \eqref{eq:compatible.pi0} guarantees
that the two actions on groups of connected components match  under the
identification $\pi_{1}(G) \cong Z(\lan{G})^{\wedge}$. The
transitivity reduces the problem of constructing the 
isomorphism \eqref{eq:global.iso} to the problem of constructing a
canonical  isomorphism of connected $Z(G)$-gerbes
\begin{equation} \label{eq:global.iso.0.over.B}
\gHiggs_{0}  \stackrel{\cong}{\to}
\left(\lan{\Higgs}\right)^{D},
\end{equation}
which intertwines the actions of all abelianized Hecke and
tensorization operators that preserve these connected components. 

Indeed, the isomorphism \eqref{eq:global.iso.0.over.B} extends
automatically to the desired isomorphism \eqref{eq:global.iso} by
Hecke equivariance. The operators $\trans^{\lambda,\tilde{x}}$ and
$\tens^{\lambda,\tilde{x}}$ are labeled by $(\tilde{x},\lambda) \in
\widetilde{\mycal{C}}\times \Lambda$. The operators preserving the
connected component of the identity on either side are the Hecke or
the tensorization operators labeled by $(\tilde{x},\lambda) \in
\widetilde{\mycal{C}}\times \crts_{\mathfrak{g}}$. Therefore to finish
the proof of Theorem~\ref{thm:gerbes} we must construct a lift of
$\paiso_{\widetilde{C}} : P_{\widetilde{C}} \to
\lan{P}_{\widetilde{C}}^{D}$ from Theorem~\ref{thm:duality} to an
isomorphism
\begin{equation} \label{eq:global.iso.0}
\gHiggs_{0,\widetilde{C}}  \stackrel{\cong}{\to}
\left(\lan{\Higgs}_{\widetilde{C}}\right)^{D},
\end{equation}
which intertwines these Hecke and tensorization operators. Explicitly
this means that for every $\tilde{x} \in \widetilde{C}$ the
isomorphism \eqref{eq:global.iso.0} fits in a commutative diagram of
group stacks
\begin{equation} \label{eq:equivariant.iso.0}
\xymatrix{
& \crts_{\mathfrak{g}} \ar[dl]_-{\hm^{'}_{\tilde{x}}}
  \ar[dr]^-{\hm_{\tilde{x}}^{D'}} & \\ 
\gHiggs_{0,\widetilde{C}}  \ar[rr] & & 
\left(\lan{\Higgs}_{\widetilde{C}}\right)^{D}.
}
\end{equation}
Here we fixed the cameral cover $\widetilde{C}$ to simplify the
exposition. However we will construct the isomorphism
\eqref{eq:global.iso.0} in such a way that the construction will
automatically globalize over $B - \Delta$. The idea is to construct
groupoid presentations of the gerbes $\gHiggs_{0,\widetilde{C}}$
and $\left(\lan{\Higgs}_{\widetilde{C}}\right)^{D}$ and then argue
that these presentations are naturally isomorphic.

To carry this out recall $\gHiggs_{0,\widetilde{C}}$ and
$\left(\lan{\Higgs}_{\widetilde{C}}\right)^{D}$ are neutralizable
$Z(G)$-gerbes with coarse moduli spaces $P_{\widetilde{C}}$ and
$\lan{P}_{\widetilde{C}}^{D}$ respectively. Our strategy is to build
atlases for the gerbes as moduli spaces of neutralizations and then
identify explicitly these moduli spaces with each other by matching
the respective geometric data.

First let us look at the right hand side of
\eqref{eq:global.iso.0}. 

\

\punkt \ {\bfseries A groupoid presentation of
  $\left(\lan{\Higgs}_{\widetilde{C}}\right)^{D}$:} \label{sssec:rhs}
  \ The  Hitchin fiber
$\lan{\Higgs}_{\widetilde{C}}$ is an algebraic group which
fits in a short exact sequence of abelian groups
\begin{equation} \label{eq:rhs.ses}
\xymatrix@1{
0 \ar[r] &  \lan{P}_{\widetilde{C}} \ar[r] &
\lan{\Higgs}_{\widetilde{C}} \ar[r] & 
\pi_{1}(\lan{G}) \ar[r] &  0.
}
\end{equation}
Since $\lan{G}$ is a simple Lie group we have that $\pi_{1}(\lan{G})$
is finite, and since $\lan{P}_{\widetilde{C}} $ is divisible, it
follows that the sequence \eqref{eq:rhs.ses} is split.

Passing to duals we see that the dual gerbe now fits in a short exact
sequence of commutative group stacks
\[
\xymatrix@R-1pc{
0 \ar[r] & \pi_{1}(\lan{G})^{D} \ar[r] \ar@{=}[d] & 
\left(\lan{\Higgs}_{\widetilde{C}}\right)^{D} \ar[r] &
\lan{P}_{\widetilde{C}}^{D} \ar[r]  & 0. \\
& BZ(G) & & &
}
\]
Since $\lan{P}_{\widetilde{C}}^{D}$ is the coarse moduli space of
$\left(\lan{\Higgs}_{\widetilde{C}}\right)^{D}$ it follows that every
splitting of \eqref{eq:rhs.ses} will induce a neutralization of the
gerbe $\left(\lan{\Higgs}_{\widetilde{C}}\right)^{D}$, i.e. a map from
$\left(\lan{\Higgs}_{\widetilde{C}}\right)^{D}$ to the gerbe
$[\lan{P}^{D}/Z(G)]$ (with $Z(G)$ acting trivially) which for every
test scheme $S$ gives rise to an equivalence
$\left(\lan{\Higgs}_{\widetilde{C}}\right)^{D}(S) \to
[\lan{P}^{D}/Z(G)](S)$ between the associated groupoids of
sections. In fact a choice of a splitting of \eqref{eq:rhs.ses} is
equivalent to choosing a neutralization which is compatible with the
group structures.

As we saw in sections \ref{sss:hecke.action} and
\ref{sss:dual.hecke.action}, the homomorphism $\lan{\hm}_{\tilde{x}} :
\lan{\Lambda} \to \lan{\Higgs}_{\widetilde{C}}$ induces a map of short
exact sequences of abelian groups:
\[
\xymatrix@M+0.3pc{ 
0 \ar[r] &  \lan{P}_{\widetilde{C}} \ar[r] &
\lan{\Higgs}_{\widetilde{C}} \ar[r] & 
\pi_{1}(\lan{G}) \ar[r] &  0 \\
0 \ar[r] & \crts_{\lan{\mathfrak{g}}} \ar[r]
\ar@{^{(}->}[u]^-{\lan{\hm}'_{\tilde{x}}}  & \lan{\Lambda} \ar[r]_-{\lan{\bq}}
\ar@{^{(}->}[u]^-{\lan{\hm}_{\tilde{x}}} & \pi_{1}(\lan{G}) \ar[r] \ar@{=}[u]
& 0,
}
\]
where $\lan{\hm}'_{\tilde{x}}$ is the restriction of
$\lan{\hm}_{\tilde{x}}$ to $\crts_{\lan{\mathfrak{g}}}$. Now every
choice of a homomorphism
\[
\bphi : \lan{\Lambda} \to \lan{P}_{\widetilde{C}} 
\]
which lifts $\lan{\hm}'_{\tilde{x}} : \crts_{\lan{\mathfrak{g}}} \to
\lan{P}_{\widetilde{C}}$ will give a splitting of
\eqref{eq:rhs.ses}. Moreover the set of all such lifts is in bijection
with the set of all splittings of \eqref{eq:rhs.ses}.

Thus the choice of a pair $(\tilde{x},\bphi)$ gives a splitting of
\eqref{eq:rhs.ses} and a neutralization of the gerbe $\left(
\lan{\Higgs}_{\widetilde{C}} \right)^{D}$ compatible with the group
structure. The space parametrizing such pairs can be used to build an
atlas of the gerbe $\left( \lan{\Higgs} \right)^{D}$. To avoid
complications with branching we will only use pairs
$(\tilde{x},\bphi)$ in which $\tilde{x}$ is not a ramification point
for the cover $\widetilde{C} \to C$. Specifically let
$\widetilde{C}^{0} \subset \widetilde{C}$ be the complement of the
ramification divisor of the map $\widetilde{C} \to C$. Define a cover
$\widetilde{C}^{0}_{\RHS} \to \widetilde{C}^{0}$ parametrizing all
pairs $(\tilde{x},\bphi)$ as above:
\[
\widetilde{C}^{0}_{\RHS} := \left\{ (\tilde{x},\bphi) \left| 
\text{
\begin{minipage}[c]{1.8in}
$\tilde{x} \in \widetilde{C}^{0}$, 
$\bphi : \lan{\Lambda} \to \lan{P}_{\widetilde{C}} $,
  so that $\bphi_{|\crts_{\lan{\mathfrak{g}}}} = \lan{\hm}'_{\tilde{x}}$.
\end{minipage}
}
\right.\right\}.
\]
Denote by $\widetilde{\mycal{C}}^{0}_{\RHS} \to (B-\Delta)$ the family of all
$\widetilde{C}^{0}_{\RHS} $ viewed as a cover of an open subset of 
the universal cameral cover.
With this notation we get
\begin{lem} \label{lem:rhs.atlas} The variety 
\[
\mathcal{U}_{\RHS} :=  \widetilde{\mycal{C}}^{0}_{\RHS} \, \times_{B-\Delta}
\left(
\lan{\Higgs}_{0} \right)^{D}
\]
is an atlas for the $Z(G)$-gerbe $\left(
\lan{\Higgs} \right)^{D}$.
\end{lem}
{\bf Proof.} \ The statement is relative over the Hitchin base so we
can argue fiber by fiber.  By definition $\widetilde{C}^{0}_{\RHS}$ is
a etale Galois cover of $\widetilde{C}^{0}$ with Galois group
$\op{Hom}(\pi_{1}(\lan{G}),\lan{P}_{\widetilde{C}})$.  The assignment
$(\tilde{x},\bphi) \to \bphi$ gives a map from
$\widetilde{C}^{0}_{\RHS}$ to the space of splittings of
\eqref{eq:rhs.ses} and so the pull back of the gerbe $\left(
\lan{\Higgs}_{\widetilde{C}} \right)^{D}$ to
$\widetilde{C}^{0}_{\RHS}\times \lan{P}_{\widetilde{C}}^{D}$ is
equipped with a universal neutralization. Thus
$\widetilde{C}^{0}_{\RHS}\times \lan{P}_{\widetilde{C}}^{D}$ is a
smooth atlas for $\left( \lan{\Higgs}_{\widetilde{C}} \right)^{D}$
with a structure morphism $\widetilde{C}^{0}_{\RHS}\times
\lan{P}_{\widetilde{C}}^{D} \to \left( \lan{\Higgs}_{\widetilde{C}}
\right)^{D}$ given by the universal neutralization. \ \hfill $\Box$

\

\medskip

\noindent
The structure map 
$\mathcal{U}_{\RHS} \to \left( \lan{\Higgs}
\right)^{D}$ gives rise to a groupoid presentation 
\[
\xymatrix@1{
\mathcal{U}_{\RHS} \times_{\left(
\lan{\Higgs}\right)^{D}} \mathcal{U}_{\RHS} \ar@<.8ex>[r]
\ar@<.1ex>[r] & \mathcal{U}_{\RHS} }
\]
of the gerbe $\left( \lan{\Higgs} \right)^{D}$.

Denote the
space $\mathcal{U}_{\RHS}\times_{\left( \lan{\Higgs}\right)^{D}}
\mathcal{U}_{\RHS}$ by $\mathcal{R}_{\RHS}$. We can describe the space 
$\mathcal{R}_{\RHS}$ explicitly. The gerbe $\left(
\lan{\Higgs}\right)^{D}$ has a coarse moduli space $\left(
\lan{\Higgs}_{0} \right)^{D}$ and so $\mathcal{R}_{\RHS}$ will be a
$Z(G)$-torsor over
\[
\mathcal{U}_{\RHS}\times_{\left( \lan{\Higgs}_{0} \right)^{D}}
\mathcal{U}_{\RHS} = \widetilde{\mycal{C}}^{0}_{\RHS}\times_{B-\Delta}
\widetilde{\mycal{C}}^{0}_{\RHS}\times_{B-\Delta} \left(\lan{\Higgs}_{0}
\right)^{D}.
\]
A point in this space mapping to $\widetilde{C} \in
B-\Delta$ consists of the data \linebreak $((\tilde{x}_{1},\bphi_{1}),
(\tilde{x}_{2},\bphi_{2}), \mathbb{L})$ where
$(\tilde{x}_{1},\bphi_{1})$, $(\tilde{x}_{2},\bphi_{2})$ are in the
cover $\widetilde{C}^{0}_{\RHS} \to \widetilde{C}^{0}$, and
$\mathbb{L}$ is a point of $\lan{P}_{\widetilde{C}}^{D}$. 

To any such point we can attach group homomorphisms
$\lan{\hm}_{\tilde{x}_{1}}, \lan{\hm}_{\tilde{x}_{2}} : \lan{\Lambda} \to
\lan{\Higgs}_{\widetilde{C}}$, and $\bphi_{1}, \bphi_{2} :
\lan{\Lambda}  \to \lan{P}_{\widetilde{C}}$. The fact that
$\widetilde{C} \in B - \Delta$ implies that $\widetilde{C}$ is
smooth and irreducible. Hence   $\widetilde{C}^{0}$ is connected and
so $\lan{\hm}_{\tilde{x}_{1}} - \lan{\hm}_{\tilde{x}_{2}}$  will map
$\lan{\Lambda}$ to the connected component of the identity of the
group  $\lan{\Higgs}_{\widetilde{C}}$. Since $\bphi_{i}$ agree with
$\lan{\hm}_{\tilde{x_{i}}}$ on $\crts_{\lan{\mathfrak{g}}}$ it follows that
  the homomorphism 
\[
(\lan{\hm}_{\tilde{x}_{1}} - \lan{\hm}_{\tilde{x}_{2}}) - (\bphi_{1} - \bphi_{2})
: \lan{\Lambda} \to \lan{P}_{\widetilde{C}}
\]
will factor through a homomorphism 
\begin{equation} \label{eq:nu}
\nu((\tilde{x}_{1},\bphi_{1}),
(\tilde{x}_{2},\bphi_{2})) : \pi_{1}(\lan{G}) \to
\lan{P}_{\widetilde{C}}. 
\end{equation}
Now since $\lan{P}_{\widetilde{C}}^{D} =
\op{Hom}_{\text{gp-stack}}(\lan{P}_{\widetilde{C}},B\mathbb{G}_{m})$
we can view $\mathbb{L}$ as an abelian group extension 
\[
0 \to \mathbb{C}^{\times} \to \mathbb{L} \to \lan{P}_{\widetilde{C}}
\to 0.
\]
Pulling back this extension by the map $\nu((\tilde{x}_{1},\bphi_{1}),
(\tilde{x}_{2},\bphi_{2}))$ gives an extension of
$\pi_{1}(\lan{G})$ by $\mathbb{C}^{\times}$. Since $\pi_{1}(\lan{G})$
is finite and $\mathbb{C}^{\times}$ is divisible this extension is
necessarily split. We have obtained the following
description of $\mathcal{R}_{\RHS}$: 

\begin{lem} \label{lem:rhs.torsor} $\mathcal{R}_{\RHS}$ is a
  $Z(G)$-torsor over   
$\mathcal{U}_{\RHS}\times_{\left( \lan{\Higgs}_{0} \right)^{D}}
\mathcal{U}_{\RHS}$ whose fiber over a point
is the set of all splittings of the group extension 
\[
0 \to \mathbb{C}^{\times} \to \nu((\tilde{x}_{1},\bphi_{1}),
(\tilde{x}_{2},\bphi_{2}))^{*}\mathbb{L} \to \pi_{1}(\lan{G}) \to 0,
\]
where $\nu((\tilde{x}_{1},\bphi_{1}),
(\tilde{x}_{2},\bphi_{2}))$ is defined in formula \eqref{eq:nu}.
\end{lem}

\

\medskip

\punkt {\bfseries A groupoid presentation of
  $\gHiggs_{0,\widetilde{C}}$:} \label{sssec:lhs} \ The gerbe
  $\gHiggs_{0,\widetilde{C}}$ fits in a short exact sequence of
  commutative group stacks
\begin{equation} \label{eq:lhs.ses}
\xymatrix@1{
0 \ar[r] & BZ(G) \ar[r] & \gHiggs_{0,\widetilde{C}} \ar[r]^-{\bpr} &
P_{\widetilde{C}} \ar[r] & 0.
}
\end{equation} 
Following the same pattern as above we will construct an atlas for
$\gHiggs_{0,\widetilde{C}}$ from the moduli space of all splittings of
\eqref{eq:lhs.ses}. To construct sections of $\bpr$ we will again use
points of the cameral cover $\widetilde{C}$. 

Since by \cite[Section~6]{ron-dennis} $\gHiggs_{0,\widetilde{C}}$ can
be identified with the moduli stack of (topologically trivial)
$G$-spectral data $(\mycal{L},\bi,\bb)$ we have a universal family of
spectral data
\[
(\underline{\mycal{L}},\underline{\bi},\underline{\bb}) \to
\gHiggs_{0,\widetilde{C}} \times \widetilde{C}.
\] 
Let now $\tilde{x} \in \widetilde{C}^{0}$ be a point away from
 ramification. Restricting $\underline{\mycal{L}}$ to
 $\gHiggs_{0,\widetilde{C}} \times \{ \tilde{x} \}$ we get a
 $T$-bundle $\underline{\mycal{L}}_{\tilde{x}}$ with square structure,
 i.e. an extension 
\[
0 \to T \to \underline{\mycal{L}}_{\tilde{x}} \to
\gHiggs_{0,\widetilde{C}} \to 0
\]
of commutative group stacks.

\

\begin{rem} \label{eq:square.on.L}
This square structure is encoded in the interpretation of
$\gHiggs_{0,\widetilde{C}}$ as moduli of spectral data. To see this
suppose that $S$ is a test scheme.  A section $\xi : S \to
\gHiggs_{0,\widetilde{C}}$ of the Higgs stack over $S$ is given by a
family of spectral data $(\mycal{L}^{\xi},\bi^{\xi},\bb^{\xi})$ on
$S\times \widetilde{C}$.  Since
$\xi^{*}((\underline{\mycal{L}},\underline{\bi},\underline{\bb})) =
(\mycal{L}^{\xi},\bi^{\xi},\bb^{\xi})$ it follows that
$\xi^{*}\underline{\mycal{L}}_{\tilde{x}} = \mycal{L}^{\xi}_{|S\times
  \{ \tilde{x} \}} =: \mycal{L}^{\xi}_{\tilde{x}}$. If now $\xi_{1}$
and $\xi_{2}$ are two sections of $\gHiggs_{0,\widetilde{C}}$ over
$S$, then we can add $\xi_{1}$ and $\xi_{2}$ in the group structure on
$\gHiggs_{0,\widetilde{C}}$ to get a new section $\xi_{1} + \xi_{2} :
S \to \gHiggs_{0,\widetilde{C}}$. But in terms of spectral data the group
structure is given by tensoring the corresponding $T$-torsors, twists,
and framings. Thus 
\[
\begin{split}
(\xi_{1} +
\xi_{2})^{*}((\underline{\mycal{L}},\underline{\bi},\underline{\bb}))
& = \left(\mycal{L}^{\xi_{1}},\bi^{\xi_{1}},\bb^{\xi_{1}}\right)\cdot
\left(\mycal{L}^{\xi_{2}},\bi^{\xi_{2}},\bb^{\xi_{2}}\right) \\
& = \left( \mycal{L}^{\xi_{1}}\otimes \mycal{L}^{\xi_{2}},
\bi^{\xi_{1}}\otimes \op{id} + \op{id}\otimes \bi^{\xi_{2}},
\bb^{\xi_{1}}\otimes \op{id} + \op{id}\otimes \bb^{\xi_{2}}\right),
\end{split} 
\]
and so 
\[
\begin{split}
(\xi_{1}+\xi_{2})^{*}\underline{\mycal{L}}_{\tilde{x}} & =
  \mycal{L}^{\xi_{1}+\xi_{2}}_{\tilde{x}} \\
& = \mycal{L}^{\xi_{1}}\otimes \mycal{L}^{\xi_{2}} \\
& = \left(\xi_{1}^{*}\underline{\mycal{L}}_{\tilde{x}}\right)\otimes
  \left(\xi_{2}^{*}\underline{\mycal{L}}_{\tilde{x}}\right).
\end{split}
\]
This identification gives the multiplication isomorphism in the square
structure. The idenitity isomorphism in the square structure is
defined analogously.
\end{rem}

\

\begin{lem} \label{lem:lhs.choce.of.lift} The following data are
canonically  equivalent:
\begin{itemize}
\item[(a)] A neutralization of the $Z(G)$-gerbe $\gHiggs_{0,\widetilde{C}} \to
P_{\widetilde{C}}$, compatible with the group structures;
\item[(b)] A splitting of \eqref{eq:lhs.ses};
\item[(c)] A lift of
$\underline{\mycal{L}}^{\, '}_{\tilde{x}}$ to a $T$-torsor (with
square structure) on $P_{\widetilde{C}}$.
\end{itemize}
\end{lem}
{\bfseries Proof.} The data (a) and (b) are tautologically the
same. 

\

\noindent
To identify data (b) and (c) note that by construction, the stabilizer
$Z(G)$ of any section $\xi : S \to \gHiggs_{0,\widetilde{C}}$ acts on
$\underline{\mycal{L}}_{\tilde{x}}$ via the canonical inclusion $Z(G)
\subset T$. In other words we can view
$\underline{\mycal{L}}_{\tilde{x}}$ as a neutralization of the
$T$-gerbe on $P_{\widetilde{C}}$ induced from the $Z(G)$-gerbe
$\gHiggs_{0,\widetilde{C}} \to P_{\widetilde{C}}$. To obtain a
neutralization of the $Z(G)$-gerbe we will have to choose additional
data. Specifically, the natural map from $\gHiggs_{0,\widetilde{C}}$
to the induced $T$-gerbe $\gHiggs_{0,\widetilde{C}}\times_{BZ(G)} BT$
is a $T/Z(G)$-torsor. Pulling back this $T/Z(G)$ torsor to
$P_{\widetilde{C}}$ via the map 
$\sigma_{\tilde{x}} : P_{\widetilde{C}} \to
\gHiggs_{0,\widetilde{C}}\times_{BZ(G)} BT$ corresponding to
$\underline{\mycal{L}}_{\tilde{x}}$ gives a $T/Z(G)$-torsor
$\underline{\mycal{L}}^{\, '}_{\tilde{x}}$ (with square structure) on
$P_{\widetilde{C}}$:
\[
\xymatrix{ \ar@{}[dr]|{\Box}
\underline{\mycal{L}}^{\, '}_{\tilde{x}} \ar[r] \ar[d] &
\gHiggs_{0,\widetilde{C}} \ar[d] \\
P_{\widetilde{C}} \ar[r]_-{\sigma_{\tilde{x}}} &
\gHiggs_{0,\widetilde{C}}\times_{BZ(G)} BT 
}
\]
In this way we obtain an identification of $\gHiggs_{0,\widetilde{C}}$
with the gerbe of all lifts of $\underline{\mycal{L}}^{\,
'}_{\tilde{x}}$ to a $T$-torsor with square structure on
$P_{\widetilde{C}}$. Therefore a splitting of \eqref{eq:lhs.ses} is
the same thing as a choice of such a lift. \ \hfill $\Box$

\

\medskip

\noindent
To package all the above data  in the most efficient way, note first that
the  fact that the stabilizer action on
$\underline{\mycal{L}}_{\tilde{x}}$ is tautological implies that the
$T/Z(G)$-torsor with square structure induced from
$\underline{\mycal{L}}_{\tilde{x}}$ will descend to a $T/Z(G)$-torzor on
$P_{\widetilde{C}}$. In terms of group extensions this means that the
  induced extension
\[
\xymatrix@1{
0 \ar[r] & T/Z(G) \ar[r] & \underline{\mycal{L}}_{\tilde{x}}/Z(G) \ar[r] &
\gHiggs_{0,\widetilde{C}} \ar[r]  & 0 
}
\]
is canonically a pullback of an extension
\[
\xymatrix@1{
0 \ar[r] & T/Z(G) \ar[r]  &  \underline{\mycal{L}}^{\, '}_{\tilde{x}} \ar[r]  &
P_{\widetilde{C}} \ar[r]  & 0
}
\] 
via the structure homomorphism $\bpr : \gHiggs_{0,\widetilde{C}} \to
P_{\widetilde{C}}$. This extension is precisely the $T/Z(G)$ torsor we
described before the statement of Lemma~\ref{lem:lhs.choce.of.lift}
and so a splitting of \eqref{eq:lhs.ses} corresponds to a pair
$(\tilde{x},\bPhi)$ where $\bPhi$ is an extension
\[
\xymatrix@1{
0 \ar[r] & T  \ar[r] & \bPhi \ar[r]  & P_{\widetilde{C}} \ar[r] & 0 
}
\]
together with a choice of a group isomorphism $\bPhi/Z(G)
\stackrel{\cong}{\to} \underline{\mycal{L}}^{\, '}_{\tilde{x}}$
inducing an isomorphism of extensions:
\[
\xymatrix@R-1pc{
0 \ar[r] & T/Z(G) \ar[r] \ar@{=}[d] &  \bPhi/Z(G) \ar[r] \ar[d]^-{\cong} & 
P_{\widetilde{C}} \ar[r] \ar@{=}[d] & 0 \\ 
0 \ar[r] & T/Z(G) \ar[r]  &  \underline{\mycal{L}}^{\, '}_{\tilde{x}} \ar[r]  &
P_{\widetilde{C}} \ar[r]  & 0.
}
\]
Note that by the five lemma the isomorphism $\bPhi/Z(G)
\stackrel{\cong}{\to} \underline{\mycal{L}}^{\, '}_{\tilde{x}}$ will
be unique if it exists. So the choice of such an isomorphism will not
be extra data. 

The space parametrizing pairs $(\tilde{x},\bPhi)$ can  be used
to build an atlas for the gerbe $\gHiggs_{0,\widetilde{C}}$. To that
end define an etale cover $\widetilde{C}^{0}_{\LHS} \to \widetilde{C}$
parametrizing all pairs $(\tilde{x},\bPhi)$, i.e. 
\[
\widetilde{C}^{0}_{\LHS} := \left\{ (\tilde{x},\bPhi) \left| 
\text{
\begin{minipage}[c]{3in}
$\tilde{x} \in \widetilde{C}^{0}$, $\bPhi$ is an extension of
  $P_{\widetilde{C}}$ by $T$, 
  so that the extension 
\[
\xymatrix@1{
0 \ar[r] & T/Z(G) \ar[r] &  \bPhi/Z(G) \ar[r]  & 
P_{\widetilde{C}} \ar[r]  & 0
}
\]
is isomorphisc to 
\[
\xymatrix@1{
0 \ar[r] & T/Z(G) \ar[r]  &  \underline{\mycal{L}}^{\, '}_{\tilde{x}} \ar[r]  &
P_{\widetilde{C}} \ar[r]  & 0.
}
\]
\end{minipage}
}
\right.\right\}.
\]
Denote by $\widetilde{\mycal{C}}^{0}_{\LHS} \to (B-\Delta)$ the family of all
$\widetilde{C}^{0}_{\LHS} $ viewed as a cover of the universal cameral cover.
With this notation we now have
\begin{lem} \label{lem:lhs.atlas} The variety 
\[
\mathcal{U}_{\LHS} :=  \widetilde{\mycal{C}}^{0}_{\LHS} \, \times_{B-\Delta}
\Higgs_{0} 
\]
is an atlas for the $Z(G)$-gerbe $\gHiggs_{0,\widetilde{C}}$.
\end{lem}
{\bf Proof.} \ Again the statement is relative over the Hitchin base
and we can argue fiber by fiber.  The assignment $(\tilde{x},\bPhi)
\to \bPhi$ gives a map from $\widetilde{C}^{0}_{\LHS}$ to the space of
splittings of \eqref{eq:lhs.ses} and so the pull back of the gerbe
$\gHiggs_{0,\widetilde{C}}$ to $\widetilde{C}^{0}_{\LHS}\times
P_{\widetilde{C}}$ is equipped with a universal neutralization. Thus
$\widetilde{C}^{0}_{\LHS}\times P_{\widetilde{C}}$ is a smooth atlas
for $\gHiggs_{0,\widetilde{C}}$ with a structure morphism
$\widetilde{C}^{0}_{\LHS}\times P_{\widetilde{C}} \to
\gHiggs_{0,\widetilde{C}}$ given by the universal neutralization. \
\hfill $\Box$

\

\medskip

\noindent
As usual the structure map $\mathcal{U}_{\LHS} \to
\gHiggs_{0,\widetilde{C}}$ gives rise to a groupoid presentation
\[
\xymatrix@1{
\mathcal{U}_{\LHS}\times_{\gHiggs_{0,\widetilde{C}}} \mathcal{U}_{\LHS}
\ar@<.8ex>[r]
\ar@<.1ex>[r] & \mathcal{U}_{\LHS}
}
\]
of the gerbe $\gHiggs_{0,\widetilde{C}}$. The relations
\[
\mathcal{R}_{\LHS} :=
\mathcal{U}_{\LHS}\times_{\gHiggs_{0,\widetilde{C}}} \mathcal{U}_{\LHS} 
\]
for this presentation form a variety, which is a
$Z(G)$-torsor over the space
\[
\mathcal{U}_{\LHS}\times_{\Higgs_{0}}
\mathcal{U}_{\LHS} = \widetilde{\mycal{C}}^{0}_{\LHS}\times_{B-\Delta}
\widetilde{\mycal{C}}^{0}_{\LHS}\times_{B-\Delta} \Higgs_{0}.
\]
A point in this space sitting over a cameral cover $\widetilde{C} \in
B-\Delta$ consists of the data \linebreak $((\tilde{x}_{1},\bPhi_{1}),
(\tilde{x}_{2},\bPhi_{2}), \mathbb{L})$ where
$(\tilde{x}_{1},\bPhi_{1})$, $(\tilde{x}_{2},\bPhi_{2})$ are in the
cover $\widetilde{C}^{0}_{\LHS} \to \widetilde{C}^{0}$, and
$\mathbb{L}$ is a point of $P_{\widetilde{C}}$.

Since the stabilizer stack of $\gHiggs_{0,\widetilde{C}}$ acts
tautologically on each of the $Z(G)$-torsors
$\underline{\mycal{L}}_{\tilde{x}_{1}}$ and
$\underline{\mycal{L}}_{\tilde{x}_{2}}$, it follows that the
$T$-torsor
$\underline{\mycal{L}}_{\tilde{x}_{1}}\otimes
\underline{\mycal{L}}^{-1}_{\tilde{x}_{2}}$ 
fits in a group stack 
extension
\[
\xymatrix@1{ 
0 \ar[r] & T \ar[r] &
\underline{\mycal{L}}_{\tilde{x}_{1}}\otimes
\underline{\mycal{L}}^{-1}_{\tilde{x}_{2}} \ar[r] &
\gHiggs_{0,\widetilde{C}} \ar[r] & 0 
}
\]
which is a pullback of an extension 
\[
\xymatrix@1{ 0 \ar[r] & T \ar[r] &
\mathcal{M}_{\tilde{x}_{1},\tilde{x}_{2}} \ar[r] & P_{\widetilde{C}}
\ar[r] & 0 }
\]
via the map $\bpr : \gHiggs_{0,P_{\widetilde{C}}} \to
  P_{\widetilde{C}}$. Therefore the $T$-torsor
\begin{equation} \label{eq:N} 
N\left((\tilde{x}_{1},\bPhi_{1}),(\tilde{x}_{2},\bPhi_{2})\right) := 
\mathcal{M}_{\tilde{x}_{1},\tilde{x}_{2}}\otimes \left(
  \bPhi_{1}^{-1}\otimes \bPhi_{2}\right)
\end{equation}
fits in a group extension
\[
\xymatrix@1{ 0 \ar[r] & T \ar[r] &
N\left((\tilde{x}_{1},\bPhi_{1}),(\tilde{x}_{2},\bPhi_{2})\right)
\ar[r] & P_{\widetilde{C}} \ar[r] & 0. }
\]
Since the $T/Z(G)$-torsors induced from  
$\mathcal{M}_{\tilde{x}_{1},\tilde{x}_{2}}$ and $\bPhi_{1}\otimes
\bPhi_{2}^{-1}$ are both equal to $\underline{\mycal{L}}^{\,
  '}_{\tilde{x_{1}}}\otimes \underline{\mycal{L}}^{\,
  '}_{\tilde{x_{2}}}$ it follows that the extension
\begin{equation} \label{eq:split.induced}
\xymatrix@1{ 0 \ar[r] & T/Z(G) \ar[r] &
N\left((\tilde{x}_{1},\bPhi_{1}),(\tilde{x}_{2},\bPhi_{2})\right)/Z(G)
\ar[r] & P_{\widetilde{C}} \ar[r] & 0. }
\end{equation}
is split. Furthermore, since  $P_{\widetilde{C}}$ is an abelian
variety, and  $T/Z(G)$ is affine it follows that
\eqref{eq:split.induced} will have a unique splitting $\sigma :
P_{\widetilde{C}}  \to
N\left((\tilde{x}_{1},\bPhi_{1}),(\tilde{x}_{2},\bPhi_{2})\right)/Z(G)$
Pulling back the $Z(G)$-torsor
$N\left((\tilde{x}_{1},\bPhi_{1}),(\tilde{x}_{2},\bPhi_{2})\right) \to
N\left((\tilde{x}_{1},\bPhi_{1}),(\tilde{x}_{2},\bPhi_{2})\right)/Z(G)$
by $\sigma$ gives a $Z(G)$-torsor on $P_{\widetilde{C}}$. Now the
definition of the relations $\mathcal{R}_{\LHS}$ immediately gives the
following 
\begin{lem} \label{lem:lhs.torsor} $\mathcal{R}_{\LHS}$ is a
  $Z(G)$-torsor over   
$\mathcal{U}_{\LHS}\times_{\Higgs_{0}}
\mathcal{U}_{\LHS}$ whose fiber over a point \linebreak
$\left((\tilde{x}_{1},\bPhi_{1}),(\tilde{x}_{2},\bPhi_{2}),\mathbb{L}\right)$
is the fiber of the $Z(G)$-torsor
\[
\sigma^{*} N\left((\tilde{x}_{1},\bPhi_{1}),(\tilde{x}_{2},\bPhi_{2})\right)
\]
at $\mathbb{L} \in P_{\widetilde{C}}$.
\end{lem} 

\

\bigskip

\

\punkt {\bfseries The construction of the isomorphism:}
\label{ssec:extend.iso} \ The duality
of Higgs stacks will follow if we can show that our isomorphism of
Prym varieties lifts to an isomorphism between the groupoid
presentations we described in Sections \ref{sssec:lhs} and
\ref{sssec:rhs}. Thus the existence  of the isomorphism
\eqref{eq:global.iso.0} will follow from the following

\begin{theo} \label{theo:iso.presentations}
{\bfseries (a)} \ For all $\tilde{x} \in \widetilde{\mycal{C}}$,
$\lambda \in \crts_{\mathfrak{g}}$, the action of the Hecke
operators $\trans^{\lambda,\tilde{x}}$, on $\gHiggs_{0}$ and the action
of the tensorizaton operators $\tens^{\lambda,\tilde{x}}$ on
$\left(\lan{\Higgs}\right)^{D}$ lift to actions on the groupoids
$\xymatrix@1{\mathcal{R}_{\LHS} \ar@<.8ex>[r] \ar@<.1ex>[r] &
\mathcal{U}_{\LHS}}$ and $\xymatrix@1{\mathcal{R}_{\RHS} \ar@<.8ex>[r]
\ar@<.1ex>[r] & \mathcal{U}_{\RHS}}$ respectively.

\

\noindent 
{\bfseries (b)} \ The isomorphism 
\[
\paiso  : \Higgs_{0} 
\stackrel{\cong}{\to} \left(\lan{\Higgs}_{0}\right)^{D}
\]
of abelian schemes over $B -\Delta$ lifts to a canonical isomorphism
of groupoids 
\[
\xymatrix{
\mathcal{R}_{\LHS} \ar@<.8ex>[d]
\ar@<.1ex>[d] \ar[r]^-{\mathfrak{r}}_-{\cong} &  \mathcal{R}_{\RHS}
\ar@<.8ex>[d] 
\ar@<.1ex>[d] \\
\mathcal{U}_{\LHS} \ar[d] \ar[r]^-{\mathfrak{u}}_-{\cong} &
\mathcal{U}_{\RHS}  \ar[d] \\
\Higgs_{0}  \ar[r]^-{\paiso}_-{\cong} &
\left(\lan{\Higgs}_{0}\right)^{D} 
}
\]
which intertwines 
$\trans^{\lambda,\tilde{x}}$ with $\tens^{\lambda,\tilde{x}}$.
\end{theo} 
{\bfseries Proof.} The statement is relative over the Hitchin base and
so it suffices to construct the isomorphism of groupoids canonically
on every Hitchin fiber.

The lifts of the Hecke and the tensorization operators claimed in part
{\bfseries (a)} are pre-built in the constructions of the two groupoid
presenations. So the only thing to do is to construct the isomorphisms 
$\mathfrak{u}$ and $\mathfrak{r}$ in part {\bfseries (b)}.

Fix $\widetilde{C} \in B-\Delta$. Write $U_{\LHS,\widetilde{C}}$,
$R_{\LHS,\widetilde{C}}$ for the restriction of the atlas and relations
on the Hitchin fiber over $\widetilde{C}$. We want to construct a
isomorphism of groupoids 
\[
\xymatrix{
R_{\LHS,\widetilde{C}} \ar@<.8ex>[d]
\ar@<.1ex>[d] \ar[r]^-{\mathfrak{r}_{\widetilde{C}}} &  R_{\RHS,\widetilde{C}}
\ar@<.8ex>[d] 
\ar@<.1ex>[d] \\
U_{\LHS,\widetilde{C}} \ar[d] \ar[r]^-{\mathfrak{u}_{\widetilde{C}}} &
U_{\RHS,\widetilde{C}}  \ar[d] \\
P_{\widetilde{C}}  \ar[r]^-{\paiso_{\widetilde{C}}} &
\lan{P}_{\widetilde{C}}^{D} 
}
\]
which lifts $\paiso_{\widetilde{C}}$ and is Hecke equivariant. 

Recall that $U_{\LHS,\widetilde{C}} = \widetilde{C}^{0}_{\LHS}\times
P_{\widetilde{C}}$ and $U_{\RHS,\widetilde{C}} =
\widetilde{C}^{0}_{\RHS}\times \lan{P}_{\widetilde{C}}^{D}$. We will define
the isomorphism of atlases $\mathfrak{u}_{\widetilde{C}} :
U_{\LHS,\widetilde{C}} \to U_{\RHS,\widetilde{C}}$ as a product 
$\mathfrak{u}_{\widetilde{C}} = \mathfrak{a}_{\widetilde{C}}\times
\paiso_{\widetilde{C}}$, where 
\[
\mathfrak{a}_{\widetilde{C}} : \widetilde{C}^{0}_{\LHS} \to
\widetilde{C}^{0}_{\RHS}
\]
is the isomorphism of
$\op{Hom}(\pi_{1}(\lan{G}),\lan{P}_{\widetilde{C}})$-Galois covers of
$\widetilde{C}^{0}$ defined as follows. 

Given a point $(\tilde{x},\bPhi) \in
\widetilde{C}^{0}_{\LHS}$ and an element $\mu \in \lan{\Lambda} =
\chr(T)$ we get an associated element $\mu_{*}\bPhi  \in
\op{Ext}^{1}(P_{\widetilde{C}}, \mathbb{C}^{\times}) =
P_{\widetilde{C}}^{D}$: 
\[
\xymatrix@R-1pc{
0 \ar[r] & T \ar[r] \ar[d]^-{\mu} & \bPhi \ar[r] \ar[d] &
P_{\widetilde{C}} \ar[r]  \ar@{=}[d] & 0 \\
0 \ar[r] & \mathbb{C}^{\times} \ar[r] & \mu_{*}\bPhi \ar[r] &
P_{\widetilde{C}} \ar[r]  & 0.
}
\]
Combined with the isomorphism
$\left(\paiso_{\widetilde{C}}^{D}\right)^{-1} :
P_{\widetilde{C}}^{D} \to \lan{P}_{\widetilde{C}}$  this construction
gives rise to a homomorphism
\[
\bphi_{\bPhi} : \lan{\Lambda} \to \lan{P}_{\widetilde{C}}, \quad \mu
\mapsto \left(\paiso_{\widetilde{C}}^{D}\right)^{-1}(\mu_{*}\bPhi).
\]
We set 
\[
\mathfrak{a}_{\widetilde{C}}(\tilde{x},\bPhi) := (\tilde{x},\bphi_{\bPhi}).
\]
It is straightforward to check that this map is an isomorphism but for
future reference it will be useful to write the inverse map
explicitly. 

Let $(\tilde{x},\bphi) \in \widetilde{C}^{0}_{\RHS}$. Let $\mycal{P}
\to P_{\widetilde{C}}\times \lan{P}_{\widetilde{C}}$ be the Poincare
$\mathbb{C}^{\times}$-torsor corresponding to the isomorphism
$\paiso_{\widetilde{C}}$. Consider the pullback $(\op{id}\times
\bphi)^{*}\mycal{P}$. It is a $\mathbb{C}^{\times}$-torsor on
$P_{\widetilde{C}}\times \lan{\Lambda}$, and since $\mycal{P}$ is a
biextension, it follows that for all $\mu \in \lan{\Lambda}$ the
restriction  $((\op{id}\times
\bphi)^{*}\mycal{P})_{|P_{\widetilde{C}}\times \{ \mu \}}$ is a
  $\mathbb{C}^{\times}$-torsor with square structure and that the map
\[
0 \to \mathbb{G}_{m} \to ((\op{id}\times
\bphi)^{*}\mycal{P} \to
P_{\widetilde{C}}\times \lan{\Lambda} \to 0,
\]
is a short exact sequence of commutative group schemes over
$P_{\widetilde{C}}$. Equivalently this means that $((\op{id}\times
\bphi)^{*}\mycal{P}$ is a $\op{Hom}(\lan{\Lambda},\mathbb{C}^{\times}) = T$
torsor on $P_{\widetilde{C}}$ with square structure. It is now
immediate to see that
\[
\mathfrak{a}_{\widetilde{C}}^{-1}(\tilde{x},\bphi) =
(\tilde{x},((\op{id}\times \bphi)^{*}\mycal{P}).
\]
\

Next we need to lift the isomorphism $\mathfrak{u}_{\widetilde{C}}$ to
an isomorphism of relations. To that end, note that if
$(\tilde{x}_{1},\bPhi_{1})$ and $(\tilde{x}_{2},\bPhi_{2})$ are two
points of $\widetilde{C}^{0}_{\LHS}$, then tautologically we have
\begin{equation} \label{eq:relation.N.nu}
\begin{split}
\bphi_{N((\tilde{x}_{1},\bPhi_{1}),(\tilde{x}_{2},\bPhi_{2}))}  & =  
(\lan{\hm}_{\tilde{x}_{1}} - \lan{\hm}_{\tilde{x}_{2}}) - (\bphi_{\bPhi_{1}} -
\bphi_{\bPhi{2}}) \\ & = \nu((\tilde{x}_{1},\bphi_{\bPhi{1}}),
(\tilde{x}_{2},\bphi_{\bPhi_{2}}))\circ \lan{\bq}
\end{split}
\end{equation}
where $N((\tilde{x}_{1},\bPhi_{1}),(\tilde{x}_{2},\bPhi_{2}))$ and
$\nu((\tilde{x}_{1},\bphi_{\bPhi{1}}),
(\tilde{x}_{2},\bphi_{\bPhi_{2}}))$ are defined by \eqref{eq:N} and
\eqref{eq:nu} respectively, and $\lan{\bq}$ is the natural projection
$\lan{\bq} : \lan{\Lambda}  \to \pi_{1}(\lan{G})$. 

In view of the identification \eqref{eq:relation.N.nu} and
Lemma~\ref{lem:lhs.torsor} and Lemma~\ref{lem:rhs.torsor} the
construction of the isomorphism $\mathfrak{r}_{\widetilde{C}}$ now
reduces to a general statement about abelian varieties which we
formulate next.

Suppose that $P$ is a polarized abelian variety and let $P^{D}$ be the
dual abelian variety. Let 
\[
0 \to T \to N \to P \to 0
\]
be a short exact sequence of commutative algebraic groups and suppose
that the induced sequence 
\begin{equation} \label{eq:N.sequence}
0 \to T/Z(G) \to N/Z(G) \to P \to 0
\end{equation}
is split. Let $\mathbb{L} \in P$ be a point. The data $(N,\mathbb{L})$
gives rise to two $Z(G)$-torsors (over a point):
\begin{itemize}
\item[{\bfseries (i)}] The splittings of the sequence
  \eqref{eq:N.sequence}  form a torsor over
the algebraic group homomorphisms $\op{Hom}(P,T/Z(G))$. Since
$T/Z(G)$ is affine it follows that we have a unique splitting of
  \eqref{eq:N.sequence} which we will denote by $\sigma
: P \to N/Z(G)$. The group $N$ is a $Z(G)$-torsor over $N/Z(G)$ and so
  the pullback $\sigma^{*}N$ is a $Z(G)$-torsor over $P$. The fiber
  $\left(\sigma^{*}N\right)_{\mathbb{L}}$ of
  this torsor at $\mathbb{L}$ is a $Z(G)$-torsor over a point. We
  denote this torsor by $\LHS_{N,\mathbb{L}}$.
\item[{\bfseries (ii)}] Let $\bphi_{N} : \lan{\Lambda} \to P^{D}$ be
the homomorphism which to each $\mu \in \lan{\Lambda} \cong \chr(T)$ 
assigns the
induced $\mathbb{C}^{\times}$ extension $\mu_{*}N$. Regarding
$\mathbb{L}$ as a group extension 
\[
0 \to \mathbb{C}^{\times} \to \mathbb{L} \to P^{D} \to 0 
\]
we get a $Z(G)$-torsor over a point defined as
\[
\RHS_{N,\mathbb{L}} := \left\{ \tilde{\nu} : \lan{\Lambda} \to
\mathbb{L} \; \left| \; \tilde{\nu} \text{ lifts } \nu
\right.\right\},
\]
where $\nu : \pi_{1}(\lan{G}) \to P^{D}$ is the unique homomorphism
for which $\nu\circ \lan{\bq} = \bphi_{N}$.
\end{itemize}

\

Now the construction of $\mathfrak{u}_{\widetilde{C}}$ follows from 

\begin{lem} \label{lem:relations}
The $Z(G)$-torsors $\LHS_{N,\mathbb{L}}$ and $\RHS_{N,\mathbb{L}}$ are
canonically isomorphic.
\end{lem}
{\bfseries Proof.} $\LHS_{N,\mathbb{L}}$ and $\RHS_{N,\mathbb{L}}$ are
the fibers at $\mathbb{L}$ of two torsors $\LHS_{N}$ and $\RHS_{N}$ on
$P$. So to prove the lemma it suffices to identify the corresponding
group extensions canonically. As explained above $\LHS_{N,\mathbb{L}}$
is the fiber at $\mathbb{L}$ of the group extension
\[
0 \to Z(G) \to \sigma^{*}N \to P \to 0,
\]
and so $\LHS_{N} = \sigma^{*}N$.

Next consider the Poincare $\mathbb{C}^{\times}$-torsor $\mycal{P} \to
P\times P^{D}$. From the biextension property it follows that we can
view $\mycal{P}$ as an extension of commutative groups schemes on $P$:
\begin{equation} \label{eq:Poincare.extension}
0 \to \mathbb{G}_{m} \to \mycal{P} \to \underline{P}^{D} \to 0,
\end{equation}
where $\underline{P}^{D}$ denotes the constant group scheme on $P$
with fiber $P^{D}$. At the same time, the extension $N$ gives a
homomorphism $\bphi_{N} : \lan{\Lambda} \to P^{D}$, $\bphi_{N}(\mu) :=
\mu_{*}N$ and the condition that the sequence \ref{eq:N.sequence} is
split implies that $\bphi_{N}$ factors through the quotient $\lan{\bq}
: \lan{\Lambda} \twoheadrightarrow \pi_{1}(\lan{G})$. In other words
there is a unique homomorphism $\nu : \pi_{1}(\lan{G}) \to P^{D}$
satisfying $\bphi_{N} = \nu\circ \lan{\bq}$. Pulling back the Poincare
extension \eqref{eq:Poincare.extension} by $\nu$ gives an extension of
commutative group schemes on $P$:
\begin{equation} \label{eq:nu.Poincare}
0 \to \mathbb{G}_{m} \to \nu^{*}\mycal{P} \to
\underline{\pi_{1}(\lan{G})} \to 0.
\end{equation}
The sequence of group schemes \eqref{eq:nu.Poincare} is split locally
on $P$ but in general is not globally split. The sheaf of local
splittings of \eqref{eq:nu.Poincare}  is representable by a space 
$\RHS_{N}$ which is a principal
$\op{Hom}(\pi_{1}(\lan{G}),\mathbb{C}^{\times}) = 
Z(G)$-bundle on $P$. By definition $\RHS_{N,\mathbb{L}}$ is the fiber of
$\RHS_{N}$ at $\mathbb{L} \in P$. 

To compare $\LHS_{N}$ and $\RHS_{N}$ as $Z(G)$-bundles on $P$ note
first that they are classified by the same element in
$H^{1}(P,Z(G))$. Indeed, the sequence $0 \to T \to N \to P \to 0$
corresponds to an element in $\op{Ext}^{1}(P,T)$. Under the natural
identifications 
\[
\op{Ext}^{1}(P,T) = 
\op{Ext}^{1}(P,\op{Hom}(\lan{\Lambda},\mathbb{C}^{\times}))
\stackrel{(\dagger)}{=} 
\op{Hom}(\lan{\Lambda},\op{Ext}^{1}(P,\mathbb{C}^{\times}))
\stackrel{(\ddagger)}{=} 
\op{Hom}(\lan{\Lambda},P^{D}) 
\]
this element is just $\bphi_{N} \in
\op{Hom}(\lan{\Lambda},P^{D})$. The identification $(\dagger)$ is just
the adjunction isomorphism, and the identification $(\ddagger)$ is the
contraction with the Poincare biextension class
$e_{\mycal{P}}(\bullet,\bullet)$ for $P\times P^{D}$. In other words,
the sequence $0 \to T \to N \to P \to 0$ corresponds to the class
$e_{\mycal{P}}(\bullet,\bphi_{N})$. By the same token the sequence $0
\to Z(G) \to \sigma^{*}N \to P \to 0$ corresponds to the class
$e_{\mycal{P}}(\bullet,\nu) \in \op{Ext}^{1}(P,Z(G))$, i.e. $\LHS_{N}$
is classified by $e_{\mycal{P}}(\bullet,\nu)$.

Similarly, the class of the torsor $\RHS_{N}$ can be computed from its
definition.  The local-to-global spectral sequence identifies the
group $H^{1}(P,\sHom(\underline{\pi_{1}(\lan{G})},\mathbb{G}_{m}))$
with the subgroup of
$\op{Ext}^{1}_{P}(\underline{\pi_{1}(\lan{G})},\mathbb{G}_{m}))$,
consisting of extension classes of locally split extensions. Since
$\RHS_{N}$ is the torsor of local splittings of \eqref{eq:nu.Poincare},
its class in $H^{1}(P,Z(G)) =
H^{1}(P,\sHom(\underline{\pi_{1}(\lan{G})},\mathbb{G}_{m}))$ will be
precisely the extension class of \eqref{eq:nu.Poincare} viewed as an
element in the subgroup
\[
H^{1}(P,\sHom(\underline{\pi_{1}(\lan{G})},\mathbb{G}_{m})) \subset
\op{Ext}^{1}_{P}(\underline{\pi_{1}(\lan{G})},\mathbb{G}_{m})).
\] 
But the class of \eqref{eq:nu.Poincare} is given by
$e_{\mycal{P}}(\bullet,\nu)$, and so $\RHS_{N}$ is also classified by
the element $e_{\mycal{P}}(\bullet,\nu)$. Since $\LHS_{N}$ and
$\RHS_{N}$ are classified by the same element of $H^{1}(P,Z(G))$ it
follows that they are isomorphic as $Z(G)$ torsors.

To exhibit a canonical isomorphism between $\LHS_{N}$ and $\RHS_{N}$
as $Z(G)$-extensions of $P$ it now suffices to identify the fibers
$\LHS_{N}$ and $\RHS_{N}$ at some point of $P$. Let $o \in P$ be
the origin. Then from the descriptions {\bfseries (i)} and {\bfseries
  (ii)} we see that $\LHS_{N,o}$ and $\RHS_{N,o}$ are both canonically
isomorphic to $Z(G)$. Since  $\LHS_{N}$ and $\RHS_{N}$  are isomorphic
as covering spaces of $P$ this isomorphism of fibers extends to a
unique canonical isomorphism of $\LHS_{N}$ and $\RHS_{N}$
as $Z(G)$-extensions of $P$. \ \hfill $\Box$

\

\noindent
Applying the previous lemma to $P = P_{\widetilde{C}}$, and $N =
N((\tilde{x}_{1},\bPhi_{1}),(\tilde{x}_{2},\bPhi_{2}))$ yields the
desired isomorphism $\mathfrak{r}_{\widetilde{C}}$. The intertwining
property of the induced isomorphism of gerbes
$\gHiggs_{0,\widetilde{C}} \stackrel{\cong}{\to} \left( \lan{\Higgs}
\right)^{D}$ follows tautologically from the construction. This
completes the proof of Theorem~\ref{theo:iso.presentations}. \ \hfill
$\Box$

\

\medskip

\begin{rem} \label{rem:another.proof} The proof of
  Lemma~\ref{lem:relations} can be streamlined somewhat and the
  isomorphism $\mathfrak{u}_{\widetilde{C}}$ can be constructed
  without computation of extension classes. As we explained
  above, specifying  $N$ is the same thing as specifying an element
  $\nu \in \op{Hom}(\pi_{1}(\lan{G}),P^{D})$. Varying $N$ we
  see that $\LHS_{N}$ and $\RHS_{N}$ fit together in $Z(G)$-torsors
  $\LHS$ and $\RHS$ on $\op{Hom}(\pi_{1}(\lan{G}),P^{D})\times P$. In
  the proof of  Lemma~\ref{lem:relations} we argued that these torsors
  are canonically isomorphic by first showing that they are
  abstractly isomorphic and then noticing that
  $\LHS_{|\op{Hom}(\pi_{1}(\lan{G}),P^{D})\times \{ o\}}$ and
  $\RHS_{|\op{Hom}(\pi_{1}(\lan{G}),P^{D})\times \{ o\}}$ are
  canonically trivial. 

We can instead show directly that the torsors are canonically
isomorphic by using the see-saw theorem. Indeed if $A$ and $B$ are
connected projective algebraic groups, and if we have two extensions
of $A\times B$ by $Z(G)$ whose restrictions on $\{ o \} \times B$ and
$A\times \{ o \}$ are identified, then the see-saw theorem implies
that the extensions themselves are identified by a unique
isomorphism. Consider $A = \op{Hom}(\pi_{1}(\lan{G}),P^{D})$ and $B =
P$. The restrictions of $\LHS$ and $\RHS$ to the slices
$\op{Hom}(\pi_{1}(\lan{G}),P^{D})\times \{ o\}$ and $\{ o \}\times P$
are identified in an obvious manner from the definitions but
unfortunately the see-saw theorem does not apply immediately since
$\op{Hom}(\pi_{1}(\lan{G}),P^{D})$ is a finite (and hence
disconnected) group. We can remedy this by embedding the finite group
$\op{Hom}(\pi_{1}(\lan{G}),P^{D})$ in an abelian variety on which the
see-saw theorem does apply.  We pull back the extension
\[
0 \to \crts_{\lan{\mathfrak{g}}} \to \lan{\Lambda} \to
\pi_{1}(\lan{G}) \to 0 
\]
via some surjective homomorphism 
$\lan{L} \to \pi_{1}(\lan{G})$:
\begin{equation} \label{eq:pullback.diagram}
\xymatrix@R-1pc@C-1pc{
&  & 0 & 0 \\
0 \ar[r] & \crts_{\lan{\mathfrak{g}}} \ar[r]  \ar@{=}[d] & 
\lan{\Lambda} \ar[r] \ar[u] & \pi_{1}(\lan{G}) \ar[r] \ar[u] & 0 \\
0 \ar[r] & \crts_{\lan{\mathfrak{g}}} \ar[r]  &
\lan{\Lambda}\times_{\pi_{1}(\lan{G})}\lan{\Lambda} \ar[r] \ar[u] & 
\lan{\Lambda} \ar[r] \ar[u] & 0 \\
& &  \crts_{\lan{\mathfrak{g}}} \ar@{=}[r] \ar[u] &
\crts_{\lan{\mathfrak{g}}} \ar[u] & \\
& & 0 \ar[u] & 0 \ar[u] & 
}
\end{equation}
In the same way we defined $\LHS$ and $\RHS$ we can use the second row
of this diagram to define two
$\op{Hom}(\lan{\Lambda},\mathbb{C}^{\times}) = T$-torsors
$\widetilde{\LHS}$ and $\widetilde{\RHS}$ on
$\op{Hom}(\lan{\Lambda},P^{D})\times P$. Note that
$\op{Hom}(\pi_{1}(\lan{G}),P^{D}) \subset
\op{Hom}(\lan{\Lambda},P^{D})$ and by construction we have
tautological identifications 
\begin{equation} \label{eq:restrict.induce}
\begin{split}
\widetilde{\LHS}_{|\op{Hom}(\pi_{1}(\lan{G}),P^{D})\times P} & =
\LHS\times_{Z(G)} T, \\
\widetilde{\RHS}_{|\op{Hom}(\pi_{1}(\lan{G}),P^{D})\times P} & =
\RHS\times_{Z(G)} T.
\end{split}
\end{equation}
In particular the restrictions of $\widetilde{\LHS}/Z$ and
$\widetilde{\RHS}/Z$ to $\op{Hom}(\pi_{1}(\lan{G}),P^{D})\times P$ are
naturally trivialized, hence naturally isomorphic. Thus the
restriction of the  $T/Z$
torsor of isomorphisms $\sIsom(\widetilde{\LHS}/Z,\widetilde{\RHS}/Z)$ to
$\op{Hom}(\pi_{1}(\lan{G}),P^{D})\times P$ admits a canonical
trivialization $\mathfrak{triv}_{0}$.

Now $\op{Hom}(\lan{\Lambda},P^{D})$ is connected and the see-saw
argument provides a canonical identification $\widetilde{\LHS} \cong
\widetilde{\RHS}$ and hence a canonical trivialization
$\mathfrak{triv}$ of the $T$-torsor of isomorphisms
$\sIsom(\widetilde{\LHS},\widetilde{\RHS})$. From the identifications
\eqref{eq:restrict.induce} it then follows that the $Z(G)$-torsor of
isomorphisms $ \sIsom(\LHS,\RHS)$ is naturally identified with the
torsor of all trivializations of $\sIsom(\widetilde{\LHS},
\widetilde{\RHS})_{|\op{Hom}(\pi_{1}(\lan{G}),P^{D})\times P}$ that
lift the trivialization $\mathfrak{triv}_{0}$ of the torsor
$\sIsom(\widetilde{\LHS}/Z,\widetilde{\RHS}/Z)_{|\op{Hom}(
\pi_{1}(\lan{G}),P^{D})\times P}$. But diagram
\eqref{eq:restrict.induce} implies $\mathfrak{triv}$ and
$\mathfrak{triv}_{0}$ are compatible trivializations, and so
$\mathfrak{triv}$ naturally trivializes $\sIsom(\widetilde{\LHS},
\widetilde{\RHS})_{|\op{Hom}(\pi_{1}(\lan{G}),P^{D})\times P}$ as a
$Z(G)$-torsor over
$\sIsom(\widetilde{\LHS}/Z,\widetilde{\RHS}/Z)_{|\op{Hom}(
\pi_{1}(\lan{G}),P^{D})\times P}$.
\end{rem}

\section{Extensions and refinements}

In this section we extend Theorem\ref{thm:gerbes} from simple groups to all
reductive groups.  
We also discuss some additional structures related to the weight filtration. 
Finally, we draw the main geometric corollaries of the duality: 
existence of a Fourier-Mukai equivalence, 
and the construction of Hecke eigensheaves.

\subsection{Extension to reductive groups} \label{ss:reductive}

\noindent
Our main duality result extends to general reductive
groups.

\begin{thm} \label{thm:reductive} Let $\mathbb{G}$ be a connected complex
reductive group, let $\lan{\mathbb{G}}$ be the Langlands dual
reductive group, and let $C$ be a smooth compact complex curve. Write
$\gHiggs_{\mathbb{G}}$ and $\gHiggs_{(\lan{\mathbb{G}})}$ for the
stacks of $K_{C}$-valued Higgs bundles on $C$ with structure group
$\mathbb{G}$ and $\lan{\mathbb{G}}$ respectively. Then there is an
isomorphism $\bl : B \; \widetilde{\to} \; \lan{B}$ of the respective
Hitchin bases which gives an identification $B - \Delta \cong \lan{B}
- \lan{\Delta}$. Under this identification one has an
isomorphism 
\[ 
\gHiggs_{\mathbb{G}}  \cong 
 \left( \gHiggs_{(\lan{\mathbb{G}})} \right)^{D}
\] 
of commutative group stacks over $B - \Delta$, intertwining the action
of translation and tensorization operators.
\end{thm}

\

\noindent
{\bfseries Remark} \ Since in this proof we need to work with the
moduli stacks of Higgs bundles for different isogenous groups we will
deviate from our normal notation and will label the corresponding
stacks with the group as a subscript. Thus we will write
$\gHiggs_{\mathbb{G}}$, $\gHiggs_{\mathbb{\lan{G}}}$, etc..

\

\noindent
{\bf Proof. } The proof is somewhat involved. We start by presenting
what seems to us to be a very natural homological approach. However,
this runs into technical issues, which we explain in
Remark~\ref{rem:tack}.  Rather than
settling these issues, we find it easier to change tack and give a
separate argument in subsections~\ref{sss:red.Prym}, \ref{sss:pi0},
\ref{sss:split.sequences}. This argument is a modification of the
proof of Theorem~\ref{thm:gerbes}.

Since $\mathbb{G}$ is connected and reductive, we can
always fit $\mathbb{G}$ in a short exact sequence 
\begin{equation} \label{eq:seqG}
\xymatrix@1{
1 \ar[r] & K \ar[r] & G\times H \ar[r] & \mathbb{G} \ar[r] & 1,
}
\end{equation}
where $K$ is a finite subgroup in the center $Z(G\times H)$ of
$G\times H$, $G = \prod_{i =1}^{a}G_{i}$ is a
product of complex simple groups, and $H \cong
(\mathbb{C}^{\times})^{b}$ is an affine complex torus. Passing to
Langlands duals gives the sequence
\begin{equation} \label{eq:seqLG}
\xymatrix@1{
1 \ar[r] & K^{\wedge} \ar[r] & \lan{\mathbb{G}} \ar[r] & \lan{G}
\times \lan{H} \ar[r] & 1.
}
\end{equation}
Next observe that the construction of the Hitchin base, the formation
of the moduli stack of Higgs bundles, the definition of the sheaf
$\mathcal{T}$, as well as the operations $\lan{(\bullet)}$ and
$(\bullet)^{D}$, all respect the operation of taking products of
groups. Combined with Theorems~\ref{thm:duality} and \ref{thm:gerbes},
and with the standard selfduality of $\gPic$ of a smooth curve, we get
an identification of the Hitchin bases for $G\times H$ and $\lan{G}
\times \lan{H}$, as well as a global duality
\begin{equation} \label{eq:product}
\left(\gHiggs_{\lan{G} 
\times \lan{H}}\right)^{D} \cong \gHiggs_{G\times H}.
\end{equation}
Also, since the Hitchin base depends only on the Lie algebra, and not
on the Lie group, it follows that the identification of the Hitchin
bases for $G\times H$ and $\lan{G}
 \times \lan{H}$ can be interpreted as the desired isomorphism $B
 \cong \lan{B}$. 

Furthermore, the definition of $\mathcal{T}$ (see \eqref{eq:sheafT})
gives short exact sequences of abelian sheaves on $B\times C$:
\begin{equation} \label{eq:TbbG}
\xymatrix@R-1pc{
0 \ar[r] & K \ar[r] & \mathcal{T}_{G}\oplus H \ar[r] &
\mathcal{T}_{\mathbb{G}} \ar[r] & 0  \\
0 \ar[r] & K^{\wedge} \ar[r] & \mathcal{T}_{\lan{\mathbb{G}}} \ar[r] &
\mathcal{T}_{\lan{G}}\oplus \lan{H} \ar[r] & 0. 
}
\end{equation}
Applying $R\pi_{*}[1]$ to the first sequence we get a distinguished
triangle in $D^{b}(B)$:
\begin{equation} \label{eq:firsttr}
\xymatrix@1{
R\pi_{*}K[1] \ar[r]^-{g} & R\pi_{*}\mathcal{T}_{G}[1]\oplus
R\pi_{*}H[1] \ar[r] & 
R\pi_{*}\mathcal{T}_{\mathbb{G}}[1] \ar[r] & R\pi_{*}K[2].
}
\end{equation}
Applying $\underline{R\op{Hom}}(R\pi_{*}(\bullet),\mathcal{O}^{\times})$
to the second sequence, and taking into account the isomorphism
\eqref{eq:product} and the isomorphism $\gHiggs_{\mathbb{G}} \cong
R\pi_{*}\mathcal{T}_{\mathbb{G}}[1]$ from Lemma~\ref{lem:section}, we
get another 
distinguished triangle  
\begin{equation} \label{eq:sectr}
\xymatrix@1@C-0.5pc{
\underline{R\op{Hom}}(R\pi_{*}K^{\wedge}[1],\mathcal{O}^{\times})
\ar[r]^-{\lan{g}} &  
  R\pi_{*}\mathcal{T}_{G}[1]\oplus R\pi_{*}H[1] \ar[r] &
  (R\pi_{*}\mathcal{T}_{\lan{\mathbb{G}}}[1])^{D} \ar[r] &
  \underline{R\op{Hom}}(R\pi_{*}K^{\wedge},\mathcal{O}^{\times}) .
}
\end{equation}
We wish to show that these two triangles are isomorphic and so we have a
quasi-isomorphism $R\pi_{*}\mathcal{T}_{\mathbb{G}}[1] \cong
(R\pi_{*}\mathcal{T}_{\lan{\mathbb{G}}}[1])^{D}$. 
Now Poincare duality on a smooth curve $C$ identifies the
cohomology of  $C$ with coefficients in $K$ with the
Pontryagin dual of the cohomology of $C$ with coefficients
in $K^{\wedge}$. In particular, it induces an isomorphism of complexes 
\[
\op{PD} :
R\pi_{*}K[1]  \stackrel{\cong}{\to}
\underline{R\op{Hom}}(R\pi_{*}K^{\wedge}[1],\mathcal{O}^{\times}) 
\]
and what we have to show is that it intertwines the maps $g$ and
$\lan{g}$. So, we must show that the 
diagram:
\begin{equation} \label{eq:PD.diagram}
\xymatrix@R-1pc{
R\pi_{*}K[1] \ar[rd]^-{g} \ar[dd]_-{\op{PD}} & \\
&  R\pi_{*}\mathcal{T}_{G}[1]\oplus
R\pi_{*}H[1] \\
\underline{R\op{Hom}}(R\pi_{*}K^{\wedge}[1],\mathcal{O}^{\times}) 
\ar[ru]_-{\lan{g}}  &
}
\end{equation}
commutes in the derived category.

Let $K_{G}$, $K_{H}$ be the projections to $G$, $H$ respectively of
(the image in $G\times H$ of) $K$. These are central subgroups. The
short exact sequence \eqref{eq:seqG} extends to a commutative diagram
with exact rows:
\[
\xymatrix@R-1pc{
1 \ar[r] & K \ar[r] \ar[d] & G\times H \ar[r] \ar@{=}[d] & \mathbb{G}
\ar[r] \ar[d] & 1 \\
1 \ar[r] & K_{G}\times K_{H} \ar[r] & G\times H \ar[r]  &
\overline{G}\times \overline{H}
\ar[r] & 1.
}
\]
Now diagram \eqref{eq:PD.diagram} factors:
\[
\xymatrix@R-1pc@C-1pc{
R\pi_{*}K[1] \ar[r] \ar[dd]_-{\op{PD}} & R\pi_{*}K_{G}[1]\oplus
R\pi_{*}K_{H}[1] 
\ar[rd] \ar[dd]_-{\op{PD}} & \\
&  & R\pi_{*}\mathcal{T}_{G}[1]\oplus
R\pi_{*}H[1] \\
\underline{R\op{Hom}}(R\pi_{*}K^{\wedge}[1],\mathcal{O}^{\times})  
 \ar[r] & 
\underline{R\op{Hom}}(R\pi_{*}(K_{G}^{\wedge}\oplus
K_{H}^{\wedge})[1],\mathcal{O}^{\times})
\ar[ru]  &
}
\]
Commutativity of the square follows from functoriality of Poincare
duality. The triangle part of this diagram is just
\eqref{eq:PD.diagram} but for the subgroup $K_{G}\times K_{H}$ instead
of $K$. Commutativity for the $G$-factor of the triangle is discussed
below.  Commutativity for the torus $H$ is obvious.

To understand better the $G$-factor of the triangle we have to examine
more closely the relationship between Poincare duality on $C$ and the
canonical isomorphism $\gerbeiso : \gHiggs_{G} \to
(\gHiggs_{\lan{G}})^{D}$ from Theorem~\ref{thm:gerbes}.  
Let us examine this commutativity statement for a fixed cameral
cover $\widetilde{C} \in B - \Delta$. Since $K_{G}
\subset Z(G) \subset G$, we again have short exact sequences of
sheaves of abelian groups on $C$
\begin{equation} \label{eq:KGses}
\xymatrix@1{
1 \ar[r] & K_{G} \ar[r] & \mathcal{T}_{G} \ar[r] &
\mathcal{T}_{G/K}  \ar[r] & 1,
}
\end{equation}
and 
\begin{equation} \label{eq:KGLses}
\xymatrix@1{
1 \ar[r] & {K_{G}}^{\wedge} \ar[r] & \mathcal{T}_{\lan{(G/K)}} \ar[r] &
\mathcal{T}_{\lan{G}}  \ar[r] & 1,
}
\end{equation}
where now $\mathcal{T}_{G}$, $\mathcal{T}_{G/K}$,
$\mathcal{T}_{\lan{(G/K)}}$, $\mathcal{T}_{\lan{G}}$ are the sheaves
corresponding to the cover $\widetilde{C} \to C$. 

Passing to cohomology in \eqref{eq:KGses} yields a map
\[
g_{G} : R\Gamma(C,K_{G})[1] \to R\Gamma(C,\mathcal{T}_{G})[1].
\]
Similarly, passing to cohomology in \eqref{eq:KGLses} yields a map
\[
\partial_{G} : R\Gamma(C,\mathcal{T}_{\lan{G}}) \to
R\Gamma(C,{K_{G}}^{\wedge})[1]. 
\]
Dualize $\partial_{G}$ by applying
$R\op{Hom}(\bullet,\mathbb{C}^{\times})$ to obtain a map 
\[
{\partial_{G}}^{\wedge} :
R\op{Hom}(R\Gamma(C,{K_{G}}^{\wedge})[1],\mathbb{C}^{\times}) \to 
R\op{Hom}(R\Gamma(C,\mathcal{T}_{\lan{G}}),\mathbb{C}^{\times}). 
\]
With this notation the commutativity of the $G$-part of the triangle
(for a fixed cameral cover $\widetilde{C}$) becomes the
statement that 
the diagram  
\begin{equation} \label{eq:G.factor}
\xymatrix{
R\Gamma(C,K_{G})[1] \ar[d]_-{\op{PD}} \ar[r]^-{g_{G}} &
R\Gamma(C,\mathcal{T}_{G})[1] & \hspace{-0.6in} = &
\hspace{-0.6in}\gHiggs_{G,\widetilde{C}}  
\ar@<-4ex>[d]^-{\gerbeiso_{\widetilde{C}}} \\
R\op{Hom}(R\Gamma(C,{K_{G}}^{\wedge})[1],\mathbb{C}^{\times})
\ar[r]^-{{\partial_{G}}^{\wedge}} &
R\op{Hom}(R\Gamma(C,\mathcal{T}_{\lan{G}}),\mathbb{C}^{\times})
& \hspace{-0.6in} = & 
\hspace{-0.6in} \left(\gHiggs_{\lan{G},\widetilde{C}}\right)^{D}
}
\end{equation}
commutes in the derived category of abelian groups. Note also that the
commutativity of \eqref{eq:PD.diagram} on the level of cohomology is a
statement about sheaves and so can be checked locally on the Hitchin
base. In other words the commutativity of \eqref{eq:G.factor} on the
level of cohomology implies the commutativity of \eqref{eq:PD.diagram}
on the level of cohomology. (This will be relevant for the discussion
in subsection~\ref{sss:red.Prym}.) 

\

\noindent
In summary, if we could establish the commutativity of
\eqref{eq:G.factor} and  \eqref{eq:PD.diagram}, this would imply
that the distinguished triangles \eqref{eq:firsttr} and \eqref{eq:sectr} are
isomorphic and so we would have a quasi-isomorphism
\begin{equation} \label{eq:thmC.complex.iso}
(R\pi_{*}\mathcal{T}_{\lan{\mathbb{G}}}[1])^{D} \cong
R\pi_{*}\mathcal{T}_{\mathbb{G}}[1]. 
\end{equation}
Passing to the associated stacks we would obtain the statement of
Theorem~\ref{thm:reductive}. 

\

\medskip

\begin{rem} \label{rem:tack}
Despite its streamlined form and homological appeal the above argument
has a couple of drawbacks. First, the duality isomorphism
\eqref{eq:thmC.complex.iso} is obtained by using the cone-filling
axiom in triangulated categories and so it is unclear whether this
construction can be made canonical.  Second, the whole construction is
based on showing that \eqref{eq:G.factor} and \eqref{eq:PD.diagram}
commute. The commutativity of \eqref{eq:G.factor} and
\eqref{eq:PD.diagram} is hard to verify, mainly because our duality
isomorphism $\gerbeiso : \gHiggs_{G} \to
\left(\gHiggs_{\lan{G}}\right)^{D}$ was constructed geometrically
rather than homologically.

\

\noindent
Instead, we will follow a slightly different approach, which is a
combination of the homological and geometric arguments.

\begin{enumerate}
\item[(1)] We will use the above ideas to analyze the connected
components of $\gHiggs_{\mathbb{G}}$ and
$(\gHiggs_{\lan{\mathbb{G}}})^{D}$ and to extend the duality of Prym
varieties $P_{G,\widetilde{C}}$ and $P_{\lan{G},\widetilde{C}}$ to a
canonical duality between the Prym varieties
$P_{\mathbb{G},\widetilde{C}}$ and
$P_{\lan{\mathbb{G}},\widetilde{C}}$. We will see that this boils down
to the case of a simple group, treated in Theorem~\ref{thm:duality},
plus the commutativity of \eqref{eq:G.factor} and of
\eqref{eq:PD.diagram} on the level of cohomology, which is much easier
than the commutativity in the derived category.
\item[(2)] The proof of Theorem~\ref{thm:reductive} (or
Theorem~\ref{thm:gerbes} for reductive groups) is reduced to the
connected case, via a computation showing that: $\pi_0(\gHiggs_{\mathbb{G}}) =
\pi_0(\Higgs_{\mathbb{G}}) = \pi_1(\mathbb{G})$.
\item[(3)] The isomorphism for stacky connected components is almost
identical to the corresponding part of Theorem~\ref{thm:gerbes},
except that the splitting of some sequences needs to be checked more
carefully, as the quotient groups involved are finitely generated rather than
finite.
\end{enumerate}

\

\noindent
We complete the argument over the next three subsections.
\end{rem}

\

\medskip

\noindent
\punkt \label{sss:red.Prym} \ First observe that that in the long
exact sequences of cohomology associated with the distinguished
triangles \eqref{eq:firsttr} and \eqref{eq:sectr} 
the groups $H^{1}(C,K)$ and
$H^{1}(C,K^{\wedge})^{\wedge}$ embed in the connected
components of the groups $\Higgs_{G\times H, \widetilde{C}}$ and
$\left(\left(\gHiggs_{\lan{G}\times \lan{H},
  \widetilde{C}}\right)_{0}\right)^{D}$. This gives 
canonical identifications
\begin{equation} \label{eq:Prym.quotients}
\begin{split}
P_{\mathbb{G},\widetilde{C}}  & \cong \left.\left(P_{G,\widetilde{C}} \times
H^{1}(C,H)_{0}\right)\right/H^{1}(C,K) \\[+1pc]
(P_{\lan{\mathbb{G}},\widetilde{C}})^{D} & \cong
\left.\left((P_{\lan{G},\widetilde{C}})^{D} \times
(H^{1}(C,\lan{H})_{0})^{D}\right)\right/H^{1}(C,K^{\wedge})^{\wedge},
\end{split}
\end{equation}
where $H^{1}(C,H)_{0}$ and $H^{1}(C,\lan{H})_{0}$ denote the
connected components of $H^{1}(C,H)$ and $H^{1}(C,\lan{H})$ respectively.

Now, in order to show that the isomorphisms of abelian varieties 
\[
\xymatrix@R-3pc@M+1pc{
\paiso_{G,\widetilde{C}} : \hspace{-5pc} &  
P_{G,\widetilde{C}} \ar[r]^-{\cong} & (P_{\lan{G},\widetilde{C}})^{D}
\\
\paiso_{H} : \hspace{-5pc} &   (H^{1}(C,H))_{0} \ar[r]^-{\cong} &
(H^{1}(C,\lan{H})_{0})^{D} 
}
\] 
induce a canonical isomoprhism of abelian varieties
\[
\xymatrix@1@M+1pc{ \paiso_{\mathbb{G},\widetilde{C}} : \hspace{-5pc} &
P_{\mathbb{G},\widetilde{C}} \ar[r]^-{\cong} &
(P_{\lan{\mathbb{G}},\widetilde{C}})^{D} }
\] 
it suffices to check the commutativity of the diagram one obtains from
\eqref{eq:PD.diagram} after passing to degree 0 cohomology. This is
much more manageable. In fact we have the following slightly
stronger statement:

\begin{lem} \label{lem:G.factor} The diagrams  \eqref{eq:G.factor} and
  \eqref{eq:PD.diagram} 
  commute on the level of cohomology. 
\end{lem}
{\bf Proof.} As we noted above the commutativity of
\eqref{eq:PD.diagram} on the level of cohomology is a statement about
sheaves. This can be checked locally on the Hitchin base and so it is
enough to check the statement for \eqref{eq:G.factor} on the level of
cohomology. The complexes in the left column of \eqref{eq:G.factor}
are concentrated in degrees $(-1)$, $0$, and $1$, while the complexes
in the right column are concentrated in degrees $(-1)$ and
$0$. Computing the cohomologies in degree $(-1)$ we get:
\[
\boxed{H^{-1}:} \qquad
\xymatrix{
H^{0}(C,K_{G}) \ar[d]_-{\op{PD}} \ar[r]^-{h^{-1}(g_{G})} &
H^{0}(C,\mathcal{T}_{G}) \ar@{=}[r] & Z(G) 
\ar[ddr]^-{h^{-1}(\gerbeiso_{\widetilde{C}})}  & \\
H^{2}(C,{K_{G}}^{\wedge})^{\wedge} 
\ar[r]^-{h^{-1}({\partial_{G}}^{\wedge})} &
\op{Hom}(H^{1}(C,\mathcal{T}_{\lan{G}}),\mathbb{C}^{\times})
\ar@{=}[d] & & \\
&
\op{Hom}(\pi_{0}(H^{1}(C,\mathcal{T}_{\lan{G}})),\mathbb{C}^{\times})
\ar@{=}[r] & 
\pi_{1}(\lan{G})^{\wedge} \ar@{=}[r] & 
Z(G)   
}
\]
Since $C$ is a smooth compact curve we have canonical identifications
$H^{0}(C,K_{G}) = K_{G}$ and $H^{2}(C,{K_{G}}^{\wedge})^{\wedge}  =
K_{G}^{\wedge\wedge} = K_{G}$, and so the diagram of cohomologies in degree
$(-1)$ becomes 
\[
\xymatrix{ K_{G} \ar[r] \ar[d]_-{\op{id}} & Z(G)
\ar[d]^-{h^{-1}(\gerbeiso_{\widetilde{C}})} \\ K_{G} \ar[r] & Z(G)
 }
\]
The two horizontal arrows here are simply the inclusion of $K_{G}$ in
$Z(G)$ and so to show that this diagram commutes we only need to check
that $h^{-1}(\gerbeiso_{\widetilde{C}})$ induces the identity when
restricted to $K_{G}$. This however is automatic since the map of
presentations in  Theorem~\ref{theo:iso.presentations} is
$Z(G)$-equivariant.

Similarly, the diagram of  cohomologies in degree $0$ reads: 
\[
\boxed{H^{0}:} \qquad  
\xymatrix@M+0.5pc{
H^{1}(C,K_{G}) \ar[r]^-{h^{0}(g_{G})} \ar[d]_-{\op{PD}} &
H^{1}(C,\mathcal{T}_{G}) & \hspace{-0.6in} = 
& \hspace{-0.6in} \Higgs_{G,\widetilde{C}} 
\ar@<-5ex>[d]^-{h^{0}(\gerbeiso_{\widetilde{C}})} \\
H^{1}(C,{K_{G}}^{\wedge})^{\wedge} \ar[r]^-{h^{0}({\partial_{G}}^{\wedge})}  &
R^{0}\op{Hom}(R\Gamma(C,\mathcal{T}_{\lan{G}}),\mathbb{C}^{\times}) &
\hspace{-0.6in} = & 
\hspace{-0.6in} 
\left(\left(\gHiggs_{\lan{G},\widetilde{C}}\right)_{0}\right)^{D}
}
\]
The top horizontal map is the map of inducing a $\mathcal{T}_{G}$
torsor from a $K_{G}$ torsor via the inclusion $K_{G} \subset
\mathcal{T}_{G}$ and thus lands in the connected component of the
identity $P_{\widetilde{C}}$ of the abelian group
$\Higgs_{G,\widetilde{C}}$. Similarly the bottom horizontal map is the
dual of the natural map
$\left(\gHiggs_{\lan{G},\widetilde{C}}\right)_{0} \to \lan{P} \to
H^{1}(C,{K_{G}}^{\wedge})[1]$ and thus also factors through the
connected component of the identity $\lan{P}_{\widetilde{C}}^{D}$ of
$\left(\left(\gHiggs_{\lan{G},\widetilde{C}}\right)_{0}\right)^{D}$.
The right vertical map is the map of commutative group schemes induced from
$\gerbeiso_{\widetilde{C}}$, i.e. coincides with the map
$\paiso_{\widetilde{C}}$ from Theorem~\ref{thm:duality}. In other
words the commutativity of the diagram of degree $0$ cohomologies
reduces to showing that the diagram
\begin{equation} \label{eq:KG.av}
\xymatrix{
H^{1}(C,K_{G}) \ar[r] \ar[d]_-{\op{PD}} & P_{\widetilde{C}}
\ar[d]^-{\paiso_{G,\widetilde{C}}} \\
H^{1}(C,{K_{G}}^{\wedge})^{\wedge} \ar[r] & \lan{P}_{\widetilde{C}}^{D}
}
\end{equation}
commutes.

Passing to cocharacter lattices for the abelian varieties
$P_{\widetilde{C}}$ and $\lan{P}_{\widetilde{C}}^{D}$ we see that the
commutativity of \eqref{eq:KG.av} is equivalent to the commutativity
of the following diagram of finitely generated abelian groups:
\begin{equation} \label{eq:KG.lattices}
\xymatrix{
\cchr(P_{\widetilde{C}})^{\vee} 
\ar[r] \ar[d] & H^{1}(C,K_{G})^{\wedge} \ar[d]^-{\op{PD}}  \\
\cchr(\lan{P}_{\widetilde{C}}) \ar[r] 
& H^{1}(C,K_{G})^{\wedge}
}
\end{equation}
However in the proof of Theorem~\ref{thm:duality} we constructed the
isomorphism $\cchr(P_{\widetilde{C}})^{\vee}  \to
\cchr(\lan{P}_{\widetilde{C}})$ by sandwiching both groups between the
two isogeneous lattices
$H^{1}(C,\jmath_{*}A^{\vee}) \subset
H^{1}(C,\jmath_{*}A)_{\op{tf}}^{\vee}$. Furthermore Lemma~\ref{lem:Hij},
Corollary~\ref{cor:image}, and Claim~\ref{claim:cochar} imply that the
sandwiching maps commute with the Poincare duality isomorphisms for
the cohomologies of $A$ and $A^{\vee}$ on $U$. This shows that
\eqref{eq:KG.lattices} commutes. 

Finally, the diagram of  cohomologies in degree $1$ reads
\[
\boxed{H^{1}:} \qquad \xymatrix{
H^{2}(C,K_{G}) \ar[r] \ar[d]_-{\op{PD}} & 0 \ar[d] \\
H^{0}(C,{K_{G}}^{\wedge})^{\wedge} \ar[r] & 0
}
\]
and so is automatically commutative. 

This completes the proof of the lemma. \ \hfill $\Box$

\

\noindent
It is immediate from the lemma that the diagram 
\[
\xymatrix{
H^{1}(C,K) \ar@{^{(}->}[r] \ar_-{\op{PD}}[d] & P_{G,\widetilde{C}}\times
H^{1}(C,H)_{0} \ar^-{\paiso_{G,\widetilde{C}}\times \paiso_{H}}[d] \\
H^{1}(C,K^{\wedge})^{\wedge} \ar@{^{(}->}[r]  & P_{G,\widetilde{C}}^{D}
\times (H^{1}(C,\lan{H})_{0})^{D} 
}
\]
commutes, and so the isomorphism $\paiso_{G,\widetilde{C}}\times
\paiso_{H}$ induces a natural isomorphism of the quotients for the top
and bottom inclusions. By \eqref{eq:Prym.quotients} this gives the
desired isomorphism:
\[
\xymatrix@1@M+1pc{ \paiso_{\mathbb{G},\widetilde{C}} : \hspace{-5pc} &
P_{\mathbb{G},\widetilde{C}} \ar[r]^-{\cong} &
(P_{\lan{\mathbb{G}},\widetilde{C}})^{D}. 
}
\]

\

\medskip
\punkt \label{sss:pi0} \
Having the isomorphism $\paiso_{\mathbb{G},\widetilde{C}}$ of abelian varieties
 at our disposal we can proceed to
construct the isomorphism of Higgs stacks as in
Section~\ref{ssec:proofThmB}. The construction of the presentations of
the two stacks and the isomorphism of presentations follows almost
verbatim the constructions in the case of a simple $G$ described in
Section~\ref{ssec:proofThmB}. We will not repeat the steps in the
proof of Theorem~\ref{thm:gerbes} but will only indicate the changes
necessary to make the arguments work for a general reductive $\mathbb{G}$.

First note that the abelianized Hecke and
tensorization  
operators were constructed in Section~\ref{ssec:AJ} and
Section~\ref{ssec:dualHecke} by using the Abel-Jacobi map and the
abelianization procedure from \cite{ron-dennis} which works for
arbitrary reductive groups. In particular the reasoning in
Section~\ref{ssec:AJ} and Section~\ref{ssec:dualHecke} applies
directly to the case of reductive groups and again gives Hecke
and tensorization automorphisms 
\[
\xymatrix@R-2pc{ \trans^{\lambda,\tilde{x}} : &
  \gHiggs_{\mathbb{G},\widetilde{C}} \ar[r]^-{\cong} &
  \gHiggs_{\mathbb{G},\widetilde{C}} \\ \tens^{\lambda,\tilde{x}} : &
  \left(\gHiggs_{\lan{\mathbb{G}},\widetilde{C}}\right)^{D}
  \ar[r]^-{\cong} &
  \left(\gHiggs_{\lan{\mathbb{G}},\widetilde{C}}\right)^{D} }
\]
labeled by characters $\lambda$ of the maximal torus of $\mathbb{G}$
and points $\tilde{x}$ of $\widetilde{C}$. This will be done in
Section~\ref{sss:split.sequences}. 

In particular, if we can check that the group of connected components
of $\gHiggs_{\mathbb{G}}$ (or $\Higgs_{\mathbb{G}}$) is 
naturally isomorphich to $\pi_{1}(\mathbb{G})$, we can proceed as in
the beginning of Section~\ref{ssec:proofThmB} and reduce the problem
of constructing the duality between $\gHiggs_{\mathbb{G}}$ and
$\gHiggs_{\lan{\mathbb{G}}}$ to constructing an isomoprphism of
connected $Z(\mathbb{G})$-gerbes
\begin{equation} \label{eq:bbG0iso}
\gHiggs_{\mathbb{G},0} \cong  \left(  \Higgs_{\lan{\mathbb{G}}} \right)^{D}
\end{equation}
which intertwines the action of the Hecke and tensorization operators
labeled by $(\lambda,\tilde{x}) \in \crts(\op{Lie}(\mathbb{G}))\times
\widetilde{\mycal{C}}$. 

To show that we indeed have a natural identification
$\pi_{0}(\Higgs_{\mathbb{G}}) = \pi_{1}(\mathbb{G})$ we will again
use the homological desctiption of the Higgs moduli space.  The short exact
sequence of groups \eqref{eq:seqG} induces a short exact sequence of
fundamental groups:
\begin{equation} \label{eq:fg}
\xymatrix@1{
0 \ar[r] & \pi_{1}(G\times H)  \ar[r] &  \ar[r] \pi_{1}(\mathbb{G}) & K \ar[r]
& 0.
}
\end{equation}
On the other hand if we take the first short exact sequence of abelian
sheaves in \eqref{eq:TbbG} and pass to cohomology, then we will get a
(piece of a) long exact sequence 
\begin{equation} \label{eq:lesbbG}
\xymatrix@R-1.5pc@C-1pc{
0 \ar[r] & H^{1}(C,K) \ar[r] & H^{1}(C,\mathcal{T}_{G}\oplus H) \ar[r]
\ar@{=}[d] &
H^{1}(C,\mathcal{T}_{\mathbb{G}}) \ar[r] \ar@{=}[d] & H^{2}(C,K)
\ar[r] & 0. \\
& & \Higgs_{G\times H,\widetilde{C}} &
\Higgs_{\mathbb{G},\widetilde{C}} & & 
}
\end{equation}
Here the injectivity at $H^1(C,K)$ follows because the long exact
sequence in cohomology begins as 
\[
0 \to K \to Z(G)\times H \to Z(\mathbb{G}) \to H^1(K) \to \cdots 
\]
and because the map  $Z(G)\times H \to Z(\mathbb{G})$ is surjective.
Similarly the surjectivity at $H^2(C,K)$ follows because the sheaf
$\mathcal{T}_{G}\oplus H$ has cohomologies only in degrees 0 and
1. For $\mathcal{T}_{G}$ this is proven in Section~\ref{s:duality} and
for $H$ this is obvious since $H$ is a product of several copies of
$\mathbb{G}_{m}$.

Separating the groups in sequence \eqref{eq:lesbbG} into their
connected and disconnected parts gives a commutative diagram with
exact rows and columns:
\begin{equation} \label{eq:pi0bbG}
\xymatrix@R-1pc@C-1pc{
&  0 \ar[d] & 0 \ar[d] & 0 \ar[d] & & \\ 
0 \ar[r] & H^{1}(C,K) \ar[r] \ar[d] & P_{G,\widetilde{C}}\times
H^{1}(C,H)_{0} \ar[r] \ar[d] &  P_{\mathbb{G},\widetilde{C}} \ar[d] \ar[r] &
0 \ar[d] & \\ 
0 \ar[r] & H^{1}(C,K) \ar[r] \ar[d] & H^{1}(C,\mathcal{T}_{G}\oplus H) \ar[r]
\ar[d] &
H^{1}(C,\mathcal{T}_{\mathbb{G}}) \ar[r] \ar[d] & H^{2}(C,K)
\ar[d] 
\ar[r] & 0 \\
& 0 \ar[r] & \pi_{1}(G\times H) \ar[r] \ar[d] &
\pi_{0}(\Higgs_{\mathbb{G},\widetilde{C}}) \ar[r]  \ar[d] & K \ar[r]
\ar[d] & 0 \\
& & 0 & 0 & 0 &   
}
\end{equation}
We have already discussed the top row of this diagram - it gives the
top isomorphism in \eqref{eq:Prym.quotients}. For the present
discussion the important part is the bottom row. The first Chern class
gives a natural map
from the bottom row of \eqref{eq:pi0bbG} to the sequence
\eqref{eq:fg}. This map  
is the identity on the two outside groups, and so the map from
$\pi_{0}(\Higgs_{\mathbb{G},\widetilde{C}})$ to $\pi_{1}(\mathbb{G})$ 
must be an isomorphism.

\

\punkt \label{sss:split.sequences} 
 \ Next we can proceed to construct the isomorphism \eqref{eq:bbG0iso} by
following the steps in Section~\ref{ssec:proofThmB}: construct atlases
for the two gerbes and then show that the corresponding groupoid
presentatins are isomorphic. This works without any modifications but
the details of the geometry are slightly altered by the fact that
$\pi_{1}(\mathbb{G})$ is no longer finite but rather is a finitely
generated abelian group. Again the atlases are constructed as moduli
spaces of splittings of extensions but in the reductive case these are
extensions of finitely generated abelian groups by abelian
varieties. Such extensions are again split since the finitely
generated abelian groups are products of free abelian groups and
finite abelian groups and the abelian varieties are divisible
groups. This fact allows us to carry out the
Section~\ref{ssec:proofThmB} constructions of atlases and relations
for the Higgs stacks in the reductive case. The resulting atlases are
no longer finite Galois covers of the cameral curves but rather are
families of abelian varieties over finite Galois covers. This does not
affect the rest of the arguments in any way. Finally for the
construction of the isomorphism of groupoid presentations we note that
the isomorphism argument given in Remark~\ref{rem:another.proof}
applies directly to 
the reductive setting. 

\

\noindent
This completes the proof of Theorem~\ref{thm:reductive}. 
\ \hfill $\Box$

\

\subsection{Generalized $1$-motives}

\noindent
In this section we point out two minor refinements of our
results. These refinements will not be used elsewhere in the paper.

The duality isomorphisms in Theorem~\ref{thm:gerbes} and
Theorem~\ref{thm:reductive} respect all the additional structures on
the stacks of Higgs bundles. For instance if $G$ is semisimple we can
view the stacks of Higgs bundles as generalized $1$-motives in the
sense of \cite{laumon}. We claim that the duality
isomorphisms respect the weight filtrations of these $1$-motives:

\

\medskip

  The isomorphism of group stacks
  $(\lan{\gHiggs})^{D} \cong \gHiggs$ in the Theorem~\ref{thm:gerbes}
  is compatible with the isomorphism of abelian schemes
\[
(\lan{\Higgs}_{0})^{D} = {\lan{P}}^{D} \cong P = \Higgs_{0},
\]
  constructed in the proof of Theorem~\ref{thm:duality}{\bf (2)}. More
  precisely, if we use $\bl$ to identify $B-\Delta$ with $\lan{B} -
  \lan{\Delta}$, then the Hitchin fibrations allow us to view
  $\gHiggs$ and $\lan{\gHiggs}$ as Beilinson $1$-motives over
  $B-\Delta$. This means that $\gHiggs$ and $\lan{\gHiggs}$ are
  commutative group stacks over $B-\Delta$, which are naturally
  filtered, with graded pieces which are either abelian varieties, or
  finite abelian groups, or classifying stacks of finite abelian
  groups. The filtrations are given as
\[
\xymatrix@R-1.5pc@C-2pc{
W_{0}\gHiggs \ar@{=}[d] &  \supset & W_{-1}\gHiggs \ar@{=}[d] & 
 \supset & W_{-2}\gHiggs \ar@{=}[d] & \supset 0 \\
\gHiggs & & \gHiggs_{0} & & BZ(G) & \\
}
\]
respectively
\[
\xymatrix@R-1.5pc@C-2pc{
W_{0}(\lan{\gHiggs}) \ar@{=}[d] & \supset & W_{-1}(\lan{\gHiggs})
\ar@{=}[d] & 
\supset & W_{-2}(\lan{\gHiggs}) \ar@{=}[d] & \supset 0 \\
\lan{\gHiggs} & & \lan{\gHiggs}_{0} & & BZ(\lan{G}) & \\
}
\]
and the duality operation $(\bullet)^{D} :=
\underline{\op{Hom}}_{\op{gp}}(\bullet, \mathcal{O}_{B}^{\times}[1])$,
is compatible with these filtrations:

\begin{lem}  \label{rem:motives}  The duality $(\bullet)^{D}$
transforms each filtered commutative group stack into a stack of the
same type, and the isomorphism $(\lan{\gHiggs})^{D} \cong \gHiggs$ 
respects the filtrations.
\end{lem}
{\bfseries Proof.}  The above filtrations give rise to short exact
sequences 
of commutative group stacks over $B - \Delta$:
\settowidth{\seqone}{$
\xymatrix@1{
0 \ar[r] & \gHiggs_{0} \ar[r] & \gHiggs \ar[r] & \pi_{1}(G) \ar[r] & 0
 & }$}
\[
\tag{$*$}
\left[
\begin{minipage}[c]{\seqone}
$
\xymatrix@R-2pc{
0 \ar[r] & \gHiggs_{0} \ar[r] & \gHiggs \ar[r] & \pi_{1}(G) \ar[r] & 0
\\ 
0 \ar[r] & BZ(G) \ar[r] & \gHiggs_{0} \ar[r] & \Higgs_{0} \ar[r] & 0 
}
$
\end{minipage}
\right]
\]
and 
\settowidth{\sone}{$\xymatrix@1{
0 \ar[r] & \Higgs_{0} \ar[r] & \Higgs \ar[r] & \pi_{1}(G) \ar[r] & 0 & }
$}
\[
\tag{$**$}
\left[
\begin{minipage}[c]{\sone}
$\xymatrix@R-2pc{
0 \ar[r] & \Higgs_{0} \ar[r] & \Higgs \ar[r] & \pi_{1}(G) \ar[r] & 0
\\ 
0 \ar[r] & BZ(G) \ar[r] & \gHiggs \ar[r] & \Higgs \ar[r] & 0.
}$
\end{minipage}
\right]
\]
Writing the same sequences for $\lan{\gHiggs}$ and applying
$(\bullet)^{D}$ we get 
\settowidth{\seqtwo}{$\xymatrix@1{
0 \ar[r] & (\lan{\Higgs}_{0})^{D} \ar[r] & (\lan{\gHiggs}_{0})^{D}
\ar[r] & Z(\lan{G})^{\wedge} \ar[r] & 0 & 0 0 0}$}
\[
\tag{$\lan{*}^{D}$}
\left[
\begin{minipage}[c]{\seqtwo}
$
\xymatrix@R-2pc{
0 \ar[r] & (\lan{\Higgs}_{0})^{D} \ar[r] & (\lan{\gHiggs}_{0})^{D}
\ar[r] & Z(\lan{G})^{\wedge} \ar[r] & 0 
\\ 
0 \ar[r] & B\pi_{1}(\lan{G})^{\wedge}  \ar[r] & (\lan{\gHiggs})^{D}
\ar[r] & (\lan{\gHiggs}_{0})^{D} \ar[r] & 0  
}
$
\end{minipage}
\right]
\]
and 
\settowidth{\stwo}{$\xymatrix@1{
0 \ar[r] & \lan{\Higgs}^{D} \ar[r] & (\lan{\gHiggs})^{D} \ar[r] &
Z(\lan{G})^{\wedge} \ar[r] & 0 & 0 0 0 }$
}
\[
\tag{$\lan{**}^{D}$}
\left[
\begin{minipage}[c]{\stwo}
$\xymatrix@R-2pc{
0 \ar[r] & (\lan{\Higgs})^{D} \ar[r] & (\lan{\gHiggs})^{D} \ar[r] &
Z(\lan{G})^{\wedge} \ar[r] & 0 
\\ 
0 \ar[r] & B\pi_{1}(\lan{G})^{\wedge} \ar[r] & (\lan{\Higgs})^{D}
\ar[r] & (\lan{\Higgs}_{0})^{D} \ar[r] & 0. 
}$
\end{minipage}
\right]
\]
The fact that the isomorphism \eqref{eq:gerbeD}  in
Theorem~\ref{thm:gerbes} respects the
filtrations is equivalent to showing that \eqref{eq:gerbeD} induces an
identification of short exact sequences $(*) \cong (\lan{**}^{D})$
(equivalently $(**) \cong (\lan{*}^{D})$). This follows by considering the
compatible isomorphisms of commutative group
stacks (or spaces) that we obtained in the proof of Theorem~\ref{thm:gerbes}:
\[
\xymatrix@R-0.5pc@C-2pc{
\gHiggs  & \cong & (\lan{\gHiggs})^{D} 
\\
\gHiggs_{0} \ar@{^{(}->}[u] \ar@{->>}[d] & \cong & (\lan{\Higgs})^{D}
\ar@{^{(}->}[u]\ar@{->>}[d]  
\\
\Higgs_{0} & \cong &(\lan{\Higgs}_{0})^{D}
}
\]
where in the top row we have the isomorphism \eqref{eq:gerbeD} from
Theorem~\ref{thm:gerbes}, and in
the bottom row we have the isomorphism from
Theorem~\ref{thm:duality}{\bf (2)}.  \ \hfill $\Box$

\

\medskip

  The proof of Theorem~\ref{thm:duality}
  and the calculation in the proof of Lemma~\ref{lem:auto}{\bf (i)}
  suggest that the duality of Hitchin systems proven in
  Theorem~\ref{thm:gerbes} admits a refinement in the case when the
  simple group $G$ is of type ${\sf{B}}_{r}$ or ${\sf{C}}_{r}$. In general, the
  inclusion of sheaves $\mathcal{T}^{o} \subset \mathcal{T} \subset
  \overline{\mathcal{T}}$ gives rise to three stacky integrable systems 
  over the Hitchin base $B$: 
\[
\xymatrix@R-1pc{
\bH^{o}_{G} \ar[r] \ar@{=}[d] & \gHiggs_{G} \ar[r] \ar@{=}[d] &
\overline{\bH}_{G}. \ar@{=}[d] \\
\gTors_{\mathcal{T}^{o}} & \gTors_{\mathcal{T}} &
\gTors_{\overline{\mathcal{T}}} 
}
\]

\

\begin{lem} \label{rem:spso}
The duality statement
$(\gHiggs_{\lan{G}})^{D} \cong \gHiggs_{G}$ extends to a
duality \linebreak $(\bH^{o}_{\lan{G}})^{D} \cong \overline{\bH}_{G}$.
\end{lem}
{\bf Proof.}
The corresponding coarse moduli spaces admit
cohomological interpretations as $H^{1}(\mathcal{T}^{o}_{G})$,
$H^{1}(\mathcal{T}_{G})$, and $H^{1}(\overline{\mathcal{T}}_{G})$ respectively.
These integrable systems coincide for all simple groups $G \neq \op{Sp}(r),
\op{SO}(2r+1)$, and 
\[
\begin{split}
\gHiggs_{\op{Sp}(r)} & \cong \overline{\bH}_{\op{Sp}(r)} \\
\bH^{o}_{\op{SO}(2r+1)} & \cong \gHiggs_{\op{SO}(2r+1)}.
\end{split}
\]
Now, the calculations in Claim~\ref{claim:components}{\bf (ii)} and
Lemma~\ref{lem:auto}{\bf (i)} give the following values for the
stabilizer groups and the groups of connected components of these
group stacks:

\

\begin{center}
\begin{tabular}{|l||c|c|c|} \hline
\begin{minipage}[c]{0.8in}
    \vspace{0.1in}
$G$ \vspace{0.1in} \end{minipage} & $H^{0}(\mathcal{T}^{o}_{G})$ & 
$H^{0}(\mathcal{T}_{G})$ &  $H^{0}(\overline{\mathcal{T}}_{G})$ \\
  \hline\hline
\begin{minipage}[c]{0.8in}
    \vspace{0.1in}$\op{Sp}(r)$ \vspace{0.1in} \end{minipage}
& $0$ & $Z(G) = \mathbb{Z}/2$ & $Z(G) = \mathbb{Z}/2$ \\
  \hline 
\begin{minipage}[c]{0.8in}
    \vspace{0.1in}$\op{SO}(2r+1)$ \vspace{0.1in} \end{minipage}
& $Z(G) = 0$ & $Z(G) = 0$  & $\mathbb{Z}/2$ \\
  \hline
\end{tabular}
\end{center}

\

and

\

\begin{center}
\begin{tabular}{|c||c|c|c|} \hline
\begin{minipage}[c]{0.8in}
    \vspace{0.1in}$G$ \vspace{0.1in} \end{minipage}
& $\pi_{0}(H^{1}(\mathcal{T}^{o}_{G}))$ &
$\pi_{0}(H^{1}(\mathcal{T}_{G}))$ &
  $\pi_{0}(H^{1}(\overline{\mathcal{T}}_{G}))$ \\ 
  \hline\hline
\begin{minipage}[c]{0.8in}
    \vspace{0.1in}$\op{Sp}(r)$ \vspace{0.1in} \end{minipage}
& $\mathbb{Z}/2$ & $\pi_{1}(G) = 0$ & $\pi_{1}(G) = 0$ \\
  \hline 
\begin{minipage}[c]{0.8in} \vspace{0.1in} $\op{SO}(2r+1)$ \vspace{0.1in}
\end{minipage} & $\pi_{1}(G) = \mathbb{Z}/2$ & $\pi_{1}(G) =
  \mathbb{Z}/2$  & $0$  \\
  \hline
\end{tabular}
\end{center}
The above tables and the calculation of the cocharacter lattices
of the Prym varieties $P^{o}$ and
$\overline{P}$ in Claim~\ref{claim:cochar} show that the isomorphism
$(\bH^{o}_{\lan{G}})^{D} \cong \overline{\bH}_{G}$ holds for the
graded pieces with respect to the weight filtrations. The full duality
of filtered objects follows from the argument in the proof of
Theorem~\ref{theo:iso.presentations}. \ \hfill $\Box$

\subsection{Equivalence of derived categories} \label{ss:cats}

Theorem~\ref{thm:gerbes}  has some immediate corollaries. First, we
get a categorical equivalence 

\begin{cor} \label{cor:FM} Over $B-\Delta$, there is a Fourier-Mukai
  type equivalence of derived categories
\[
\mathfrak{c} : D^{b}_{c}\left({\gHiggs} \right)
\widetilde{\to} D^{b}_{c}\left( \lan{\gHiggs}\right). 
\]
Moreover, for every $\alpha \in \pi_{0}(\Higgs) = \pi_{1}(G) =
Z(\lan{G})^{\wedge}$, and every $\beta\in \pi_{0}(\lan{\Higgs}) =
\pi_{1}(\lan{G}) = Z(G)^{\wedge}$, the functor $\mathfrak{c}$ gives
rise to
a Fourier-Mukai equivalence 
\[
D^{b}_{c}\left({}_{\beta}{\Higgs}_{\alpha}\right) \widetilde{\to} 
D^{b}_{c}\left({}_{\alpha}\lan{\Higgs}_{\beta} \right)
\]
for the derived categories of the induced $\mathcal{O}^{\times}$-gerbes.
\end{cor}
{\bf Proof.} The isomorphism \eqref{eq:gerbeD} implies that the
$\mathcal{O}^{\times}$-gerbes ${}_{\alpha}\lan{\Higgs}_{\beta}$ and
${}_{\beta}{\Higgs}_{\alpha}$ are compatible, in the sense of
\cite{dp}. In particular, the categorical equivalence statement from
\cite{dp} implies the equivalence of derived categories
$D^{b}_{c}(\gHiggs_{0}) = D^{b}_{c}(\lan{\Higgs})$. To get the full
categorical duality $D^{b}_{c}\left( \lan{\gHiggs}\right) \cong
D^{b}_{c}\left({\gHiggs} \right)$, one can combine \eqref{eq:gerbeD}
with the duality for representations of commutative group stacks
described in Arinkin's appendix to \cite{dp} (see also \cite{bb}), or
invoke the recent result \cite{oren} of O.Ben-Bassat. In fact,
Ben-Bassat's proof works in a much more general context and will imply
the full categorical duality even over the discriminant $\Delta$, as
long as one can show that the Poincare sheaf on the cameral Pryms
extends across $\Delta$. \ \hfill $\Box$



\subsection{Hecke eigensheaves} \label{ss:hecke}

\noindent
As observed in Lemma~\ref{rem:motives}, the duality in
Theorem~\ref{thm:gerbes} and \ref{thm:reductive} respects the weight
filtrations on $\gHiggs$ and $\lan{\gHiggs}$. In particular we have
$\gHiggs_{0} \cong (\lan{\Higgs})^{D}$ and so $\mathfrak{c}$ restricts
to a well defined equivalence
\[
\mathfrak{c}_{0} : D^{b}_{c}(\gHiggs_{0}) \widetilde{\to}
D^{b}_{c}(\lan{\Higgs}).
\]
Finally, we have that the natural orthogonal spanning class of the
category $D^{b}_{c}(\gHiggs_{0})$ is transformed by $\mathfrak{c}$
into 
the class of automorphic sheaves on $\lan{\Higgs}$. 
This is precisely the sense in which the categorical equivalence
$\mathfrak{c}$ can be thought of as a classical limit of the geometric
Langlands correspondence. To spell this out, recall that in the proof
of Theorem~\ref{thm:gerbes} (see also
Appendix~\ref{appendix:spectral}), we introduced abelianized Hecke 
maps 
\[
\lan{\trans}^{\mu} : \lan{\Higgs}_{\widetilde{C}}\times \widetilde{C}
\to \lan{\Higgs}_{\widetilde{C}} 
\]
labeled by characters $\mu \in \lan{\Lambda} = \Lambda^{\vee} =
\chr(T)$ of $T$.  These maps were constructed at the stack
level. Here we use specifically the induced maps on the
moduli spaces. Recall that the map $\lan{\trans}^{\mu}$ gives rise to  an abelianized Hecke
operator $\lanab{\mathbb{H}}^{\mu} :=
\left(\lan{\trans}^{\mu}\right)^{*}  :
D^{b}_{c}(\lan{\Higgs}_{\widetilde{C}}) \to
D^{b}_{c}(\lan{\Higgs}_{\widetilde{C}}\times \widetilde{C})$.

\begin{thm} \label{cor:hecke} A topologically trivial
  $G$-Higgs bundle $(V,\varphi)$ on $C$ determines an eigensheaf for
  the abelianized Hecke operators.
  Explicitly let $p : \widetilde{C} \to C$ be a cameral cover
  corresponding to a point in $B - \Delta$, and let
  $\mathcal{T}_{\widetilde{C}}$ be 
  the corresponding sheaf of regular centralizers on $C$. The choice of
  $(V,\varphi)$ gives:
\begin{itemize}
\item A  $\mathcal{T}$-torsor $\mycal{L}_{(V,\varphi)}$ on $\widetilde{C}$.
\item A representable structure morphism 
$\boldsymbol{\iota} : B\op{Aut}((V,\varphi)) \to \gHiggs_{0}$.
\end{itemize} 
Write $\mathfrak{o}_{(V,\varphi)} :=
\boldsymbol{\iota}_{*}\mathcal{O}_{B\op{Aut}((V,\varphi))}$ for the
corresponding sheaf on $\gHiggs_{0}$. (This is nothing but the
structure sheaf of the stacky point of $\gHiggs_{0}$ corresponding to
$(V,\varphi)$.) Then for every character $\mu \in
\Lambda^{\vee}$ we have a functorial isomorphism 
\[
\lanab{\mathbb{H}}^{\mu}\left(
\mathfrak{c}_{0}(\mathfrak{o}_{(V,\varphi)})\right) \cong  
\mathfrak{c}_{0}(\mathfrak{o}_{(V,\varphi)})\boxtimes \mu\left(
\mycal{L}_{(V,\varphi)}\right), 
\]
i.e. $\mathfrak{c}_{0}(\mathfrak{o}_{(V,\varphi)})$ is an abelianized Hecke
eigensheaf with eigenvalue 
$\mycal{L}_{(V,\varphi)}$. 
\end{thm}
{\bf Proof.} This is automatic from the definition of the Hecke
correspondences, the abelianization procedure of
\cite[Theorem~6.4]{ron-dennis}, and the fact that the categorical
equivalence $\mathfrak{c}$ is compatible with the usual Fourier-Mukai
equivalence of $\Higgs_{0}$ and $\lan{\Higgs}_{0}$ as discussed in
Theorem~\ref{theo:iso.presentations}. \ \hfill $\Box$

\section{The topological structure of a cameral Prym} 
\label{s:Prym_structure}

In this section we discuss the cohomology groups describing the
cocharacter lattices of cameral Prym varieties and the behavior of
those groups under Poincare duality.  
Most of the material in section \ref{ss:ls_cohomology} is
well known, but we couldn't find it in the literature, in the form
needed for the proof of Theorem~\ref{thm:duality}. 
We include here the necessary statements;
the proofs are left to the reader 
(or can be found at arXiv:math/0604617 v1).
The results in Section~\ref{ss:topo} are new. 
We give an explicit description of the
cocharacter lattice of a cameral Prym in terms of the local
monodromies of a cameral cover.  We used this result in the proof of
Claim~\ref{claim:components} to analyze the connected components of
the Hitchin fiber, but it may also be of independent interest.

\subsection{Remarks on local system cohomology} \label{ss:ls_cohomology}

Let $C$ be a smooth compact complex curve of genus $g$ and let $S = \{
s_{1}, \ldots, s_{b} \} \subset C$ be a finite set of points. We write
$U := C - S$ for the complement of $S$ and denote by $\imath : S
\hookrightarrow C$ and $\jmath : U \hookrightarrow C$ the
corresponding closed and open inclusions. We will also fix a base
point $\bpo \in U$.

Let $A$ be a local system on $U$ of free abelian groups of rank
$r$. Let $A^{\vee} :=
\underline{\op{Hom}}_{\mathbb{Z}_{U}}(A,\mathbb{Z}_{U})$ denote the
dual local system. We want to understand the cohomology of $\jmath_{*}A$
in concrete terms and to find the precise relationship between the
cohomology of $\jmath_{*}A$ and $\jmath_{*}(A^{\vee})$. This is all
standard for local systems of vector spaces, see
e.g. \cite{looijenga}, but it requires some care for local system of
free abelian groups.

Suppose $s_{i} \in S$ and let $s_{i} \in \bD_{i} \subset C$ be a small disc
centered at $s_{i}$ and not containing any other point of $S$. Fix a point
$\bpo_{i} \in \partial \bD_{i}$ and let $c_{i}$ denote the loop starting
and ending at $\bpo_{i}$ and traversing $\partial \bD_{i}$ once in the
positive direction. Write $\op{mon}(c_{i}) : A_{\bpo_{i}} \to
A_{\bpo_{i}}$ for the monodromy operator associated with $c_{i}$. Now,
from
the definition of the direct image and the fact that $A$ is locally
constant we get the following description of the stalk of
$\jmath_{*}A$ at $s_{i}$: 
\[
\begin{split}
(\jmath_{*}A)_{s_{i}} & = \lim_{\longleftarrow} \left\{
    \left. H^{0}(V\cap U, A) 
    \right| \  s_{i} \in V, \
    V \subset C 
    \text{ - open } \right\} \\
& = H^{0}(\bD_{i}-\{s_{i} \}, A) \\
& = (A_{\bpo_{i}})^{\op{mon}(c_{i})}.
\end{split}
\]
Here, as usual $(A_{\bpo_{i}})^{\op{mon}(c_{i})} := \{ a \in A_{\bpo_{i}} |
\op{mon}(c_{i})(a) = a \}$ denotes the invariants of the
$\op{mon}(c_{i})$-action. 

\begin{figure}[!ht]
\begin{center}
\psfrag{s}[c][c][1][0]{{$s_{i}$}}
\psfrag{o}[c][c][1][0]{{$\bpo$}}
\psfrag{os}[c][c][1][0]{{$\bpo_{i}$}}
\psfrag{Ds}[c][c][1][0]{{$\bD_{i}$}}
\psfrag{as}[c][c][1][0]{{$a_{i}$}}
\epsfig{file=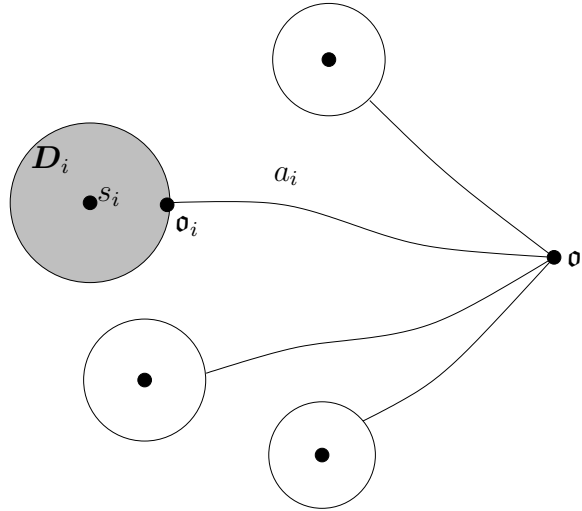,width=3in} 
\end{center}
\caption{An arc system for $S\subset C$.}\label{fig:arcs} 
\end{figure}

To organize things better, we choose 
an ordered system of arcs $\{ a_{i}
\}_{i = 1}^{b}$ in $C - \cup_{i=1}^{b} \bD_{i}$ which connect the base
point $\bpo$ with each of the points $\bpo_{i}$ as in
Figure~\ref{fig:arcs}.

\

\noindent
A choice of an arc system yields a collection of elements $\gamma_{i} \in
\pi_{1}(U,\bpo)$. Geometrically $\gamma_{i}$ is the $\bpo$-based
loop in $U$ obtained by tracing $a_{i}$, followed by tracing $c_{i}$
and then tracing back $a_{i}$ in the opposite direction. Since
parallel transport along $a_{i}$ identifies the stalks $A_{\bpo}$ and
$A_{\bpo_{i}}$ and conjugates the mondromy transformation
$\op{mon}(\gamma_{i})$ into the monodromy transformation
$\op{mon}(c_{i})$, it follows that we also have the
identification 
\[
(\jmath_{*}A)_{s_{i}} = (A_{\bpo})^{\op{mon}(\gamma_{i})} 
\]
for all $s_{i} \in S$. To simplify notation we set $\rho_{i} :=
\op{mon}(\gamma_{i})$.

With this notation in place we are ready to analyze the cohomology of
the sheaves $A$ and $\jmath_{*}A$. First note that $U$ is a smooth
$2$-manifold and so has homological dimension $2$ with respect to
compactly supported cohomology \cite[Section~III.9]{iversen} or
\cite[Section~3.1]{dimca-sheaves}. Also, since $A$ is a locally constant
sheaf, it can not have any compactly supported sections and so the
only compactly supported cohomology groups of $A$ that can be
potentially non-zero are $H^{1}_{c}(U,A)$ and
$H^{2}_{c}(U,A)$. On the other hand, since $A$ is a local system, its
cohomology is homotopy invariant. Taking into account the fact that
$U$ is homotopy equivalent to a bouquet of circles, we conclude that
the only cohomology groups of $A$ that are potentially non-zero are
$H^{0}(U,A)$ and $H^{1}(U,A)$. Furthermore the following version of
Poincare duality holds for these groups:

\begin{lem} \label{lem:duality} The cup product pairing
\[ 
\label{eq:cup}
\tag{\text{\bfseries\sf cup}}
\xymatrix@1{H^{k}_{c}(U,A)\otimes H^{2-k}(U,A^{\vee}) \ar[r]^-{\cup} & 
  H^{2}_{c}(U,A\otimes A^{\vee}) \ar[r]^-{\int \op{tr}} & \mathbb{Z}},
\]
induces a perfect pairing between the free abelian groups 
$H^{1}(U,A)_{\op{tf}}$ and 
$H^{1}_{c}(U,A^{\vee})_{\op{tf}}$. Moreover $H^{0}(U,A)$ has
no torsion and the cup product pairing \eqref{eq:cup} induces a
perfect pairing between  $H^{0}(U,A)$ and $H^{2}_{c}(U,A^{\vee})_{\op{tf}}$. 
\end{lem}

The proof, using Verdier duality and the universal coefficient theorem, is omitted.

\

\hfill

\noindent
Next we will use this information to compute the cohomology of $A$ and
$\jmath_{*}A$ explicitly in terms of the monodormy.
We begin with a  standard lemma:

\begin{lem} \label{lem:Rj} The direct image sheaves $R^{k}\jmath_{*}A$ can
  be described as follows:
\begin{itemize}
\item[{\bf (a)}] The sheaf $\jmath_{*}A$ fits in a short exact sequence
\[
0 \to \jmath_{!}A \to \jmath_{*}A \to \oplus_{i=1}^{b}
A_{\bpo}^{\rho_{i}} \to 0, 
\]
where $\jmath_{!}$ denotes the pushforward with compact supports.
\item[{\bf (b)}] The sheaf $R^{1}\jmath_{*}A$ satisfies 
\[
R^{1}\jmath_{*}A = \oplus_{i=1}^{b} (A_{\bpo})_{\rho_{i}},
\]
where $(A_{\bpo})_{\rho_{i}} := A_{\bpo}/(1-\rho_{i})A_{\bpo}$ denotes
the group of coinvariants of the $\rho_{i}$-action on $A_{\bpo}$. 
\item[{\bf (c)}] $R^{k}\jmath_{*}A = 0$ for all $k \geq 2$.
\end{itemize}
\end{lem}
{\bf Proof.} 
One can easily compute $H^{k}(c_{i},A)$
  via group cohomology. Since $c_{i} \cong S^{1} = K(\mathbb{Z},1)$ we
  have
\[
H^{k}(c_{i}, A) =H^{k}(\pi,M),
\]
where $\pi := \pi_{1}(c_{i},\bpo_{i}) \cong \mathbb{Z}$, and $M$
denotes the $\pi$ module $(\mathfrak{A},\mathfrak{r})$. Now the group
ring $\mathbb{Z}[\pi] = \mathbb{Z}[t]$ is a polynomial ring in one
variable over $\mathbb{Z}$ and so $\mathbb{Z}$ has a two step
resolution 
\[
\xymatrix@1{
0 \ar[r] &  \mathbb{Z}[t] \ar[r]^-{\partial} & \mathbb{Z}[t]
\ar[r]^-{\varepsilon} & \mathbb{Z} \ar[r] & 0
}
\]
by free $\mathbb{Z}[\pi]$-modules. Here $\varepsilon : \mathbb{Z}[t]
\mathbb{Z}$ is the augmentation map $\varepsilon(p(t)) := p(0)$, and
$\partial : \mathbb{Z}[t] \to \mathbb{Z}[t]$ is the map $\partial p(t)
:= (t-1)p(t)$ of multiplication by $t - 1$. In particular, for any
$\pi$ module $M = (\mathfrak{A},\mathfrak{r})$ we can compute
$H^{\bullet}(\pi,M)$ as the cohomology of the complex
\[
\xymatrix@R-1pc{
(\text{degree } 0) & (\text{degree } 1) \\
\op{Hom}_{\pi}(\mathbb{Z}[t],M) \ar[r]^-{\delta} \ar@{=}[d] &
\op{Hom}_{\pi}(\mathbb{Z}[t],M) \ar@{=}[d] \\
\mathfrak{A} \ar[r]_-{\mathfrak{r} - 1} & \mathfrak{A}.
}
\]
In other words $H^{0}(\pi,M) = \mathfrak{A}^{\mathfrak{r}}$,
$H^{1}(\pi,M) = \mathfrak{A}_{\mathfrak{r}}$, and $H^{k}(\pi,M) = 0$
for $k \geq 2$. We leave the remaining details to the reader. 

Alternatively, one can give a very elementray argument based on 
the usual gluing \cite[Section~1.4]{bbd}, 
\cite[Chapter~4,Exercise~3]{gelfand-manin} for the inclusions $U
\stackrel{\jmath}{\hookrightarrow} C \stackrel{\imath}{\hookleftarrow} S$,
which yields a distinguished triangle 
\[
\jmath_{!}\jmath^{-1}F^{\bullet} \to F^{\bullet} \to
\imath_{*}\imath^{-1}F^{\bullet} \to \jmath_{!}F^{\bullet}[1]
\]
defined for every $F^{\bullet} \in D(\mathbb{Z}_{U}-\op{mod})$. 
Our lemma follows by taking $F^{\bullet} := R\jmath_{*} A$. 
\ \hfill $\Box$

\

\bigskip

\bigskip

\noindent
As a consequence of the calculation of $R^{k}\jmath_{*}A$ one
immediately gets the following:

\begin{lem} \label{lem:Hij}
The cohomology of the sheaf $\jmath_{*} A$ satisfies:
\[
\begin{split}
H^{0}(C,\jmath_{*}A) & = A_{\bpo}^{\pi_{1}(U,\bpo)} \\
H^{1}(C,\jmath_{*}A) & = \op{im}\left[ H^{1}_{c}(U,A) \to H^{1}(U,A)
  \right] \\
H^{2}(C,\jmath_{*}A)_{\op{tf}} & =
((A_{\bpo}^{\vee})^{\pi_{1}(U,\bpo)})^{\vee}  \\
H^{2}(C,\jmath_{*}A)_{\op{tor}} & = H^{1}(U,A^{\vee})^{\wedge}_{\op{tor}}.
\end{split}
\]
\end{lem} 

\

\begin{cor} \label{cor:image} For any local system $A$ of finite rank
  free abelian 
  groups on $U$ we have a natural identification
\[
H^{1}(C,\jmath_{*}A^{\vee})_{\op{tf}} = \op{im} \left[ H^{1}(U,A)^{\vee} \to
  H^{1}(C,\jmath_{*}A)^{\vee} \right] 
\]
\end{cor}
{\bf Proof.} By the previous lemma we have an identification
\[
H^{1}(C,\jmath_{*}A^{\vee}) = \op{im} \left[ H^{1}_{c}(U,A^{\vee})
\to H^{1}(U,A^{\vee}) \right].
\]
In particular we have
$H^{1}(C,\jmath_{*}A^{\vee})_{\op{tf}} = \op{im} \left[
H^{1}_{c}(U,A^{\vee})_{\op{tf}} \to H^{1}(U,A^{\vee})_{\op{tf}}
\right]$ On the other 
hand, by Lemma~\ref{lem:duality} the natural map 
$H^{1}_{c}(U,A^{\vee})_{\op{tf}} \to H^{1}(U,A^{\vee})_{\op{tf}}$ 
is equal to minus the transpose of the map 
$H^{1}_{c}(U,A)_{\op{tf}} \to H^{1}(U,A)_{\op{tf}}$. Since by
Lemma~\ref{lem:Hij} we have
\[
H^{1}(C,\jmath_{*}A) = \op{im}\left[H^{1}_{c}(U,A) \to
  H^{1}(U,A)\right],
\] 
we get that 
\[
\op{im} \left[ H^{1}(U,A)^{\vee} \to H^{1}_{c}(U,A)^{\vee}\right] = 
\op{im} \left[ H^{1}(U,A)^{\vee} \to
  H^{1}(C,\jmath_{*}A)^{\vee}\right],
\] 
which
yields the lemma. \
\hfill $\Box$

\subsection{The cocharacters of a cameral Prym} \label{ss:topo} 

Let $p : \widetilde{C} \to C$ be a generic Galois cameral cover as in
the proof of Theorem~\ref{thm:duality}. Let $\{ s_{1}, \ldots, s_{b}
\} \subset C$ be the branch points of this cover. We write $\jmath : U
\hookrightarrow C$ for the inclusion of the complement, and $p^{o} :
p^{-1}(U) \to U$ for the unramified part of $p$. Define a local system
$A$ on $U$ by $A := (p^{o}_{*}\Lambda)^{W}$. The canonical identification
$\lan{\Lambda} = \Lambda^{\vee} = \op{Hom}(\Lambda,\mathbb{Z})$ gives
also an identification $\lan{A} = A^{\vee} =
\underline{\op{Hom}}(\Lambda,\mathbb{Z})$. 

Fix a base point $\bpo
\in U$ and choose an arc system as in the previous section. We choose
once and for all an identification $A_{\bpo} \cong \Lambda$. By
definition, the
monodromy  $\op{mon}_{\bpo} : \pi_{1}(U,\bpo) \to GL(\Lambda)$ factors
through $W \subset SO(\Lambda,\langle
\bullet, \bullet \rangle) \subset GL(\Lambda)$. By the genericity
assumption on  $p : \widetilde{C} \to C$, it follows that the
monodromy image of each generator
$\gamma_{i} \in \pi_{1}(U,\bpo)$ is a reflection $\rho_{i} : \Lambda
\to \Lambda$ corresponding to some root $\alpha_{i}$ of
$\mathfrak{g}$.  Explicitly we have $\rho_{i}(\lambda) = \lambda -
(\alpha_{i},\lambda)\cdot \alpha_{i}^{\vee}$, where $\alpha_{i}^{\vee}
\in \Lambda$ is the coroot corresponding to $\alpha_{i}$. Since for
each root $\alpha$  the
divisor $D^{\alpha} \subset \op{tot}(K_{C}\otimes \mathfrak{t})$ is
ample, it follows that the
collection of roots $\{ \alpha_{1},
\ldots, \alpha_{b} \}$ contains both long and short roots of
$\mathfrak{g}$. Let
$\varepsilon_{i} := \varepsilon_{\alpha_{i},G}$ and
$\varepsilon_{i}^{\vee} := \varepsilon^{\vee}_{\alpha_{i},G}$. 
 
In the proof of Theorem~\ref{thm:duality} we described the cocharacter
lattice of the cameral Prym $P$ corresponding to $p : \widetilde{C}
\to C$ and $G$ in terms of the first cohomology of the sheaves $\jmath_{*}A$
and $\jmath_{*}A^{\vee}$ on $C$. We now give explicit formulas for
these cohomology groups. Note that without a loss of generality we may
assume that $S$, the arcs $a_{i}$ and the discs $\bD_{i}$ are all
contained in the interior of a disc $\bD \subset C$ for which $\bpo \in
\partial \bD$.  In particular we can choose a collection of $\bpo$ based loops
$\delta_{1}, \delta_{2}, \ldots, \delta_{2g} \subset C-\bD$
which intersect only at $\bpo$ and form a system of standard  $a$-$b$
generators for the fundamental group $\pi_{1}(C,\bpo)$ of the compact
curve $C$. Choosing the orientation of the  loops $\delta_{j}$ 
and $\gamma_{i}$ appropriately we get a presentation
of the fundamental group of $U$:
\[
\pi_{1}(U,\bpo) = \left\langle \left. \begin{minipage}[c]{0.8in}
  $\delta_{1}, \ldots, \delta_{2g}$, \\
$\gamma_{1}, \ldots, \gamma_{b}$ \end{minipage}
\right| \prod_{j= 1}^{g} [\delta_{i},\delta_{g+i}] \prod_{i=1}^{b} \gamma_{i}
  = 1
  \right\rangle 
\]
To simplify notation we set $\bw_{i} := \op{mon}_{\bpo}(\delta_{i})
\in W$. Finally, it will be convenient to add to $S$ an extra point
$s_{0} \neq \bpo \in U$ and a loop $\gamma_{0}$ arount $s_{0}$,
oriented so that the relation defining $\pi_{1}(U - \{s_{0}\},\bpo)$
is $\gamma_{0} = \prod_{j =1}^{g}[\delta_{j}, \delta_{g+j}] \prod_{i =
1}^{b} \gamma_{i}$. Note that we have $\rho_{0} =
\op{mon}_{\bpo}(\gamma_{0}) = 1 \in W$. Note also that the deletion of
$s_{0}$ from  $U$ will not affect our interpretation of
$H^{1}(C,\jmath_{*}A)$ as a kernel of a homomorphism. That is, we
still have 
\[
\begin{split}
H^{1}(C,\jmath_{*}A) & = \ker\left[ H^{1}(U,A) \to \oplus_{i =
1}^{b} \Lambda/(1 - \rho_{i})\Lambda \right] \\
& = \ker\left[ H^{1}(U-\{s_{0}\},A) \to \oplus_{i =
0}^{b} \Lambda/(1 - \rho_{i})\Lambda \right],
\end{split}
\]
since the deletion adds a copy of $\Lambda$ to both
$H^{1}(\bullet,A)$ and $R^{1}\jmath_{*} A$. 

\begin{prop} \label{prop:explicit} 
\begin{itemize}
\item[{\bf (a)}] There is a natural isomorphsim 
\[
H^{1}(U-\{ s_{0} \},A) \cong \frac{\Lambda^{2g + b}}{(1 - \bw_{1},
  \ldots, 1 - 
  \bw_{2g}, 1 - \rho_{1}, \ldots, 1- \rho_{b})\Lambda} 
\]
which depends only on the choice of an arc system, the $a-b$ loops
  $\delta_{j}$, and the identification $A_{\bpo} \cong \Lambda$.
In addition, there is a non-canonical isomorphism 
\[
H^{1}(U-\{ s_{0} \},A) = H^{1}(C,\Lambda)\oplus  \frac{\Lambda^{b}}{(1
  - \rho_{1}, 
  \ldots, 1- \rho_{b})\Lambda}.
\]
\item[{\bf (b)}] Under the isomorphism $H^{1}(U-\{ s_{0} \},A) =
  H^{1}(C,\Lambda)\oplus (\Lambda^{b}/(1 - \rho_{1}, \ldots, 1-
  \rho_{b})\Lambda)$, the subgroup $H^{1}(C,\jmath_{*}A) \subset
  H^{1}(U,A)$ can be identified as
\[
H^{1}(C,\jmath_{*}A) = H^{1}(C,\Lambda)\oplus
  \frac{\ker\left[ \sum_{i = 1}^{b}
  \prod_{k = 1}^{i-1} \rho_{k} : 
\oplus_{i=1}^{b}
  \mathbb{Z}\varepsilon_{i}\alpha_{i}^{\vee}  \to \Lambda \right]}{(1
  - \rho_{1}, 
  \ldots, 1- 
  \rho_{b})\Lambda} 
\]
\end{itemize}
\end{prop}
{\bf Proof.} {\bf (a)} The surface $U - \{ s_{0} \} = C - \{ s_{0},
s_{1}, \ldots, s_{b} \}$ is homotopy equivalent to the bouquet
consiting of the $2g + b$ oriented circles $\delta_{1}, \ldots,
\delta_{2g}, \gamma_{1}, \ldots, \gamma_{b}$, where all circles are
attached to each other at the point $\bpo$. The fundamental group
$\pi$ of this bouquet of circles is a free group on the generators
$\delta_{1}, \ldots, \delta_{2g}, \gamma_{1}, \ldots, \gamma_{b}$, and
the local system $A$ corresponds to the action on of this free group
on $\Lambda$ specified by the monodromy transformations $(\bw_{1},
\ldots, \bw_{2g}, \rho_{1}, \ldots, \rho_{b} ) \in W^{2g +
  b}$. As in the proof of Lemma~\ref{lem:Rj}, the trivial $\pi$-module
$\mathbb{Z}$  has a free
$\mathbb{Z}[\pi]$-resolution  given by 
\[
\xymatrix@1{ 0 \ar[r] & \left(\oplus_{j = 1}^{2g}
\mathbb{Z}[\pi]e_{\delta_{j}}\right) \oplus \left(\oplus_{i = 1}^{b}
\mathbb{Z}[\pi]e_{\gamma_{i}}\right) \ar[r]^-{\partial} &
\mathbb{Z}[\pi] \ar[r] & \mathbb{Z} \ar[r] & 0,  }
\]
where $\left(\oplus_{j = 1}^{2g}
\mathbb{Z}[\pi]e_{\delta_{j}}\right) \oplus \left(\oplus_{i = 1}^{b}
\mathbb{Z}[\pi]e_{\gamma_{i}}\right)$ is the free $\mathbb{Z}[\pi]$
module on generators $e_{\delta_{j}}$, $e_{\gamma_{i}}$ and $\partial
 e_{\delta_{j}} = 1 - \delta_{j}$, $\partial
 e_{\gamma_{i}} = 1 - \gamma_{i}$.  

Applying $\op{Hom}_{\mathbb{Z}[\pi]}(\bullet, \Lambda)$ and computing
cohomology we get the identification 
\[
H^{1}(U-\{ s_{0} \},A) = H^{1}(\pi,\Lambda) \cong \frac{\Lambda^{2g +
  b}}{(1 - \bw_{1}, 
  \ldots, 1 - 
  \bw_{2g}, 1 - \rho_{1}, \ldots, 1- \rho_{b})\Lambda}.
\]
Next we claim that by making appropriate choices, the topological
description  description of the cameral cover can be brought into a
particualraly simple form. Recall \cite{ddp}, that there is a natural
inclusion
\begin{equation} \label{eq:inclusion}
(H^{0}(C,K_{C})\otimes \mathfrak{t})/W \hookrightarrow B. 
\end{equation}
A cameral cover $p_{b} : \widetilde{C}_{b} \to C$ corresponding to a
generic point $b$ in the image of \eqref{eq:inclusion} is reduced but
completely reducible:
\[
\widetilde{C}_{b} = \bigcup_{w \in W} \widetilde{C}_{b,w},
\] 
with each irreducible component $\widetilde{C}_{b,w}$ isomorphic to
$C$. Let $\bD \subset C$ be a disc containing the image of all the
singular points (= intersection of components) of $\widetilde{C}_{b}$.
We get that $p^{-1}(C - \bD)$ is completely disconnected:
\begin{equation} \label{eq:disconnected}
p^{-1}(C - \bD) = \coprod_{w \in W} [C-\bD]_{w}, \quad [C-\bD]_{w} :=
\widetilde{C}_{b,w}\cap p^{-1}(C - \bD),
\end{equation}
with each connected components $[C-\bD]_{w}$ isomorphic to $C - \bD$. 

A general cameral cover $\widetilde{C}_{b'}$ with $b' \in B - \Delta$
near $b \in B$
will be smooth and will still satisfy \eqref{eq:disconnected}. By taking
all the $\gamma_{i}$ in $\bD$ and all the $\delta_{j}$ in 
$C-\bD$ we have $\bw_{j} = 1$, $j = 1, \ldots, 2g$. Consequently
\begin{equation} \label{eq:noncan}
H^{1}(U-\{ s_{0} \},A) = H^{1}(C,\Lambda)\oplus  \frac{\Lambda^{b}}{(1
  - \rho_{1}, 
  \ldots, 1- \rho_{b})\Lambda},
\end{equation}
for the cover $\widetilde{C}_{b'}$. Since $B - \Delta$ is connected,
it follows that \eqref{eq:noncan} will hold for any
$\widetilde{C}_{b''}$, $b'' \in B - \Delta$ and an appropriate choice
of $\gamma_{i}$'s and $\delta_{j}$'s.

\

\medskip

\noindent
{\bf (b)} As we argued in the previous section, the group
$H^{1}(C,\jmath_{*}A)$ is the kernel of the natural map
$H^{1}(U-\{ s_{0} \},A) \to H^{0}(C,R^{1}\jmath_{*}A) = \oplus_{i =
  0}^{b} \Lambda/(1-\rho_{i})\Lambda$. Under the identification
\eqref{eq:noncan}, it is immediate that $H^{1}(C,\Lambda)$ is
contained in
the kernel of this map and that the restriction of the map to the
summand $\Lambda^{b}/(1-\rho_{1}, \ldots, \ldots, 1-\rho_{b})\Lambda$
is given by 
\begin{equation} \label{eq:maponsummand}
(\lambda_{1}, \ldots, \lambda_{b}) \mapsto
(\varphi(\lambda_{1}, \ldots, \lambda_{b}), \lambda_{1} +
(1-\rho_{1})\Lambda, \ldots, \lambda_{b} + (1-\rho_{b})\Lambda),
\end{equation}
where the map $\varphi : \Lambda^{b} \to \Lambda/(1-\rho_{0})\Lambda =
\Lambda$ corresponds to the relation $\prod_{j
=1}^{2g}[\delta_{j},\delta_{g+j}]\prod_{i =1}^{b} \gamma_{i} =
\gamma_{0}$. 

In general, suppose that we are given a bouquet of circles $V =
c_{1}\vee \ldots \vee c_{n}$ and suppose we have a word $d =
\prod_{i = 1}^{p} d_{i}$ in $\pi := \pi_{1}(V)$ for which all
$d_{i}$'s are in $\{ c_{1}, \ldots, c_{n}, c_{1}^{-1}, \ldots,
c_{n}^{-1} \} \subset \pi$. Consider the cyclic subgroup in $\pi$
generated by $d$ and let $M$ be some $\pi$-module. The inclusion 
$\langle d \rangle \subset \pi$ induces a map on cohomology
$H^{1}(\pi,M) \to H^{1}(\langle d \rangle, M)$ which can
be explicitly calculated. For this we only need to consider the tree which is
the universal cover of $V$ and follow the branches of that tree
labeled by the letters $d_{i}$ in the word $d$. More invariantly this
corresponds to a map of resolutions 
\[
\xymatrix@R-1pc{
0 \ar[r] & 
\mathbb{Z}[\langle d \rangle ]e_{d}  \ar[r]^-{\partial} \ar[d]^-{(**)} &
\mathbb{Z}[\langle d \rangle] \ar[r] \ar[d]^-{(*)} & \mathbb{Z} \ar@{=}[d]
\ar[r] & 0 \\ 
0 \ar[r] & \oplus_{i = 1}^{n}
\mathbb{Z}[\pi]e_{c_{i}}  \ar[r]^-{\partial} &
\mathbb{Z}[\pi] \ar[r] & \mathbb{Z} \ar[r] & 0,
}
\]
where $(*)$ is the natural inclusion and $(**)$ sends $e_{d}$ to the
  sum $\sum_{j = 1}^{p} \op{sgn}(d_{j}) \left(\prod_{i = 1}^{j-1}
  d_{i}\right) e_{d_{j}}$, where $\op{sgn}(d_{j}) = \pm 1$ depending on
  whether $d_{j}$ is one of the $c_{i}$'s or one of the
  $c_{i}^{-1}$'s. 

Combining this formula with the observation that
$(1-\rho_{i})(\lambda) = (\alpha_{i},\lambda)\alpha_{i}^{\vee}$, we
see that the kernel of \eqref{eq:maponsummand} is precisely 
$\ker\left[ \sum_{i = 1}^{b}
  \prod_{k = 1}^{i-1} \rho_{k} : 
\oplus_{i=1}^{b}
  \mathbb{Z}\varepsilon_{i}\alpha_{i}^{\vee}  \to \Lambda \right]$.  
  \ \hfill $\Box$

\

\bigskip

\noindent
In particular, for the torsion subgroups of $H^{1}(C,\jmath_{*}A)$ and
$H^{1}(U,A)_{\op{tor}}$  we get

\begin{cor} \label{cor:torsion}
\[
\begin{split}
H^{1}(C,\jmath_{*}A)_{\op{tor}} & = \left(\frac{\oplus_{i = 1}^{b}
  \mathbb{Z} \varepsilon_{i}\alpha_{i}^{\vee}}{(1 - \rho_{1}, \ldots, 1-
  \rho_{b})\Lambda}\right)_{\op{tor}} \\ 
& \\
H^{1}(U,A)_{\op{tor}} & = H^{1}(U-\{s_{0}\},A)_{\op{tor}} =
  \left(\frac{\Lambda^{b}}{(1 - \rho_{1}, \ldots, 1- 
  \rho_{b})\Lambda}\right)_{\op{tor}}. 
\end{split}
\]
\end{cor}
{\bf Proof.} Since $H^{1}(C,\Lambda)$ is torsion free,
  Proposition~\ref{prop:explicit} implies that 
\[
H^{1}(C,\jmath_{*}A)_{\op{tor}} = \left(\frac{\ker\left[ \sum_{i = 1}^{b}
  \prod_{k = 1}^{i-1} \rho_{k} : \oplus_{i = 1}^{b}
  \mathbb{Z} \varepsilon_{i}\alpha_{i}^{\vee} \to \Lambda \right]}{(1
  - \rho_{1}, \ldots, 1- 
  \rho_{b})\Lambda}\right)_{\op{tor}}.
\]
The corollary now follows by noticing that the saturation of $(1 -
  \rho_{1}, \ldots, 1- \rho_{b})\Lambda$ inside the lattice 
$\oplus_{i = 1}^{b}
  \mathbb{Z} \varepsilon_{i}\alpha_{i}^{\vee}$ is the same as the
  saturation of $(1 - \rho_{1}, \ldots, 1- \rho_{b})\Lambda$ inside
  the lattice $\ker\left[ \sum_{i = 1}^{b} \prod_{k = 0}^{i-1} \rho_{k} :
  \oplus_{i=1}^{b} \mathbb{Z}\varepsilon_{i}\alpha_{i}^{\vee} \to
  \Lambda \right]$, since the latter lattice is a kernel to a map to
  $\Lambda$ which is torsion free. 

Finally, $H^{1}(U,A)_{\op{tor}} = H^{1}(U-\{s_{0}\},A)_{\op{tor}}$
  since the deletion of $s_{0}$ adds a copy of $\Lambda$ as a direct
  summand. 
\ \hfill $\Box$

\

\bigskip

\bigskip

\appendix

\Appendix{Hecke functors and spectral
  data} \label{appendix:spectral}

\setcounter{equation}{0}

\

\noindent
In section~\ref{s:cl} we defined the classical limit Hecke functors
\[
\lan{\mathbb{H}}^{\mu,x} :
D^{b}_{\op{qcoh}}(\lan{\gHiggs},\mathcal{O}) \to
D^{b}_{\op{qcoh}}(\lan{\gHiggs},\mathcal{O}),
\]
which were labeled by pairs $(\mu,x)$ with $\mu \in \chr(G)$ and
$x \in C$. On the other hand the derived category of quasi-coherent
sheaves on $\lan{\gHiggs}$ is equipped with another collection of
endo-functors: the abelianized Hecke functors
$\lanab{\mathbb{H}}^{\mu,\tilde{x}}$. Here again $\mu \in \chr(G)$, but
$\tilde{x}$ is a point in some cameral cover $\widetilde{C} \to C$
which we will assume smooth, i.e. $\widetilde{C} = \widetilde{C}_{b}$,
where $b \in B - \Delta$. By definition the functor
$\lanab{\mathbb{H}}^{\mu,\tilde{x}}$ is the integral transform on
quasi-coherent sheaves on $\lan{\gHiggs}$ whose kernel is the
structure sheaf of the abelianized Hecke correspondence
\[
\xymatrix{
& {\lanab{\sHeck}}^{\mu,\tilde{x}}_{\widetilde{C}} \ar[dl]_-{p^{\mu,\tilde{x}}}
  \ar[dr]^-{q^{\mu,\tilde{x}}} & \\ 
\lan{\gHiggs}_{\widetilde{C}} & & \lan{\gHiggs}_{\widetilde{C}}
}
\]
Here $\lan{\gHiggs}_{\widetilde{C}}$ denotes the moduli stack of
$\lan{G}$-Higgs bundles on $C$ with cameral cover $\widetilde{C}$,
i.e. the fiber of the (stacky) Hitchin map over $\widetilde{C}$. The
correspondence ${\lanab{\sHeck}}^{\mu,\tilde{x}}_{\widetilde{C}}$ is
the graph of a specific translation isomorphism
\[
\lan{\trans}^{\mu,\tilde{x}} : \lan{\gHiggs}_{\widetilde{C}} \to
\lan{\gHiggs}_{\widetilde{C}},
\]
that we describe next. The  essential ingredient in the definition of
$\lan{\trans}^{\mu,\tilde{x}}$ is  the description
\cite[Theorem~6.4]{ron-dennis} of
the fiber of the Hitchin map $\lan{\gHiggs}_{\widetilde{C}}$ in terms
of geometric data on $\widetilde{C}$: the so called {\em spectral
  data}.  

To make things explicit we recall this description next.

\subsection{Spectral data for principal Higgs bundles} \label{sec:spectral}

\punkt {\bf Bundles and cocycles.} \ Suppose $\widetilde{C} \to C$ is
a fixed abstract cameral cover for $\lan{G}$. For every root $\alpha$
of $\lan{G}$ denote by $D^{\alpha} \subset \widetilde{C}$ the divisor
fixed by the reflection $\rho_{\alpha} \in W$ corresponding to
$\alpha$.  Consider the following principal bundles:
\begin{itemize}
\item[$\sR_{\alpha}$:] the $\mathbb{C}^{\times}$-bundle
  corresponding to $D_{\alpha}$, i.e. $\sR_{\alpha} :=
  \mathcal{O}_{\widetilde{C}}(D_{\alpha})^{\times}$. 
\item[$\mycal{R}_{\rho_{\alpha}}$:] the $\lan{T}$ bundle defined by
$\mycal{R}_{s_{\alpha}} := \alpha^{\vee}(\sR_{\alpha}) \in
\sBun_{\widetilde{C},\lan{T}}$. Here $\alpha^{\vee}$ is the
coroot corresponding to $\alpha$, viewed as a cocharacter
$\alpha^{\vee} : \mathbb{C}^{\times} \to \lan{T}$.
\end{itemize}

In \cite[Lemma~I.5.4, Proposition~I.5.5]{ron-dennis} it is shown that
the assignment 
$\alpha \mapsto \mycal{R}_{\rho_{\alpha}}$ extends uniquely to a
map $\mycal{R} : W \to \sBun_{\widetilde{C},\lan{T}}$
which is multiplicative in the sense that for all $w, w' \in W$ we
have a canonical isomorphism
\[
\varpi(w,w') : \mycal{R}_{w\cdot w'} \stackrel{\cong}{\longrightarrow}
\dia_{w'}\left(\mycal{R}_{w}\right)\otimes\mycal{R}_{w'}. 
\]
Here $\otimes$ denotes the natural tensor product of $\lan{T}$-bundles
on $\widetilde{C}$, $\dia_{w}(\bullet) :=
w^{*}\left((\bullet)\times_{w} \lan{T}\right)$ denotes the action of
$w$ on the groupoid of $\lan{T}$-bundles which is the combination of
the  pullback $w^{*}(\bullet)$ of
$\lan{T}$-bundles via the automorphism $w : \widetilde{C} \to
\widetilde{C}$, and the pushout $(\bullet)\times_{w} \lan{T}$ of
$\lan{T}$-bundles via the automorphism $w \in W \subset \op{Aut}(\lan{T})$. 

\begin{rem} \label{rem-cocycle} \ {\bfseries (i)} Consider the
  commutative group stack $\sBun_{\widetilde{C},\lan{T}}$ of all
  $\lan{T}$-bundles on $\widetilde{C}$. The group $W$ acts on
  $\sBun_{\widetilde{C},\lan{T}}$, i.e. a $w \in W$ acts on
  $\sBun_{\widetilde{C},\lan{T}}$ via $\ell \to \dia_{w}(\ell)$ for
  any $\lan{T}$-bundle  $\ell$. In these terms $\mycal{R}$ is just a
  1-cocycle of $W$ with values in the stacky $W$-module
  $\sBun_{\widetilde{C},\lan{T}}$. 

\

\noindent
{\bfseries (ii)} \ In our case $\widetilde{C}$ is not just an abstract
cameral cover but is a cameral cover for $K_{C}$-valued Higgs
bundles. For such $\widetilde{C}$ we can describe the line bundles
$\sR_{\alpha}$ and the cocycle $\mycal{R}$ more explicitly.

In this case we have a natural map $\widetilde{C} \to
\op{tot}(\lan{\mathfrak{t}}\otimes K_{C})$ and we can use this map to
describe the divisor $D^{\alpha}$. Indeed, since $\alpha$ is a linear
functional on $\mathfrak{t}$, we can view it as a map
\[
\mathfrak{t}\otimes K_{C} \to K_{C}
\]
of vector bundles on $C$ or equivalently, if we write $p :
\op{tot}(\mathfrak{t}\otimes K_{C}) \to C$ for the natural projection,
we can view $\alpha$ as a section 
\[
\alpha \in H^{0}(\op{tot}(\mathfrak{t}\otimes K_{C}), p^{*}K_{C}).
\]
Restricting this section to $\widetilde{C}$ we get a section 
$\alpha_{\widetilde{C}} \in H^{0}(\widetilde{C},\pi^{*}K_{C})$ whose
divisor is precisely $D^{\alpha}$. In particular
$\alpha_{\widetilde{C}}$ gives an isomorphism 
\[
\alpha_{\widetilde{C}} : \sR_{\alpha} \widetilde{\to}
 \pi^{*}K_{C}.
\]
Similarly
we get identifications 
\[
\mycal{R}_{\rho_{\alpha}} =
\alpha^{\vee}(\pi^{*}K_{C})
\]
and more generally $\mycal{R}_{w} = (\rho^{\vee} -
w\rho^{\vee})(\pi^{*}K_{C})$, where
$\rho^{\vee}$ denotes the half-sum of all positive coroots. 
\end{rem}

\

\punkt {\bf Twisted automorphisms.} \
Next suppose we are given a $\lan{T}$-bundle $\mycal{L}$. We can use the
cocycle $\mycal{R}$ to define the group of $\mycal{R}$-twisted
automorphisms of $\mycal{L}$ covering the action of $W$ on
$\widetilde{C}$. Namely we set 
\[
\op{Aut}_{\mycal{R}}(\mycal{L}) := \left\{ (w,i) \left| \;
\begin{minipage}[c]{1.7in} $w \in W$ \\
$i : \, w^{*}(\mycal{L}^{w}) \; \widetilde{\to}  \; \mycal{L}\otimes
  \mycal{R}_{w}^{-1}$  
\end{minipage}\right. \right\}
\]
with the group law given by ordinary composition in $W$ and the
isomorphisms \linebreak 
$\varpi(w,w') : \mycal{R}_{w\cdot w'} \to \mycal{R}_{w}\otimes
w^{*}\left(\left(\mycal{R}_{w'}\right)^{w}\right)$ which are part of
the cocycle data. 

Note also that for every $\mycal{L}$ we get a 
complex of groups 
\begin{equation} \label{eq:complexforl}
\lan{T} \to \op{Aut}_{\mycal{R}}(\mycal{L}) \to W. 
\end{equation}
As explained in \cite[Section~I.6]{ron-dennis} this complex is exact for a
connected compact $\widetilde{C}$ but this need not be the case in
general.

\

\punkt {\bf A $\mathbb{C}^{\times}$-torsor from a principal $SL_{2}$.}
\ Suppose now that $\alpha$ is a simple root for 
$\lan{G}$, and let $M^{\alpha} = P^{\alpha}/\op{Rad}_{u} P^{\alpha}$
be the maximal reductive quotient of
the minimal  parabolic $P^{\alpha}$ given by $\alpha$. Let $G_{\alpha} \to
[M^{\alpha},M^{\alpha}]$ be the universal cover of the semi-simple
part of $M^{\alpha}$. Note that $G_{\alpha} \cong SL_{2}(\mathbb{C})$,
and that 
\[
[M^{\alpha},M^{\alpha}] = \begin{cases} G_{\alpha}, & \quad
  \text{if } \alpha^{\vee} \text{ is primitive}, \\ 
G_{\alpha}/(\pm 1), & \quad \text{if } \alpha^{\vee} \text{ is not
  primitive}.\end{cases} 
\]
Also, if we take $T_{\alpha} \subset G_{\alpha}$ to be the maximal
torus which is the preimage of the torus \linebreak $(\text{image of
  $\lan{T}$ in 
  $M^{\alpha}$}) 
\cap [M^{\alpha},M^{\alpha}]$,
then the coroot map $\alpha^{\vee} : \mathbb{C}^{\times} \to \lan{T}$
lifts uniquely to an isomorphism $\mathbb{C}^{\times}
\to T_{\alpha}$ which by abuse of notation we will also denote by
$\alpha^{\vee}$.  

Now consider the normalizer subgroup  $N(T_{\alpha}) \subset G_{\alpha}$ 
of $T_{\alpha}$ in $G_{\alpha}$ 
and set 
\[
L_{\alpha} := N(T_{\alpha}) - T_{\alpha}. 
\]
By construction $L_{\alpha}$ is naturally a
$T_{\alpha}$-torsor and hence a $\mathbb{C}^{\times}$-torsor via the
identification $\alpha^{\vee} : \mathbb{C}^{\times} \to T_{\alpha}$. 

\

\

\medskip

\punkt {\bf Spectral data.} \ 
With all this at hand we can now define the notion of spectral data 
which is adapted to a given cameral cover. Let $\pi :
\widetilde{C} \to C$ be an abstract cameral cover for $\lan{G}$, 
$\sR := \{\sR_{\alpha}\}_{\alpha}$ denote the
collection of line bundles given by the divisors $D^{\alpha} \subset
\widetilde{C}$, and let  $\mycal{R} \in
Z^{1}(W,\sBun_{\widetilde{C},\lan{T}})$ be the corresponding cocycle.

\begin{defi} \label{defi-spectral} An {\bfseries $\lan{G}$  spectral datum
 of type $(\sR,\mycal{R})$} is a triple $(\mycal{L},\bi,\bb)$, where
\begin{description}
\item[{\bfseries\em (bundle)}] $\mycal{L}$ is a principal $\lan{T}$-bundle on
  $\widetilde{C}$;
\item[{\bfseries\em (twist)}] $\bi : N(\lan{T}) \to
  \op{Aut}_{\mycal{R}}(\mycal{L})$ is a homomorphism from the normalizer of
  $\lan{T}$ in $\lan{G}$ to the group of $\mycal{R}$-twisted
  endomorphisms of $\mycal{L}$. The homomorphim $\bi$ should fit in a
  commutative diagram
\[
\xymatrix@R-1pc{
0 \ar[r] & \lan{T} \ar[r] \ar@{=}[d] & N(\lan{T}) \ar[r] \ar[d]^-{\bi}
& W \ar[r] \ar@{=}[d] & 0 \\
& \lan{T} \ar[r] & \op{Aut}_{\mycal{R}}(\mycal{L}) \ar[r] & W,
}
\]
where the second row is the complex of groups \eqref{eq:complexforl}.
\item[{\bfseries\em (framing)}] $\bb = \{ \bb_{\alpha} \}_{
\begin{subarray}{l}\alpha \in \Delta(\lan{G}) \\ \alpha \text{ - simple}
\end{subarray}}$ is a collection of isomorphisms 
\[
\bb : \alpha(\mycal{L})_{|D^{\alpha}} \widetilde{\to}
\sR_{\alpha|D^{\alpha}}\otimes L_{\alpha}^{-1}
\]
of principal
  $\mathbb{C}^{\times}$-bundles on $D^{\alpha}$. Here as above
  $L_{\alpha}$ is viewed as a $\mathbb{C}^{\times}$-torsor (on a point) via
  the isomorphism $\alpha^{\vee} : \mathbb{C}^{\times} \to
  T_{\alpha}$, and $\otimes$ denotes the group law on
  $\mathbb{C}^{\times}$-torsors.  
\end{description}
The triple $(\mycal{L},\bi,\bb)$ should satisfy the following compatibility
condition. Let $n \in L_{\alpha}$, let $\overline{n} \in
N(\lan{T})$ denote the image of $n$ under the map
$\xymatrix@1{
L_{\alpha} & \hspace{-0.5in} \subset & \hspace{-0.7in}
N_{\alpha} \ar[r] & N(\lan{T})}$, 
and let $\bi(\overline{n}) :
\rho_{\alpha}^{*}\left(\mycal{L}^{\rho_{\alpha}}\right) \widetilde{\to}
\mycal{L}\otimes \mycal{R}_{\rho_{\alpha}}^{-1}$ be the corresponding
isomorphism of $\lan{T}$-bundles. When restricted to $D^{\alpha}$ this
gives an isomorphism $\bi(\overline{n})_{|D^{\alpha}} : \, 
\mycal{L}^{\rho_{\alpha}}_{|D^{\alpha}} \; \widetilde{\to} \;
\mycal{L}_{|D^{\alpha}}\otimes
\mycal{R}_{\rho_{\alpha}|D^{\alpha}}^{-1}$,  
which can be rewritten as an isomorphism $\bj(\overline{n},D^{\alpha})$:
\[
\xymatrix@R-1pc@C+1pc{
\left(\mycal{L}^{\rho_{\alpha}}\otimes \mycal{L}^{-1}\right)_{|D^{\alpha}}
\ar[r]^-{\cong} \ar@{=}[d] & \mycal{R}_{\rho_{\alpha}|D^{\alpha}}^{-1}
\ar@{=}[d] \\
-\alpha^{\vee}\left(\alpha(\mycal{L})_{|D^{\alpha}}\right)
\ar[r]_-{\bj(\overline{n},D^{\alpha})} &
\alpha^{\vee}\left(\sR_{\alpha|D^{\alpha}}\right).
}
\]
Then the compatibility condition on $(\mycal{L},\bi,\bb)$ is 
\begin{itemize}
\item[{\bfseries\em [C]}] For every simple root $\alpha$ of $\lan{G}$
  and every $n \in L_{\alpha}$ we have 
\[
\bj(\overline{n},D^{\alpha}) = - \alpha^{\vee}(\bb(n)).
\]
\end{itemize}
\end{defi}

\

\noindent
The main fact we will use is \cite[Theorem~6.4]{ron-dennis} which
idenitifies the Hitchin fiber $\lan{\gHiggs}_{\widetilde{C}}$ with the
  moduli stack of $\lan{G}$ spectral data of type $(\sR,\mycal{R})$.

\subsection{Abelianized Hecke functors} \label{ss:abelianized.heckes}

We want to define the abelianized Hecke functors
\begin{equation} \label{eq:abel.hecke}
\lan{\trans}^{\mu,\tilde{x}} := (-)\otimes
\mathcal{S}^{\mu,\tilde{x}} : \lan{\gHiggs}_{\widetilde{C}} \to
\lan{\gHiggs}_{\widetilde{C}}
\end{equation}
for $\mu \in \chr(G)$ and $\tilde{x} \in \widetilde{C}$. They are
translations which sends a given Higgs bundle, or equivalently a
spectral datum $(\mycal{L},\bi,\bb) \in
\lan{\gHiggs}_{\widetilde{C}}$, to its tensor product with a certain
$\lan{\mathcal{T}}$-torsor
$\mathcal{S}^{\mu,\tilde{x}}$. (\cite[Theorem~4.4]{ron-dennis} states
that $\lan{\gHiggs}_{\widetilde{C}}$ is a gerbe over
$\gTors_{\widetilde{C},\lan{\mathcal{T}}}$.  So a Higgs bundle can be
tensored with a $\mathcal{T}$-torsor.) To define
$\mathcal{S}^{\mu,\tilde{x}}$  we will use the norm map:
\[
\op{\sf{Nm}} : \gTors_{\widetilde{C},\lan{T}} \to \gTors_{C,\lan{\mathcal{T}}}
\]
associating a natural $\lan{\mathcal{T}}$-torsor on $C$ with each
$\lan{T}$-torsor on $\widetilde{C}$. The norm map $\op{\sf{Nm}}$ is
defined as follows. Averaging over $w \in W$ sends a $\lan{T}$-bundle
$F$ to a $W$-equivariant $\lan{T}$-bundle $S = \otimes_{w \in W}
\dia_{w}(F)$ on $\widetilde{C}$. Consider the subsheaf of sections
in $S$ whose value at each point in $\widetilde{C}$ is fixed under the
stabilizer of that point. This subsheaf automatically descends to a
sheaf of $\lan{\overline{\mathcal{T}}}$-modules on $C$. This sheaf is
representable by a $\lan{\overline{\mathcal{T}}}$-torsor
$\overline{\mathcal{S}}$. The associated line bundle $\alpha(S)$ comes
with a preferred frame along $D_{\alpha}$. The subsheaf of sections in
$S$ which over $D_{\alpha}$ map to the the preferred frame in
$\alpha(S)$ and are fixed under the reflection corresponding to
$\alpha$ will descend to a sheaf of $\lan{\mathcal{T}}$-modules on
$C$. This sheaf is representable by $\lan{\mathcal{T}}$-torsor which
we define to be $\op{\sf{Nm}}(F)$.

Now starting with the $\lan{T}$ torsor
$\mu(\mathcal{O}_{\widetilde{C}}(\tilde{x}))$ on $\widetilde{C}$, we
can apply $\op{\sf{Nm}}$ to obtain a $\lan{\mathcal{T}}$-torsor
$\mycal{S}^{\mu,\tilde{x}}$.  This provides the definition of our abelianized
Hecke functors $\lan{\trans}^{\mu,\tilde{x}}$. (We thank the referee
for suggesting improvements to an earlier version of this discussion.)

We can also rewrite the  operation \eqref{eq:abel.hecke}
as an operation on spectral data. For instance, given
$(\mycal{L},\bi,\bb) \in \lan{\gHiggs}_{\widetilde{C}}$ the action of
$\lan{\trans}^{\mu,\tilde{x}}$ on $(\mycal{L},\bi,\bb)$ results in a
new spectral datum $(\mycal{L}\otimes
S^{\mu,\tilde{x}},\bi,\bb\otimes \sss^{\mu,\tilde{x}})$ where
$S^{\mu,\tilde{x}}$ 
is a $W$-equivariant $\lan{T}$ bundle on $\widetilde{C}$ given by:
\[
S^{\mu,\tilde{x}} := \bigotimes_{w \in W}\displaylimits (w\mu)\left(
\mathcal{O}_{\widetilde{C}}(w\tilde{x})\right). 
\]
More precisely, note that the notion of spectral data defined above
depends on the collection of line bundles $\sR$ on the $D^{\alpha}$'s
and on the cocycle $\mycal{R}$. But the same definition makes sense
for any pair $(\sS,\mycal{S})$, where $\sS$ is a collection of line
bundles on the $D^{\alpha}$'s and $\mycal{S}$ is a
$\sBun_{\widetilde{C},\lan{T}}$-valued cocycle for $W$ compatible with
$\sS$, i.e. is equipped with $\lan{T}$-torsor isomorphisms
$\mycal{S}_{\rho_{\alpha}|D^{\alpha}} =
\alpha^{\vee}\left(\sS_{\alpha}\right)$ for all $\alpha$. In
particular we can take the pair $(\mathcal{O}^{\times},0)$ to consist
of the trivial line bundles and the zero cocycle. A spectral datum for
this pair will be a triple $(S,\bi,\sss)$, where $S$ is
just a $W$-equivariant bundle and $\bi$ is the
idenity.

 From the definition it is clear that the collection of spectral data of type
$(\sR,\mycal{R})$ is a torsor over the collection of spectral data of type
$(\mathcal{O}^{\times},0)$.  In particular
 $\mathcal{S}^{\mu,\tilde{x}}$ will correspond to spectral data
 $(S^{\mu,\tilde{x}},\op{id},\sss^{\mu,\tilde{x}})$ and so 
$\lan{\trans}^{\mu,\tilde{x}}$
will be given as
\[
\lan{\trans}^{\mu,\tilde{x}}(\mycal{L},\bi,\bb) :=
(\mycal{L},\bi,\bb)\otimes (S^{\mu,\tilde{x}},\op{id},\sss^{\mu,\tilde{x}}) =
(\mycal{L}\otimes S^{\mu,\tilde{x}},\bi,\bb\otimes \sss^{\mu,\tilde{x}}).
\]
The associated abelianized Hecke functor is the integral transform
corresponding to the structure sheaf of the graph of
$\lan{\trans}^{\mu,\tilde{x}}$, i.e. we have 
\[
\lanab{\mathbb{H}}^{\mu,x}_{\widetilde{C}} :=
\left(\lan{\trans}^{\mu,\tilde{x}}\right)^{*} :
D^{b}_{\op{qcoh}}(\lan{\gHiggs}_{\widetilde{C}},\mathcal{O}) \to
D^{b}_{\op{qcoh}}(\lan{\gHiggs}_{\widetilde{C}},\mathcal{O}).
\]
These abelianized Hecke functors appear in the
proof of Theorem~\ref{thm:gerbes} and in the statements in
Section~\ref{ss:hecke}.
Their Langlands dual versions
${\trans}^{\lambda,\tilde{x}}$,
${}_{\op{ab}}{\mathbb{H}}^{\lambda,\tilde{x}}_{\widetilde{C}}$ are used in 
section~\ref{ss:global}.


\begin{thebibliography}{DDD{\etalchar{+}}05}

\bibitem[AKS05]{aks}
P.~Argyres, A.~Kapustin, and N.~Seiberg.
\newblock {\em On {S}-duality for non-simply-laced gauge groups}.
\newblock  J. High Energy Phys.  2006,  no. 6, 043.

\bibitem[Ari02]{arinkin}
D.~Arinkin.
\newblock {\em Fourier transform for quantized completely integrable systems}.
\newblock PhD thesis, Harvard University, 2002.

\bibitem[BB09]{oren}
O.~Ben-Bassat.
\newblock {\em Twisting derived equivalences}.
\newblock  Trans. Amer. Math. Soc.  {\bfseries 361}  (2009),  no. 10,
5469--5504.  


\bibitem[BeBra07]{bb}
R.~Bezrukavnikov and A.~Braverman.
\newblock {\em Geometric {L}anglands correspondence for {D}-modules in prime
  characteristic: the {G}{L}(n) case}.
\newblock  Pure Appl. Math. Q.  {\bfseries 3}  (2007),  no. 1, part 3,
153--179.  

\bibitem[BBD82]{bbd}
A.~Beilinson, J.~Bernstein, and P.~Deligne.
\newblock Faisceaux pervers.
\newblock In {\em Analysis and topology on singular spaces, I (Luminy, 1981)},
  volume 100 of {\em Ast\'erisque}, pages 5--171. Soc. Math. France, Paris,
  1982.

\bibitem[BD03]{beilinson-drinfeld-langlands}
A.~Beilinson and V.~Drinfeld.
\newblock Quantization of {H}itchin's integrable system and {H}ecke
  eigensheaves.
\newblock {B}ook, in preparation, 2003.

\bibitem[BJSV95]{bjsv}
M.~Bershadsky, A.~Johansen, V.~Sadov, and C.~Vafa.
\newblock {\em Topological reduction of {$4$}{D} {SYM} to {$2$}{D}
  {$\sigma$}-models}.
\newblock {\em Nuclear Phys. B}, 448(1-2):166--186, 1995.

\bibitem[Br87]{breen-cube} L.~Breen 
\newblock The cube structure on the
  determinant bundle. 
\newblock Theta functions---Bowdoin 1987, Part 1
  (Brunswick, ME, 1987), 663--673, Proc. Sympos. Pure Math., 49, Part
  1, Amer. Math. Soc., Providence, RI, 1989.

\bibitem[Br90]{breen-bitorsors} L.~Breen 
\newblock {\em Bitorseurs et cohomologie non ab\'{e}lienne}. 
\newblock The Grothendieck Festschrift, Vol. I,  401--476,
Progr. Math., 86, Birkhäuser Boston, Boston, MA, 1990.  

\bibitem[DDD{\etalchar{+}}06]{d3hp}
D.-E.~Diaconescu, R.~Donagi, R.~Dijkgraaf, C.~Hofman, and T.~Pantev.
\newblock {\em Geometric transitions and integrable systems}.
\newblock  Nuclear Phys. B  {\bfseries 752}  (2006),  no. 3, 329--390. 

\bibitem[DDP07]{ddp}
 D.E.~Diaconescu, R.~Donagi, and T.~Pantev.
\newblock {\em Intermediate Jacobians and $ADE$ Hitchin systems}.  
\newblock Math. Res. Lett.  {\bfseries 14}  (2007),  no. 5, 745--756.

\bibitem[DG02]{ron-dennis}
R.~Donagi and D.~Gaitsgory.
\newblock {\em The gerbe of {H}iggs bundles}.
\newblock Transform. Groups, {\bfseries 7} (2):109--153, 2002.

\bibitem[Dim04]{dimca-sheaves}
A.~Dimca.
\newblock Sheaves in topology.
\newblock Universitext. Springer-Verlag, Berlin, 2004.

\bibitem[Don93]{donagi}
R.~Donagi.
\newblock {\em Decomposition of spectral covers}.
\newblock {\em Ast\'erisque}, 218:145--175, 1993.

\bibitem[Don95]{donagi-msri}
R.~Donagi.
\newblock {\em Spectral covers}.
\newblock In {\em Current topics in complex algebraic geometry (Berkeley, CA,
  1992/93)}, volume~28 of {\em Math. Sci. Res. Inst. Publ.}, pages 65--86.
  Cambridge Univ. Press, Cambridge, 1995.

\bibitem[DP08]{dp}
R.~Donagi and T.~Pantev.
\newblock Torus fibrations, gerbes, and duality.
\newblock with an appendix by D.Arinkin, Memoirs of the AMS, volume
          {\bf 193} (2008), 90 pp.


\bibitem[Dr83]{drinfeld} V.~Drinfeld.  
\newblock {\em Two-dimensional
$l$-adic representations of the fundamental group of a curve over a
finite field and automorphic forms on ${\rm GL}(2)$}.
\newblock Amer. J. Math.  {\bf 105}  (1983), no. 1, 85--114.

\bibitem[Fal93]{faltings}
G.~Faltings.
\newblock {\em Stable {$G$}-bundles and projective connections}.
\newblock J. Algebraic Geom., {\bfseries 2} (3):507--568, 1993.

\bibitem[FGV02]{fgv-glc} 
E.~Frenkel, D.~Gaitsgory and K.~Vilonen.
\newblock {\em On the geometric Langlands conjecture}.
\newblock J. Amer. Math.
Soc.  {\bf 15} (2002), no. 2, 367-417

\bibitem[G02]{dennis-glc} 
D.~Gaitsgory.
\newblock {\em Geometric Langlands correspondence for $GL_{n}$}, 
\newblock Proceedings of the International Congress of Mathematicians, Vol. II
(Beijing, 2002), 571-582.

\bibitem[GM03]{gelfand-manin}
S.~Gelfand and Y.~Manin.
\newblock Methods of homological algebra.
\newblock Springer Monographs in Mathematics. Springer-Verlag, Berlin, second
  edition, 2003.

\bibitem[God73]{godement}
R.~Godement.
\newblock Topologie alg\'ebrique et th\'eorie des faisceaux.
\newblock Hermann, Paris, 1973.
\newblock Troisi\`eme \'edition revue et corrig\'ee, Publications de l'Institut
  de Math\'ematique de l'Universit\'e de Strasbourg, XIII, Actualit\'es
  Scientifiques et Industrielles, No. 1252.

\bibitem[Gro57]{grothendieck-tohoku}
A.~Grothendieck
\newblock {\em Sur quelques points d'alg\'{e}bre homologique}.
\newblock  Tohoku
  Mathematical Journal. Vol. 119, 1957.

\bibitem[Hit87]{hitchin}
N.~Hitchin.
\newblock {\em Stable bundles and integrable systems}.
\newblock Duke Math. J., {\bfseries 54} (1):91--114, 1987.

\bibitem[Hit92]{Teich}
N.~Hitchin.
\newblock {\em Lie groups and {T}eichm\"uller space}.
\newblock Topology, {\bfseries 31} (3):449--473, 1992.

\bibitem[HT03]{hausel-thaddeus}
T.~Hausel and M.~Thaddeus.
\newblock {\em Mirror symmetry, {L}anglands duality, and the {H}itchin system}.
\newblock Invent. Math., {\bfseries 153} (1):197--229, 2003.

\bibitem[Ive86]{iversen}
B.~Iversen.
\newblock Cohomology of sheaves.
\newblock Universitext. Springer-Verlag, Berlin, 1986.

\bibitem[Kos63]{kostant}
B.~Kostant.
\newblock {\em Lie group representations on polynomial rings}.
\newblock Amer. J. Math., 85:327--404, 1963.

\bibitem[KW07]{kw}
A.~Kapustin and E.~Witten.
\newblock {\em Electric-magnetic duality and the geometric {L}anglands
  program}. 
\newblock  Commun. Number Theory Phys.  {\bfseries 1}  (2007),  no. 1, 1--236.

\bibitem[Lau96]{laumon} 
G.~Laumon. 
\newblock {\em Transformation de {F}ourier g\'{e}n\'{e}ralis\'{e}e}.
\newblock preprint, aXiv:alg-geom/9603004.

\bibitem[Loo97]{looijenga}
E.~Looijenga.
\newblock {\em Cohomology and intersection homology of algebraic varieties}.
\newblock In {\em Complex algebraic geometry (Park City, UT, 1993)}, volume~3
  of {\em IAS/Park City Math. Ser.}, pages 221--263. Amer. Math. Soc.,
  Providence, RI, 1997.

\bibitem[Mum09]{mumford-abelian} D.~Mumford.  
\newblock Abelian varieties.  
\newblock With appendices by C. P. Ramanujam and Yuri
Manin. Corrected reprint of the second (1974) edition. Tata Institute
of Fundamental Research Studies in Mathematics, 5. Published for the
Tata Institute of Fundamental Research, Bombay; by Hindustan Book
Agency, New Delhi, 2008. xii+263 pp.



\bibitem[Sai90]{saito-mhm} M.~Saito 
\newblock {\em Mixed Hodge modules}.
\newblock Publ. Res. Inst. Math. Sci. {\bf 26} (1990), no. 2, 221--333. 


\bibitem[Sco98]{scog}
R.~Scognamillo.
\newblock {\em An elementary approach to the abelianization of the
  {H}itchin system 
  for arbitrary reductive groups.}
\newblock Compositio Math., {\bfseries 110} (1):17--37, 1998.

\bibitem[SGA]{sga4}
M. {A}rtin, {A}. {G}rothendieck and {J.-L. Verdier}. 
\newblock {T}h\'{e}orie des topos et
  cohomologie \'{e}tale des sch\'{e}mas. 
\newblock Lecture Notes in Math. 269, 270 and 305, Springer-Verlag (1972 and
  1973).

\bibitem[Sim94]{simpson-moduli1}
C.~Simpson.
\newblock {\em Moduli of representations of the fundamental group of a smooth
  projective variety - {I}}.
\newblock Publications Math\'{e}matiques de l'{I.}{H.}{E.}{S.},
  {\bfseries 79}:47--129, 1994.

\bibitem[Sim95]{simpson-moduli2}
C.~Simpson.
\newblock {\em Moduli of representations of the fundamental group of a smooth
  projective variety - {II}}.
\newblock Publications Math\'{e}matiques de l'{I.}{H.}{E.}{S.},
	  {\bfseries 80}:5--79,
  1995.


\bibitem[Sim97]{simpson-hodge.filtration}
C.~Simpson.
\newblock {\em The Hodge filtration on  nonabelian cohomology}.
\newblock Algebraic geometry---Santa Cruz 1995, 217--281,
Proc. Sympos. Pure Math., 62, Part 2, 1997. 


\bibitem[Ver95]{verdier}
J.-L. Verdier.
\newblock Dualit\'e dans la cohomologie des espaces localement compacts.
\newblock In {\em S\'eminaire Bourbaki, Vol.\ 9}, pages Exp.\ No.\ 300,
  337--349. Soc. Math. France, Paris, 1995.

\end{thebibliography}

\newcommand{\etalchar}[1]{$^{#1}$}

\end{document}